\documentclass[12pt, reqno]{amsart}
\usepackage[T1]{fontenc}
\usepackage{textcomp} 
\renewcommand{\rmdefault}{ptm} 
\usepackage[margin=1 in]{geometry}

\pdfmapfile{+mtpro2.map}

\usepackage[subscriptcorrection,nofontinfo,mtpscr,mtpfrak,mtpbb]{mtpro2}

\usepackage{amsthm,mathtools, extarrows,textcmds}
\usepackage{array}

\usepackage{tikz-cd} 
\usepackage{xcolor, xspace} 

\usepackage[nodisplayskipstretch]{setspace}
\newcounter{numquote}

\newcommand{\defi}[1]{{\emph{\textsf{#1}}}}
\usepackage[scaled=0.87]{helvet}  

\usepackage[shortlabels]{enumitem}

\tikzset{%
  symbol/.style={
    draw=none,
    every to/.append style={
      edge node={node [sloped, allow upside down, auto=false]{$#1$}}
    },
  },
}

\usepackage{float}
\restylefloat{table}

\usepackage{aliascnt} 
\usepackage[colorlinks=true,linkcolor={red!70!black},citecolor={green!50!black},urlcolor={blue!70!black},citebordercolor={0 0 1},urlbordercolor={0 0 1},linkbordercolor={0 0 1},pagebackref]{hyperref}
\usepackage[normalem]{ulem} 
\usepackage[nameinlink]{cleveref}
\usepackage[alphabetic,msc-links]{amsrefs}

      \newtheorem{thm}{Theorem}[section]
      \newtheorem{lem}[thm]{Lemma}
      \newtheorem{prop}[thm]{Proposition}
      
      \newtheorem{cor}[thm]{Corollary}
      
      \newtheorem*{lem*}{Lemma}
\newtheorem*{thm*}{Theorem}
\newtheorem*{prop*}{Proposition}
      \theoremstyle{definition}
      \newtheorem{emp}[thm]{}
      
      \newtheorem{exmp}[thm]{Example}
      \newtheorem{rem}[thm]{Remark}

\theoremstyle{definition}
{
\newtheorem*{exmp*}{Example}
\newtheorem*{defn*}{Definition}
\newtheorem*{rem*}{Remark}
\newtheorem*{ans*}{Answer}
\newtheorem*{quest*}{Question}
}

\numberwithin{equation}{section}
\numberwithin{figure}{section}
\numberwithin{table}{section}
\makeatletter
    \let\c@equation\c@thm
    \let\c@figure\c@thm
    \let\c@table\c@thm
\makeatother

\tikzset{
    labl/.style={anchor=south, rotate=90, inner sep=.5mm}
}

\newcommand{\ol}[1]{\overline{#1}}

\usepackage{accents}
\usepackage{trimclip}

\newcommand{\wt}[1]{\mathchoice
{\accentset{\scalebox{1.1}[0.85]{\trimbox{0pt .3ex}{\textasciitilde}}}{#1}}
{\accentset{\scalebox{1}[1]{\trimbox{0pt .4ex}{\textasciitilde}}}{#1}}
{\accentset{\scalebox{.8}[.8]{\trimbox{0pt .5ex}{\textasciitilde}}}{#1}}
{\accentset{\scalebox{.7}[.6]{\trimbox{0pt .5ex}{\textasciitilde}}}{#1}}
}

\usepackage{adjustbox}
\newcommand{\inv}{^{-1}}
\newcommand{\norm}[1]{\left\lVert#1\right\rVert}
\newcommand{\Q}{\mathbf{Q}}

\newcommand{\Z}{\mathbf{Z}}
\newcommand{\C}{\mathbf{C}}
\newcommand{\N}{\mathbf{N}}
\newcommand{\F}{\mathbf{F}}
\newcommand{\p}{\mathbf{P}}
\newcommand{\X}{\mathbb{X}}
\newcommand{\Y}{\mathbb{Y}}

\newcommand{\mm}{\mathbb{M}}

\let\L\relax
\newcommand{\L}{\mathbb{L}}
\newcommand{\B}{\mathbb{B}}
\newcommand{\J}{\mathbb{J}}

\newcommand{\kk}{\mathbb{K}}

\newcommand{\calp}{\mathcal{P}}

\newcommand{\calz}{\mathcal{Z}}
\newcommand{\cald}{\mathcal{D}}

\newcommand{\calb}{\mathcal{B}}

\newcommand{\calc}{\mathcal{C}}

\newcommand{\calj}{\mathcal{J}}

\newcommand{\zhat}{\what{\mathbf{Z}}}
\newcommand{\A}{\mathbf{A}}

\newcommand{\scry}{\mathscr{Y}}
\newcommand{\mf}{\mathfrak}

\newcommand{\inj}{\lhook\joinrel\longrightarrow}
\newcommand{\surj}{\twoheadrightarrow}

\newcommand{\gal}{\textup{Gal}}

\newcommand{\Hom}{\textup{Hom}}

\newcommand{\iso}{\xrightarrow{\raisebox{-0.7ex}[0ex][0ex]{$\;\sim\;$}}}

\newcommand{\disc}{\textup{disc}}

\newcommand{\gmk}{\mbf{G}_{m,k}}

\newcommand*{\DashedArrow}[1][]{\mathbin{\tikz [baseline=-0.25ex,-latex, dashed,#1] \draw [#1] (0pt,0.5ex) -- (1.3em,0.5ex);}}

\newcommand{\scrx}{\mathscr{X}}

\newcommand{\scrj}{\mathscr{J}}

\newcommand{\scrv}{\mathscr{V}}

\newcommand{\kbar}{\ol{k}}
\renewcommand{\k}{k}
\newcommand{\g}{g}

\newcommand{\ceq}{\coloneqq}

\newcommand{\too}{\longrightarrow}
\newcommand{\mtoo}{\longmapsto}

\newcommand{\units}[1]{{#1}^{\times}}
\newcommand{\bb}[1]{\mathbb{#1}}
\newcommand{\st}{\,\vert\,}

\newcommand{\pgl}{\mbf{PGL}}

\newcommand{\cyc}[2]{\textup{${#1}_{#2}^{\textsf{cyc}}$}}
\newcommand{\dih}[2]{\textup{${#1}_{#2}^{\textsf{dih}}$}}

\let\Re\relax
\let\Im\relax
\let\bul\relax
\newcommand{\Re}[2]{\textup{${#1}_{#2}^{\textsf{re}}$}}
\newcommand{\Im}[2]{\textup{${#1}_{#2}^{\textsf{im}}$}}
\newcommand{\bul}[2]{{#1}_{#2}^{\bullet}}

\newcommand{\fcyc}[1]{\cyc{f}{#1}}
\newcommand{\fdih}[1]{\dih{f}{#1}}

\newcommand{\fRe}[1]{\Re{f}{#1}}
\newcommand{\fIm}[1]{\Im{f}{#1}}

\newcommand{\Xcyc}[1]{\cyc{\X}{#1}}
\newcommand{\Xdih}[1]{\dih{\X}{#1}}
\newcommand{\Xre}[1]{\Re{\X}{#1}}
\newcommand{\Xim}[1]{\Im{\X}{#1}}
\newcommand{\Xbul}[1]{\bul{\X}{#1}}

\newcommand{\Xwtcyc}[1]{\cyc{\wt{\X}}{#1}}
\newcommand{\Xwtdih}[1]{\dih{\wt{\X}}{#1}}
\newcommand{\Xwtre}[1]{\Re{\wt{\X}}{#1}}
\newcommand{\Xwtim}[1]{\Im{\wt{\X}}{#1}}
\newcommand{\Xwtbul}[1]{\bul{\wt{\X}}{#1}}

\newcommand{\Jcyc}[1]{\cyc{\J}{#1}}
\newcommand{\Jdih}[1]{\dih{\J}{#1}}
\newcommand{\Jre}[1]{\Re{\J}{#1}}
\newcommand{\Jim}[1]{\Im{\J}{#1}}
\newcommand{\Jbul}[1]{\bul{\J}{#1}}

\newcommand{\Acyc}[1]{\cyc{\bb{A}}{#1}}

\newcommand{\Aim}[1]{\Im{\bb{A}}{#1}}
\newcommand{\Abul}[1]{\bul{\bb{A}}{#1}}

\newcommand{\Bbul}[1]{\bul{\B}{#1}}

\newcommand{\Jwtdih}[1]{\dih{\wt{\J}}{#1}}
\newcommand{\Jwtre}[1]{\Re{\wt{\J}}{#1}}
\newcommand{\Jwtim}[1]{\Im{\wt{\J}}{#1}}
\newcommand{\Jwtbul}[1]{\bul{\wt{\J}}{#1}}

\newcommand{\kkcyc}[1]{\cyc{\kk}{#1}}

\newcommand{\kkre}[1]{\Re{\kk}{#1}}
\newcommand{\kkdih}[1]{\dih{\kk}{#1}}

\newcommand{\bbfbul}[1]{\bul{\bb{f}}{#1}}

\newcommand{\bbfcyc}[1]{\cyc{\bb{f}}{#1}}
\newcommand{\bbfRe}[1]{\Re{\bb{f}}{#1}}
\newcommand{\bbfIm}[1]{\Im{\bb{f}}{#1}}
\newcommand{\wtbbfbul}[1]{\bul{\wt{\bb{f}}}{#1}}

\newcommand{\wtbbfRe}[1]{\Re{\wt{\bb{f}}}{#1}}
\newcommand{\wtbbfIm}[1]{\Im{\wt{\bb{f}}}{#1}}

\newcommand{\gcyc}[1]{\cyc{g}{#1}}
\newcommand{\gdih}[1]{\dih{g}{#1}}
\newcommand{\gre}[1]{\Re{g}{#1}}
\newcommand{\gbul}[1]{\bul{g}{#1}}
\newcommand{\gim}[1]{\Im{g}{#1}}

\newcommand{\gwtcyc}[1]{\cyc{\tilde{g}}{#1}}
\newcommand{\gwtdih}[1]{\dih{\tilde{g}}{#1}}
\newcommand{\gwtre}[1]{\Re{\tilde{g}}{#1}}
\newcommand{\gwtbul}[1]{\bul{\tilde{g}}{#1}}
\newcommand{\gwtim}[1]{\Im{\tilde{g}}{#1}}

\newcommand{\Qim}[1]{\Im{Q}{#1}}
\newcommand{\Qre}[1]{\Re{Q}{#1}}
\newcommand{\Qdih}[1]{\dih{Q}{#1}}
\newcommand{\Qcyc}[1]{\cyc{Q}{#1}}
\newcommand{\Qbul}[1]{\bul{Q}{#1}}

\newcommand{\wtQdih}[1]{\dih{\wt{Q}}{#1}}

\newcommand{\wtQbul}[1]{\bul{\wt{Q}}{#1}}

\newcommand{\Wbul}[1]{\bul{W}{#1}}

\newcommand{\wtWbul}[1]{\bul{\wt{W}}{#1}}

\newcommand{\phiim}[1]{\Im{\phi}{#1}}
\newcommand{\phire}[1]{\Re{\phi}{#1}}
\newcommand{\phidih}[1]{\dih{\phi}{#1}}
\newcommand{\phicyc}[1]{\cyc{\phi}{#1}}
\newcommand{\phibul}[1]{\bul{\phi}{#1}}

\newcommand{\wtphiim}[1]{\Im{\wt{\phi}}{#1}}
\newcommand{\wtphire}[1]{\Re{\wt{\phi}}{#1}}
\newcommand{\wtphidih}[1]{\dih{\wt{\phi}}{#1}}
\newcommand{\wtphicyc}[1]{\cyc{\wt{\phi}}{#1}}
\newcommand{\wtphibul}[1]{\bul{\wt{\phi}}{#1}}

\newcommand{\mbQim}[1]{\Im{\mathbold{Q}}{#1}}
\newcommand{\mbQre}[1]{\Re{\mathbold{Q}}{#1}}
\newcommand{\mbQdih}[1]{\dih{\mathbold{Q}}{#1}}
\newcommand{\mbQcyc}[1]{\cyc{\mathbold{Q}}{#1}}
\newcommand{\mbQbul}[1]{\bul{\mathbold{Q}}{#1}}

\newcommand{\wtmbQre}[1]{\Re{\wt{\mathbold{Q}}}{#1}}

\newcommand{\wtmbQbul}[1]{\bul{\wt{\mathbold{Q}}}{#1}}

\newcommand{\mbWim}[1]{\Im{\mathbold{W}}{#1}}
\newcommand{\mbWre}[1]{\Re{\mathbold{W}}{#1}}
\newcommand{\mbWdih}[1]{\dih{\mathbold{W}}{#1}}

\newcommand{\mbWbul}[1]{\bul{\mathbold{W}}{#1}}

\newcommand{\wtmbWbul}[1]{\bul{\wt{\mathbold{W}}}{#1}}

\newcommand{\mbv}{\mathbold{v}}
\newcommand{\mbu}{\mathbold{u}}
\newcommand{\mbm}{\mathbold{m}}
\newcommand{\mbh}{\mathbold{h}}
\newcommand{\mbq}{\mathbold{q}}

\newcommand{\mbP}{\mathbold{P}}
\newcommand{\mbV}{\mathbold{V}}

\newcommand{\xcyc}[1]{\cyc{x}{#1}}
\newcommand{\xdih}[1]{\dih{x}{#1}}
\newcommand{\xRe}[1]{\Re{x}{#1}}
\newcommand{\xbul}[1]{\bul{x}{#1}}

\newcommand{\ycyc}[1]{\cyc{y}{#1}}
\newcommand{\ydih}[1]{\dih{y}{#1}}
\newcommand{\yRe}[1]{\Re{y}{#1}}
\newcommand{\ybul}[1]{\bul{y}{#1}}
\newcommand{\yim}[1]{\Im{y}{#1}}

\newcommand{\ywtcyc}[1]{\cyc{\wt{y}}{#1}}
\newcommand{\ywtdih}[1]{\dih{\wt{y}}{#1}}
\newcommand{\ywtRe}[1]{\Re{\wt{y}}{#1}}
\newcommand{\ywtbul}[1]{\bul{\wt{y}}{#1}}
\newcommand{\ywtim}[1]{\Im{\wt{y}}{#1}}

\newcommand{\ucyc}[1]{\cyc{u}{#1}}
\newcommand{\udih}[1]{\dih{u}{#1}}
\newcommand{\ure}[1]{\Re{u}{#1}}

\newcommand{\uim}[1]{\Im{u}{#1}}

\newcommand{\mbucyc}[1]{\cyc{\mathbold{u}}{#1}}
\newcommand{\mbudih}[1]{\dih{\mathbold{u}}{#1}}
\newcommand{\mbure}[1]{\Re{\mathbold{u}}{#1}}
\newcommand{\mbuim}[1]{\Im{\mathbold{u}}{#1}}

\renewcommand{\arraystretch}{1.15}

\usepackage{float}
\restylefloat{table}

\tikzset{
labl1/.style={anchor=north, rotate=90, inner sep=1.2mm}
}

\DeclareMathOperator{\chr}{\textup{char}}
\DeclareMathOperator{\spec}{\textup{Spec}} 
\DeclareMathOperator{\pic}{\textup{Pic}}
\DeclareMathOperator{\Div}{\textup{Div}}

\DeclareMathOperator{\rk}{\textup{rk}}

\DeclareMathOperator{\aut}{\textup{Aut}}

\DeclareMathOperator{\End}{\textup{End}}

\DeclareMathOperator{\tr}{\textup{Tr}}
\DeclareMathOperator{\Isog}{\textsf{Isog}}
\DeclareMathOperator{\lcm}{\textup{lcm}}
\DeclareMathOperator{\spn}{\textup{span}}

\newcommand{\trm}{\tr_{-}}

\let\deg\relax
\DeclareMathOperator{\deg}{\textup{deg}}
\DeclareMathOperator{\tors}{\textup{tors}}

\let\dim\relax
\DeclareMathOperator{\dim}{\textup{dim}}

\usepackage[old]{old-arrows}
\usepackage{xpatch}
\usepackage{multirow}
\usepackage{hhline}

\makeatletter
\xpatchcmd{\varinjlim}{\rightarrowfill@}{\varrightarrowfill@}{}{}
\xpatchcmd{\varinjlim}{\textstyle}{\scriptstyle}{}{}
\xpatchcmd{\varprojlim}{\leftarrowfill@}{\varleftarrowfill@}{}{}
\xpatchcmd{\varprojlim}{\textstyle}{\scriptstyle}{}{}
\makeatother

\newcommand{\mattwo}[4]{\begin{bsmallmatrix}
        #1 & #2 \\ #3 & #4
    \end{bsmallmatrix}}

\usepackage{microtype}
\begin{document}

\title[Constructing curves of large rank using composite polynomials]{Constructing curves of large rank using composite polynomials}
\author[Arvind Suresh]{Arvind Suresh}

\begin{abstract}
Let $\k$ be a number field.  We refine a construction of Mestre--Shioda to construct (infinite) families of hyperelliptic curves $X/{\k}$ having a record number of rational points and record Mordell--Weil rank relative to the genus of $g$ of $X$. Over $k=\Q$, we obtain modest improvements on the current published records, and these improvements become more significant as $k$ gets larger. For example,  we obtain curves over the real cyclotomic field ${\k} \ceq \Q(\cos(\pi/(g+1)))$ having at least $16(g+1)$ $\k$-points and rank at least $8g$. The defining equations for the curves are closely related to the classical Chebyshev polynomials, and in special cases, we recover families studied (for example) by Mestre, Shioda, Brumer, and Tautz--Top--Verberkmoes. 
\end{abstract}

\date{\today}
\maketitle

\tableofcontents

\section{Introduction} 
Let $\k$ be a number field, let $X/{\k}$ be a nice (i.e. smooth, projective, geometrically integral) curve of genus $\g$, and let $J_X/{\k}$ denote the Jacobian of $X$.  The Mordell-Weil Theorem says that $J_X({\k})$ is a finitely generated abelian group. The \defi{rank} of $X/{\k}$ is the integer
\begin{equation*}
    \rk X/{\k} \ceq \dim_{\Q} (J_X({\k}) \otimes_{\Z} \Q).
\end{equation*}
Faltings' Theorem~\cite{faltings}*{Satz~7}, formerly known as the Mordell Conjecture, says that if ${\g} \geqslant 2$, then the set $X({\k})$ of ${\k}$-rational points  of $X$ is finite. Given these finiteness theorems, it is natural to wonder, for a fixed genus ${\g}$ and number field $\k$, how large the integers $\#X({\k})$ and $\rk X/{\k}$ can be. We define two constants, taking values in $\N \cup \{\infty \}$, as follows:
\begin{align}\label{eq:constants}
\begin{split}
    R({\k},{\g}) & \ceq \sup \{ R \st \rk X/{\k} \geqslant R \textup{ for infinitely many genus ${\g}$ curves } X/{\k}\},\\
    N({\k},{\g}) & \ceq \sup \{ N \vert\, \#X({\k}) \geqslant N \textup{ for infinitely many genus ${\g}$ curves } X/{\k}\}.
\end{split}
\end{align}
Above, and for the rest of the article, by ``infinitely many curves over $k$,'' we mean ``infinitely many nice curves that are pairwise non-isomorphic over an algebraic closure ${\kbar.}$''\\

For any number field $\k$ and positive ${\g}$, whether $R({\k},{\g}) < \infty$ is an open question. On the other hand, it is a theorem of Caporaso, Harris, and Mazur~\cite{CHM}*{Theorem~1.2} that if the Weak Lang Conjecture~\cite{CHM}*{Conjecture~A} is true (which is not currently known) then $N({\k},{\g})< \infty$ for any number field $\k$ and ${\g}\geqslant 2$ (i.e. $\#X(k) \leqslant N(k,g)$ for all but finitely many ($\ol{k}$-isomorphism classes of) genus $g$ curves $X/k$).

A genus $1$ curve $X/{\k}$ with a rational point is an elliptic curve, and $X \cong J_X$. The  question of whether $R(\Q,1)<\infty $ has been the subject of multiple folklore conjectures; the papers \cites{Park-wood,discursus} give heuristics that suggest $R(\Q,1) \leqslant 21$ (in other words, they predict that only finitely many elliptic curves over $\Q$ have rank exceeding $21$). For ${\g} \geqslant 2$, to the best of our knowledge, there exist no such heuristics, nor a folklore conjecture in either direction.\\

Producing an infinite set of genus ${\g}$ curves $X/{\k}$ satisfying $\rk X/{\k} \geqslant R$ and $\#X({\k}) \geqslant N$ would yield lower bounds $R({\k},{\g}) \geqslant R$ and $N({\k},{\g}) \geqslant N$. The standard way to do this is by the method of \defi{specialization}; one constructs over the function field ${\k}(\scrv)$ of a $\k$-variety $\scrv$ a curve $X/{\k}(\scrv)$ that has the desired properties (genus ${\g}$, many rational points, large rank), and then one takes specializations of $X$ that preserve these properties. 
The use of the specialization method to construct curves of high rank originated in the work of N\'eron~\cites{neronhilbert,neronhilbert2}, who proved a theorem~\cite{neronhilbert}*{Th\'eor\`eme~6} on the injectivity of the specialization map (cf. \ref{emp:specialization}) for Jacobians of curves $X/{\k}(\scrv)$, with $\scrv$ rational over $\k$, and applied it to construct infinitely many elliptic curves of rank at least $10$ (i.e. $R(\Q,1)\geqslant 10$), and for ${\g}\geqslant2$, infinitely many genus ${\g}$ curves of rank at least $3{\g}+6$ (see also~\cite{umezu1} for a correction of N\'eron's original claim of $3{\g}+7$).\\

Mestre~\cites{mestre11,mestre12} constructed a genus ${\g}$ curve $\X/\kk$, with $\kk$ a  rational function field over $\Q$, having at least $8{\g}+12$ rational points. His construction yields $N(\Q,{\g}) \geqslant 8{\g}+12$, and he used it in the case ${\g}=1$ to prove $R(\Q,1) \geqslant 12$. Shioda~\cite{shiodasymmetry} refined and generalized Mestre's construction and proved $R(\Q,{\g}) \geqslant 4{\g}+7$ and $N(\Q,{\g}) \geqslant 8{\g}+16$ for all ${\g} \geqslant 2$; we refer to this as the \defi{Mestre--Shioda construction} (see \Cref{sec:MS}).

In this article, we present a refinement of the Mestre--Shioda construction which allows us to improve on the known lower bounds for the constants $N({\k},{\g})$ and $R({\k},{\g})$, for various pairs $({\k},{\g})$. In fact, our main result, \Cref{thm:main} below, gives bounds for the following natural generalization of the constant $R(k,g)$- for a finite abelian group $B$, we define the constant
\begin{equation}\label{eq:Rkg-B}
    R({\k},{\g}; B) \ceq \sup \{ R \st J_X(K) \supset \Z^R \oplus B \textup{ for infinitely many genus ${\g}$ curves } X/{\k}\}.
\end{equation}
Note that $R(K,g;0) = R(k,g)$. Below, and in the rest of the article, $Z_n$ denotes the cyclic group of order $n$, $V_4$ denotes the Klein-$4$ group $Z_2 \times Z_2$, $\zeta_n$ as usual denotes the primitive $n$-th root of unity $e^{2\pi i /n}$, and $\lambda_n$ denotes the real number $\zeta_n + \zeta_n\inv = 2\cos(2\pi/n)$ (so that $\Q(\lambda_n)$ is the maximal real subfield of the $n$-th cyclotomic field $\Q(\zeta_n)$).

\begin{thm}[Proved in \ref{emp:proof-thm-main}]\label{thm:main}
Let $n,d,$ and ${\g}$ be positive integers such that $nd$ equals one of the values in the first column of \Cref{tab:thm:main} below, and let ${\k}$ be a number field that contains the number $\alpha$ in the second column. Then, we have the bounds
\begin{equation*}
    N({\k},{\g}) \geqslant N \qquad \textup{ and } \qquad R({\k},{\g};B)\geqslant R,
\end{equation*}
where $B$ denotes the group in the third column and $N$ (resp. $R$) denotes the integer in the fourth column (resp. last three columns, depending on the condition satisfied by $d$ in the last row).
\begin{table}[H]
\centering
\renewcommand{\arraystretch}{1.4}
\setlength{\tabcolsep}{5pt}
\adjustbox{max width=\textwidth}{
\caption{Bounds for $k$ containing $\lambda_n,\lambda_{2n},$ or $\zeta_n$}\label{tab:thm:main}
\begin{tabular}{|c|c||c|c|l|l|l|} 
\hline 
$nd$    & $\alpha$          & $B$       & $N$       & \multicolumn{3}{c|}{$R$} \\
\hhline{|==||=====|}
$2g+2$  & $\lambda_n$       & $0$       & $8g+4n+8$ & $8g$      & $6g$  & $4g+2n+3$ \\     
\hhline{|--||-----|}
$2g+1$  & $\lambda_n$       & $0$       & $8g+4n+5$ & $8g$      &       & $4g+2n+2$ \\
\hhline{|--||-----|}
$2g$    & $\lambda_{2n}$    & $Z_2$     & $8g+4n+2$ & $8g$      & $6g$  & $4g+2n$ \\
\hhline{|--||-----|}
$2g-1$  & $\zeta_n$         & $V_4$     & $8g+4n-1$ & $8g-4$    &       & $4g+2n-2$ \\
\hhline{|--||-----|}
$g+1$   & $\zeta_n$         & $0$       & $8g+8n+8$ & $8g-4$    & $6g$  & $4g+4n+3$ \\
\hline
 \multicolumn{4}{c|}{\phantom{-}}       & $d=1$     & $d=2$     & {$d\geqslant 3$} \\
\cline{5-7}
\end{tabular}}
\end{table}
\end{thm}

\begin{emp}
\textbf{The case $k=\Q$.} For a given number field $k$, there are finitely many integers $n$ for which ${\k}$ contains $\lambda_{n},\zeta_n,$ or $\lambda_{2n}$. For each such $n$, by allowing $d$ (and therefore also $g$) to vary in \Cref{thm:main}, we obtain bounds for all sufficiently large genera $g$ comprising certain congruence classes. 

For example, $\Q$ contains $\lambda_n$ for $n=2$, so putting $2d = 2g+2$ in the first row of \Cref{tab:thm:main} we recover Shioda's bounds $N(\Q,g) \geqslant 8g+16$ and $R(\Q,g) \geqslant 4g+7$ (for $g\geqslant 2$) from \cite{shiodasymmetry}. This is not a coincidence; our method is a generalization of \cite{shiodasymmetry} (see \Cref{exmp:shioda-x2}).  

Observe now that $\Q$ also contains $\lambda_n$ if $n=3,4,$ or $6$. By setting $$6d = 2g+2,\; 4d= 2g+2,\;3d=2g+1,\;3d=2g, \;\textup{ and } \;d=2g-1$$ in \Cref{tab:thm:main}, respectively, we obtain the following bounds that improve on \cite{shiodasymmetry}.
\end{emp}


\begin{cor}\label{cor:bounds-Q}
For any positive integer $g$, we have the following bounds.
\begin{table}[H]
\centering
\renewcommand{\arraystretch}{1.3}
\setlength{\tabcolsep}{5pt}
\adjustbox{max width=\textwidth}{
\caption{Bounds for $k=\Q$}\label{tab:Q}
\begin{tabular}{|c|c|l|lc|} 
\hline 
Condition on $g$ & $B$ & $N(\Q,g)\geqslant $ & {$R(\Q,g;B) \geqslant $} &\\
\hline
$g \equiv 2 \pmod{3}$  & $0$ & $8g+32$ & $4g+15$ & if $g\geqslant 8$ \\     
\hline
$g \equiv 1 \pmod{2}$  & $0$ & $8g+24$ & $4g+11$ & if $g\geqslant 5$ \\     
\hline
$g \equiv 1 \pmod{3}$  & $0$ & $8g+17$ & $4g+8$ & if $g\geqslant 4$ \\
\hline
$g \equiv 0 \pmod{3}$  & $Z_2$ & $8g+14$ & $4g+6$ & if $g\geqslant 6$ \\
\hline
$g \equiv 1 \pmod{2}$  & $V_4$ & $8g+3$ & $4g$ & if $g\geqslant 1$ \\
\hline
\end{tabular}}
\end{table}
\end{cor}

\begin{emp}\label{emp:records-Q}
\textbf{New records for $\Q$.} We now summarize the state of the art for ${\k}=\Q$ and see where we obtain new records. 

We have $N(\Q,2)\geqslant 150$~\cite{elkiescover} and $N(\Q,4)\geqslant 126$~\cite{elkiesgenus4}*{Page~9}, both due to Elkies. For ${\g}=3$ and ${\g}\geqslant 5$, the current record is $N(\Q,{\g})\geqslant 8{\g}+16$, witnessed by Shioda's curves~\cite{shiodasymmetry}. Thus, \emph{\Cref{cor:bounds-Q} gives a new record for $N(\Q,{\g})$ for ${\g}=3$, and for all ${\g}\geqslant 5$ not divisible by $6$.} 

We have $R(\Q,1)\geqslant 19$ \cite{Elkies2006} and $R(\Q,3)\geqslant 26$~\cite{kulesz}*{Thm.~3.2.1}. The current records $R(\Q,2) \geqslant 32$ and $R(\Q,5)\geqslant 32$ have not explicitly appeared in print; the former follows by applying \cite{kulesz}*{Th\'eor\`eme~2.2.1} to the elliptic curve $E/\Q$ with $\rk E/\Q \geqslant 28$ found by Elkies~\cite{Elkies2006}, which  gives genus $2$ hyperelliptic curves $y^2 = f(x)$, with $\deg f = 6$. The record for ${\g}=5$ then follows by considering the set of genus $5$ curves given by $y^2 = f(x^2 + bx)$, with $b \in \Q$; we leave the details to the reader. For ${\g}=4$ and ${\g}\geqslant 6$, the existing record is Shioda's bound $R(\Q,{\g})\geqslant 4{\g}+7$~\cite{shiodasymmetry}. Thus, \emph{\Cref{cor:bounds-Q} gives a new record for $R(\Q,{\g})$ for ${\g}=4$, and for all ${\g}\geqslant 7$ not divisible by $6$.}

Consider next the constant $R(\Q, g;Z_2)$. The current record for $g = 1$
(elliptic curves) is $R(\Q, 1;Z_2) \geqslant 11$, due to Elkies~\cite{elkiesthreelectures}. For $g > 2$, the existing record is $R(\Q, g; Z_2) > 4g + 2$, due to Shioda~\cite{shiodasymmetry}*{Theorem~7}. Thus, \emph{\Cref{cor:bounds-Q} gives a new record for $R(\Q,g;Z_2)$ for
all g divisible by 3}.

Finally, we have $R(\Q,1;V_4) \geqslant 8$ due to Elkies~\cite{elkiesthreelectures}*{Page~13}. For $g\geqslant 2$, the bound $R(\Q,g;V_4) \geqslant 4g$ given by \Cref{cor:bounds-Q} appears to be entirely new; we do not know of any construction of large rank curves (even positive rank curves) of genus $g\geqslant 2$ over $\Q$ with Mordell--Weil group containing $V_4$. 
\end{emp}

\begin{rem}
The constants \eqref{eq:constants} and \eqref{eq:Rkg-B} may only increase as $k$ gets larger, so the bounds in \Cref{tab:Q} also hold for any number field $k$. Of course, \Cref{tab:thm:main} improves these bounds if $k$ contains cosines or roots of unity not already in $\Q$; we illustrate this in the case of the quadratic fields $\Q(\sqrt{3}) = \Q(\lambda_{12})$ and $\Q(\sqrt{-3}) = \Q(\zeta_6)$. We set $2g+2=12d$ in the former case and $g+1 = 6d$ in the latter (i.e. we look at the first and last row, respectively, of \Cref{tab:thm:main}). Then, \Cref{thm:main} yields in the case $d=1$ the bounds $ R(\Q(\sqrt{3}),5)\geqslant 40$ and $R(\Q(\sqrt{-3}),5)\geqslant 36$; compare these with the published record $R(k,5) \geqslant 32$ for any number field $k$ (cf. \ref{emp:records-Q}). For both of these quadratic fields, the case $d=2$ yields $R(k,11)\geqslant 66$, and the case $d\geqslant 3$ yields $R(k,g)\geqslant 4g+27$ for all $g\geqslant 17$ congruent to $5$ modulo $6$.
\end{rem}

When applying \Cref{thm:main}, instead of fixing $k$ and allowing $g$ to vary, we can equally well fix $g$ and allow $k$ to vary. From this perspective, \Cref{thm:main} almost \emph{doubles} the existing records for $R(k,g)$ once $k$ contains the appropriate cosines or roots of unity. Indeed, taking $d=1$ in \Cref{thm:main} yields the following. 

\begin{cor}\label{cor:8g}
For any  number field $k$ and positive integer $g$, we have 
\begin{equation*}
    R(k,g) \geqslant 
    \ccases{
        8g\; & \; \textup{ if $k$ contains $\lambda_{2g+2},\lambda_{2g+1},$ or $\lambda_{4g}$,}\\
        8g-4 \; & \;\textup{ if $k$ contains $\zeta_{g+1}$ or $\zeta_{2g-1}$}
    }
\end{equation*}
\end{cor}

\begin{exmp}
If $g$ is even then $\Q(\lambda_{2g+2}) = \Q(\lambda_{g+1})$, so in this case \Cref{cor:8g} gives $R(k,g) \geqslant 8g$ if $k$ contains $\lambda_{g+1}$. For example, it gives the bounds $R(k,6) \geqslant 48$ and $R(k,8) \geqslant 64$ if $k$ contains the real cubic fields $\Q(\lambda_7)$ and $\Q(\lambda_9)$, respectively. For comparison, note that the current published records for any number field $k$ are $R(k,6) \geqslant 31$ and $R(k,8) \geqslant 39$ (cf. \ref{emp:records-Q}). 
\end{exmp}

\begin{emp}
\textbf{Twist families.} 
The quest for infinitely many genus $g$ curves over $k$ with large rank admits a natural counterpart in which we ask for the curves to be pair-wise non-isomorphic over $k$ but isomorphic over $\ol{k}$ (recall that two such curves $X/k$ and $Y/k$ are then said to be \defi{twists} of one another). In this direction, the strongest result we know of is due to Mestre. He showed (by way of a refinement of his own original construction, see \Cref{exmp:shioda-x2}) that if $k$ contains $\zeta_{2g+2}$, with $g$ any positive integer, then there exist infinitely many (upto $k$-isomorphism) genus $g$ hyperelliptic curves
\begin{equation*}\label{eq:calx-ab}
   T_{a,b}/k: \, y^2 = ax^{2g+2} + b
\end{equation*}
having at least $16(g+1)$ $k$-rational points. 
Our refinement of the Mestre--Shioda construction is in fact a generalization of Mestre's refinement, and we are also lead to consider these curves $T_{a,b}/k$, which arise as specializations of the curve $\Xcyc{n}$ from \Cref{tab:data}. We are able to go further and show that these points generate rank $8g$ in the Jacobian. To state the full result, recall that an abelian variety $A/{\k}$ has \defi{complex multiplication} by some number field $E$ if $\dim A = [E:\Q]/2$ and there is a $\Q$-algebra embedding $E \inj \End^0_k(A) \ceq \End_k(A)\otimes_{\Z} \Q$. 
\end{emp}

\begin{thm}[Proved in \ref{emp:proof-thm-twist}]\label{thm:cm-twist}
Let $g$ be a positive integer, and let $n=2g+1$ or $2g+2$. Let $k$ be a number field which contains $\zeta_{n}$. Then, there exists a set $S \subset \A^2(k)$ which is not thin (see \ref{emp:thin sets}) such that for all $(a,b) \in S$, the curve $T_{a,b}/k:\, y^2 = ax^n + b$ is of genus $g$ and witnesses the following statements ($J_{a,b}/k$ denotes the Jacobian of $T_{a,b}/k$): 
\begin{enumerate}[(i)]
    \item There exists 
    a $k$-isogeny of abelian varieties
    \begin{equation*}
        \Theta : \prod_{m \mid n, \, m\geqslant 3} A_{m} \too J_{a,b},
    \end{equation*}
    where each $A_m/k$ is an abelian variety with complex multiplication by $\Q(\zeta_m)$. 
    \item For all $m$ dividing $n$, with $m\geqslant 3$, we have
    \begin{equation*}
        \rk A_m(k) \geqslant 4\phi(m) = 8 (\dim A_m),
    \end{equation*}
    and consequently, we have
    \begin{equation*}
        \rk \,T_{a,b}/k \geqslant 8g.
    \end{equation*}
\end{enumerate}
\end{thm}

\begin{emp}\label{emp:discussion}
The curves that witness \Cref{thm:main} are specializations of the curves defined in \Cref{tab:data}, which are equipped with many rational points over certain rational function fields (see \Cref{tab:data}). Let us here say a few words about these curves in the case $d=1$. 
\begin{enumerate}[(a),itemindent=1cm, leftmargin=0pt]
    \item Let $\fRe{n}(x) \in \Z[x]$ be the polynomial determined by the identity $\fRe{n}(x + x\inv) = x^n + x^{-n}$ (cf. \ref{emp:dih-mono}). That is, $\fRe{n}(x) = 2T_n(x/2)$, where $T_n(x)$ is the classical \defi{Chebyshev polynomial} (of the first kind) which is determined by the functional equation $T_n(\cos \theta) = \cos(n \theta)$. 
    The bounds in the first three rows of \Cref{tab:thm:main} are witnessed by specializations of the hyperelliptic curves $\Xre{n}$ and $\Xwtre{n}$ in \Cref{tab:data}. These specializations are given by equations of the form
    \begin{equation*}
        \calc_{a,b}/k: \, y^2 = a\fRe{n}(x) + b, \quad \textup{ and } \quad \wt{\calc}_{a,b}/k: (x+2)(a\fRe{n}(x) + b),
    \end{equation*}
    where $n = 2g+2$ or $2g+1$.
    The curves $\calc_{1,t}/k$ and $\wt{\calc}_{1,t}/k$, for $t\in k$, are denoted $\calc_t$ in the article \cite{tautz-real-mult} of Tautz, Top, and Verberkmoes. They show that if $n=p$ is an odd prime and $k$ contains $\lambda_p$,  then the Jacobian $\calj_t$ of $\calc_t$ has \emph{real multiplication} by $\Q(\lambda_p)$ (i.e. the endomorphism algebra $\End_k(\calj_t) \otimes_{\Z} \Q$ contains the field $\Q(\lambda_p)$). In our situation, one can show (but we do not, since it is not directly related to our main goal in this article) that for any $n\in \N$, the Jacobian of $\calc_{a,b}$ has endomorphism algebra containing the product of the fields $\Q(\lambda_m)$, where $m$ ranges over all divisors of $n$ greater than $2$. 
    \item Let $D_n$ denote the dihedral group of order $2n$. The bounds in the last row of \Cref{tab:thm:main} are witnessed by specializations of the hyperelliptic curves $\Xdih{n}$ (for $n=g+1$ even) and $\Xwtdih{n}$ (for $n=g+1$ odd) in \Cref{tab:data}. These specializations are given by equations of the form
    \begin{equation*}
        \cald_{a,b}/k: \, y^2 = a(x^{2g+2} + 1) + bx^{g+1}.
    \end{equation*}
    This curve is remarkable because, if $k$ contains $\zeta_{g+1}$, then $\aut_k(\cald_{a,b})$ contains a subgroup $H \cong Z_2 \times D_{g+1}$, with the $Z_2$ generated by $(x,y) \mapsto (x,-y)$, and the $D_n$ generated by the involution $(x,y) \mapsto (1/x,y/x^{g+1})$ together with the $n$-cycle $(x,y) \mapsto (\zeta_n x,y)$. Brumer~\cite{capratpts}*{Section~5} established the bound $N(\Q(\zeta_{g+1}),{\g})\geqslant 16(g+1)$ by constructing infinitely many $\cald_{a,b}/{\k}$ having $4$ full $H$-orbits ${\k}$-points. Our method similarly yields infinitely many such curves with four full $H$-orbits (note: putting $n=g+1$ in the last row yields $N = 16g+16$), but we are able to go one step further and determine that these points generate rank $8g-4$ in the Jacobian. We remark that, although Brumer's construction appears to be different from our (i.e. the Mestre--Shioda) construction, we suspect that a suitable change of variables might transform one into the other. 
    \item The curve $\Xre{n}$ (resp. $\Xwtre{n}$) mentioned in (a) is closely related to the curve $\Xdih{n}$ (resp. $\Xwtdih{n}$) from (b): the latter arise as a double cover of the former. These covers each give rise to so-called \emph{partner curves} (see \Cref{lem:hyp-partner} for the definition, and also \Cref{prop:partner-isog}), which are the curves $\Xim{n}$ and $\Xwtim{n}$ in \Cref{tab:data}. The specializations of $\Xim{n}$ witness the bounds in the second-last row of \Cref{tab:thm:main}; they are given by equations of the form
    \begin{equation*}
        \calb_{a,b}/k: \, y^2 = (x^2 - 4) (a\fRe{n}(x) + b),
    \end{equation*}
    where $n=2g-1$. The copy of $V_4$ in the Jacobian is generated by the divisor classes $[(2,0) - \infty]$ and $[(-2,0) - \infty]$. 
\end{enumerate}
\end{emp}

\begin{rem}
The analogues of Faltings' Theorem and the Mordell-Weil Theorem also hold for curves $X/K$, with $K$ the function field of a variety over a finite field $\F_q$. It is proved in \cite{ulmerpoints} (resp. \cite{ulmerranks}) that $N(\F_q(t),g) = \infty$ (resp. $R(\F_q(t),g) = \infty$); this automatically implies that the same is true if we replace $\F_q(t)$ with any function field over $\F_q$.
\end{rem}

\begin{emp}\label{emp:outline-article}
\textbf{Outline of the article.} In \Cref{sec:tools} we develop some tools for determining when a collection of points on a hyperelliptic curve generate a group of large rank in the Jacobian. In \Cref{sec:specialization} we review the specialization method alluded to in the introduction. In \Cref{sec:MS} we outline the Mestre--Shioda Construction, and in \Cref{sec:strategy} our basic strategy for refining the construction to get better bounds. 

The strategy requires us to understand when the splitting fields $\kk_f/\kk$ of certain composite polynomials $\bb{m}(f(x)) \in \kk[x]$ are rational function fields over $k$ (here, $\kk$ is a field showing up in Mestre's construction). We answer this question after giving an explicit construction of $\kk_f$ in \Cref{sec:splitting-field}. The answer implies, among other things, that the monodromy group of $f(x)$, defined in \ref{emp:Cf}, must be a finite subgroup of $\pgl_2(k)$.  It is this requirement that ultimately leads to the various conditions $\lambda_n \in k$ and $\zeta_n \in k$ present in \Cref{thm:main}; see \ref{emp:fin-subgp-pgl2}. The curves which witness the theorem are all hyperelliptic curves given by an equation of the form $y^2 = q(x)h(f(x))$, where $h$ is of degree $d$, $f$ is of degree $n$, and $q$ is of degree at most $2$ (hence the title of the article). In \Cref{thm:main}, the restrictions on the genera are essentially a consequence of the genus-degree formula for hyperelliptic curves.  

In \Cref{sec:mono}, we consider $\fRe{n}(x)$ (mentioned in \ref{emp:discussion}) and other closely related polynomials for which $\kk_f$ (constructed explicitly in \Cref{sec:kkn} following the method in \Cref{sec:splitting-field}) is a rational function field over $k$. Then, in \Cref{sec:curves} we define eight families of curves following the idea of \Cref{sec:strategy}. These curves come equipped with many rational points over the associated fields $\kk_f$. We define certain subgroups $\mbWbul{n}$ (in the Jacobian) associated to the points, and state in \Cref{thm:point-rank-bounds} a ``main'' theorem giving point and rank bounds for the curves. We utilize this main theorem to give the proof of \Cref{thm:main} in \ref{emp:proof-thm-main} (deferring to \Cref{sec:non-iso} the proof that certain curves over function fields are non-isotrivial). 

We outline in \ref{emp:rank-strategy} our overall strategy for bounding the ranks of the eight families of curves. The first step is a reduction step which allows us to reduce our consideration to three of the families; this is carried out in \Cref{sec:partner}. The next step consists first in fitting the curves together, by making use of some beautiful (and fairly intricate) polynomial identities, into inverse systems indexed by $n\in \N$ (see \Cref{sec:zhat}). These inverse systems give rise to various isogeny decompositions and relations between the respective Jacobians (which are possibly of independent interest and worth further exploration). In \Cref{sec:Wn-to-Vn}, we exploit these decompositions to reduce our task (of determining the ranks of the groups $\mbWbul{n}$) to determining the rank of certain distinguished subgroups $\bul{\mbV}{n} \subset \mbWbul{n}$. We accomplish this in \Cref{sec:rank-proof} by bringing to bear the tools from \Cref{sec:tools}. Finally, in \Cref{sec:cm} we use \Cref{thm:point-rank-bounds} to prove \Cref{thm:cm-twist}.
\end{emp}

\begin{emp}
\textbf{Acknowledgements.} This work is an extension of Chapter One of the author's PhD Thesis at the University of Georgia. The author is indebted to his PhD advisor Dino Lorenzini for his careful reading and thoughtful feedback on an earlier draft of this work. The author would like to thank Doug Ulmer for several useful discussions on the present draft, and in particular, for his help in clarifying some ideas behind \Cref{sec:zhat}. The author gratefully acknowledges the support of the Department of Mathematics at the University of Georgia as well as at the University of Arizona in providing a stimulating environment for research. 
\end{emp}






\section{Tools for bounding ranks of hyperelliptic curves}\label{sec:tools}


\begin{emp}\label{emp:conventions} \textbf{Notation and Conventions.}
Let $k$ be a field. 
We call a scheme $X/k$ a \defi{$k$-variety} if $X$ is geometrically integral, separated, and of finite type over $k$. In this case the function field of $X$ is denoted $k(X)$, and for a finite extension $\ell/k$, the function field of the base-change $X \times_k \ell$ is denoted $\ell(X)$. We call a morphism $\pi:X \too Y$ of $k$-varieties a \defi{cover} if it is finite and surjective. The \defi{degree} of the cover is as usual the degree of the induced function field extension $k(X)/k(Y)$. 

We use bold font as short-hand to denote tuples of indeterminates; for instance, if we say ``$\mbf{t}$ denotes the tuple of indeterminates $t_1,\dotsc,t_n$,'' then ${\k}[\mbf{t}]$ denotes ${\k}[t_1,\dotsc,t_n]$, $K(\mbf{t})$ the function field $K(t_1,\dotsc,t_n)$, and $\A^n_k{n}$ the affine space $\spec {\k}[\mbf{t}]$.

We call a field extension $\ell/k$ \defi{Galois} if it is algebraic, separable, and normal, and we denote its Galois group $\aut_k(\ell)$ by $\gal(\ell/k)$. We call $\ell/k$ a \defi{splitting field} of a polynomial $f(x) \in {\k}[x]$ if $f$ factors completely in $\ell[x]$ and $\ell$ is generated as $k$-algebra by the roots of $f$.
\end{emp} 

\begin{emp}
\textbf{Curves and Jacobians.} We call a curve $X/k$ a \defi{nice curve} if it is smooth, projective, and geometrically integral.  We denote by $\Div(X)$ the group of divisors on $X$, $\Div^0(X)$ the divisors of degree $0$, and ${\pic^0(X)}$ the group of degree $0$ divisors modulo principal divisors. We denote the Jacobian of $X/k$ by $J_X/k$;  recall that this is an abelian variety of dimension ${\g}$ (the genus of $X$) such that $J_X(\ell) \cong {\pic^0(X_{\ell})}$ for any field extension $\ell/k$ that verifies the condition $X(\ell) \neq \emptyset$. 
\end{emp}

\begin{emp}\label{emp:hyperelliptic}
\textbf{Hyperelliptic curves.} By a \defi{hyperelliptic curve} over a field $K$,  we mean a nice curve $X/K$ which admits an involution $\iota \in \aut_K(X)$ such that the quotient by $\iota$ is a double cover $\pi: X \too\p^1_K$; contrary to the convention in the literature, we do not require that the genus of $X$ be at least $2$. This introduces a slight ambiguity in that the double cover $\pi$ is unique (up to choosing a coordinate on $\p^1_K$) if and only if $g\geqslant 2$. In practice, we will typically work with hyperelliptic curves given by an explicit equation $X/K: y^2 = f(x)$, for some separable polynomial $f(x) \in {\k}[x]$ of degree $d\geqslant 1$ (this notation means $X/K$ is the smooth proper model of the affine plane curve defined by $y^2 = f(x)$). In this case, $\iota$ will always denote the involution $(x,y) \mapsto (x,-y)$, and $\pi:X \too\p^1_K$ the double cover $(x,y) \mapsto x$. We denote by $D_{\infty}$ the degree $2$ divisor $\pi^* \infty \subset X$. The genus ${\g}$ of $X$ is determined by the formula $d = 2g+1$ or $d=2g+2$, depending on the parity of $d$. 
\end{emp}

\begin{emp}\label{emp:tools-intro}
\textbf{Setup.} Let $L/K$ be an extension of perfect fields of characteristic different from $2$, let $g$ be a positive integer, and let $X/K$ be a genus $g$ hyperelliptic curve with Jacobian denoted $J/K$. We fix an involution $\iota \in \aut_K(X)$ whose quotient is a double cover $\pi:X \too\p^1_K$. This data will remain fixed until the end of this section.

In \Cref{lem:Vn-phin} below, we give sufficient conditions for a collection of $L$-points on $X$ to generate a group of large rank in the Jacobian. The main geometric tool we use is \Cref{lem:cast-sec-ineq}. This is an immediate consequence of the \emph{Castelnuovo--Severi Inequality} (see~\cite{griff-harris}*{Page~366} or~\cite{stichtenoth}*{Satz~1}), which says the following. 

\emph{Suppose $C,C_1,$ and $C_2$ are  nice curves over $K$ and $\varphi_1: C \too C_1$ and $\varphi_2:C\too C_2$ are covers which do not factor through a mutual cover $C \too \ol{C}$. Then, writing  $d_i$ for the degree of $\varphi_i$, the genera of the curves satisfy the inequality }
\begin{equation}\label{eq:cast}
     g(C) \,\leqslant\, g(C_1)(d_2 - 1) + g(C_2)(d_1 - 1) + (d_2 - 1)(d_1 - 1).
\end{equation}
\end{emp}

\begin{lem}\label{lem:cast-sec-ineq}
Let $\mf{D} \in \Div^0(X_{L})$ be a reduced divisor such that (i) $\mf{D}\neq \iota \mf{D}$ and (ii) the support of $\mf{D}$ consists of at most $2g$ geometric points of $X$. 
Then $[\mf{D}]$ is a non-zero point of $J(L)$. 
\end{lem}

\begin{proof}
Assume for a contradiction that $\mf{D}$ is the divisor of a rational function $\psi$ on $X_{L}$. Assumption (i) implies that the cover $\psi:X_{L} \too\p^1_{L}$ does not factor through the double cover $\pi:X_{L} \too\p^1_{L}$, and assumption (ii) implies that ${\g} \geqslant \deg \psi$. If we take $C$ to be $X$, $C_1$ and $C_2$ to both be $\p^1_{L}$, and the $\varphi_1$ and $\varphi_2$ to be $\pi$ and $\psi$, then the Castelnuovo-Severi Inequality \eqref{eq:cast} reads ${\g}\leqslant (\deg \psi)-1$, which is a contradiction.
\end{proof}


\begin{emp}\label{emp:tools-setup}
Fix an integer $n\geqslant 2$, and let $G \subset \aut_K(L)$ be a cyclic subgroup generated by an  automorphism $\sigma$ of order $n$.  Let $P$ be an $L$-rational point of $X$ which \defi{generates a $G$-torsor}; by this, we mean that the $G$-orbit $\{P,\sigma P,\dotsc, \sigma^{n-1}P\} \subset X(L)$ is of size $n$ (and is hence a $G$-torsor). Viewing $\Div X_{L}$ as a $\Z G$-module, this amounts to saying that the submodule $$W \ceq \spn_{\Z G} \{ P\} \subset \Div (X_{L})$$ is free of rank one over $\Z G$.  
Let $\Phi_n(x) \in \Z[x]$ denote the $n$-th cyclotomic polynomial; recall that $\Phi_n(x)$ is irreducible over $\Z$, since its roots are the primitive $n$-th roots of unity, and the degree of $\Phi_n$ is the Euler Totient $\phi(n)$. The \defi{$n$-th inverse cyclotomic polynomial} is defined by 
\begin{equation}\label{eq:Psi}
\Psi_{n}(x) \ceq \dfrac{x^n -1}{\Phi_n(x)} \in \Z[x]. 
\end{equation}
Any divisor $D \in W$ can be written in the form $D=q(\sigma)P$ for some polynomial $q(x) \in \Z[x]$, and the degree of such a divisor is the integer $q(1)$. Since $\Psi_n(1) = 0$, the divisor $\mf{P} \ceq \Psi_{n}(\sigma)P \in W$ is of degree $0$, and we may consider the $\Z G$-modules
\begin{equation}\label{eq:V-[V]}
\begin{alignedat}{2}
    V &\ceq \spn_{\Z G} \{ \mf{P} \} && \subset \Div^0 X_{L},\\
    [V] &\ceq \spn_{\Z G} \{ [\mf{P}]\} && \subset J({L}).
\end{alignedat}
\end{equation}
Identifying the group ring $\Z G$ with the commutative ring $\Z[\sigma]/(\sigma^n-1)$, we have the following.
\end{emp}

\begin{lem}\label{lem:V-[V]-basic}
The map $\Z G \too\End(V)$ factors through the natural quotient 
    \begin{equation*}
        \Z G \surj \Z[\zeta_n], \qquad {\sigma} \mtoo \zeta_n \ceq e^{2\pi i/n}
    \end{equation*}
to give an embedding $\Z[\zeta_n] \inj \End(V)$. In this way, $V$ is a free $\Z[\zeta_n]$-module of rank one, and $[V]$ verifies one of the following two statements, depending on the order of $[\mf{P}] \in J(L)$.
\begin{enumerate}[(a)]
    \item $[\mf{P}]$ is torsion; in this case $[V]_{\Q} \ceq [V] \otimes_{\Z} \Q = 0.$
    \item $[\mf{P}]$ is of infinite order; in this case $[V]_{\Q}$ is a free $\Q(\zeta_n)$-module of rank one, and hence, $\dim_{\Q} [V]_{\Q} =\phi(n)$. 
\end{enumerate}
\end{lem}

\begin{proof}
Using that $W$ is free of rank one over $\Z G$,  for any element $f({\sigma})\mf{P} \in V$, we have $f({\sigma})\mf{P} = f({\sigma})\Psi_{n}({\sigma})P = 0 \in V$ if and only if $x^n-1$ divides $f(x)\Psi_n(x)$ in $\Z[x]$, which holds if and only if $\Phi_n(x)$ divides $f(x)$ in $\Z[x]$. It follows that the kernel of the ring homomorphism $\Z G \too\End(V_n)$ is generated by $\Phi_n({\sigma})$, so  it factors through the quotient $\Z[\zeta_n] = \Z G/\Phi_n(\sigma)$ to give an embedding $\Z [\zeta_n] \inj \End(V)$. We have a natural surjection of $\Z G$-modules $[-]:V \surj [V]$ defined by $E \mapsto [E]$, and since $\Phi(\sigma)$ annihilates $V$, it also annihilates $[V]$. Thus, $[V]_{\Q}$ is a $\Q(\zeta_n)$-module generated by $[\mf{P}]$, which immediately implies that one of the statements, (a) or (b), holds. 
\end{proof}

In the next lemma, we give sufficient conditions for the point $[\mf{P}]\in J(L)$ to be of infinite order. The main idea is to produce a distinguished divisor $\mf{D} \in V$ satisfying the conditions of \Cref{lem:cast-sec-ineq} such that, if $[\mf{P}]$ were  torsion, then $[\mf{D}] = 0 \in J(L)$, contradicting \Cref{lem:cast-sec-ineq}. To make this work, we require that the torsion subgroup of $J(L)$ be bounded (in the sense of \ref{item:tors}), and that $g$ be sufficiently large relative to $n$. Define the integer
\begin{equation}\label{eq:norm-n}
    \norm{n} \ceq  \textup{ the number of distinct prime factors of } n.
\end{equation}

\begin{lem}\label{lem:Vn-phin}
Continuing with \ref{emp:tools-setup} above, assume that the following hold. 
\begin{enumerate}[label={\textup{(a\arabic*)}}]
    \item \label{item:P-neq-iotaP} $P$ and $\iota P$ generate $G$-torsors which are disjoint in $X(L)$.
    \item \label{item:tors} We have an equality of torsion subgroups $J(L)_{\tors} = J(K)_{\tors}.$
    \item \label{item:g-norm} The genus ${\g}$ of $X$ satisfies ${\g}\geqslant \max \{2, 2^{\norm{n}-1} \}$ (in particular, $g\geqslant 2$).
\end{enumerate}
Then $[V] \subset J(L)$ is a free abelian group of rank $\phi(n)$ (i.e. statement (b) of \Cref{lem:V-[V]-basic} holds). 
\end{lem}


\begin{proof}
Assume for a contradiction that $[\mf{P}]$ is a torsion point. We proceed in two separate cases. 

\begin{enumerate}[(i), leftmargin=0pt,itemindent=1cm]
    \item \emph{The case $\norm{n}=1$.} In this case, \ref{item:g-norm} says that $g \geqslant 2$, and hence, that $4 \leqslant 2g$. We have $n=p^e$ for some prime $p$ positive integer $e$, so $\mf{P}$ is the divisor $({\sigma}^{n/p}-1)P$. Define 
    \begin{equation*}
        \mf{D} \ceq ({\sigma}-1)\mf{P} = ({\sigma}^{\frac{n}{p}+1}P + P) -({\sigma}^{\frac{n}{p}}P + {\sigma}P) \in V.
    \end{equation*}
    Then \ref{item:P-neq-iotaP} implies that $\mf{D}$ is a reduced divisor with support consisting of $4$ points, and that $\mf{D} \neq \iota \mf{D}$, so \Cref{lem:cast-sec-ineq} implies $[\mf{D}]\neq 0 \in J({L})$. But \ref{item:tors} forces $[\mf{D}] = ({\sigma}-1)[\mf{P}] = 0 \in J({L})$, a contradiction. 
    \item \emph{The case $\norm{n}>1$.} In this case \ref{item:g-norm} says that $2^{\norm{n}} \leqslant 2g$. 
    Assumption \ref{item:tors} implies that $(1-\zeta_n)$ annihilates the $\Z[\zeta_n]$-module $[V]$. Since $n$ is composite, $(1-\zeta_n)$ is a unit in $\Z[\zeta_n]$, so $[V] = \{0\} \subset J(L)$. Define now the polynomial and divisor
    \begin{equation*}
        \wt{\Psi}_n(x) \ceq \prod_{p \textup{ prime}, \,p\mid n} (x^{\frac{n}{p}}-1) \in \Z[x], \qquad \mf{D} \ceq \wt{\Psi}_n(\sigma) P \in \Div^0 (X_L).
    \end{equation*}
    Since $\Psi_{n}(x)$ divides $\wt{\Psi}_n(x)$ in $\Z[x]$, the point $[\mf{D}] \in J(L)$ lies in $[V]$, so we have $[{\mf{D}}] = 0 \in J({L}).$ 
    Let $\calp$ denote the set of primes dividing $n$. For any subset $S\subset \calp$, define 
    \begin{equation*}
        {\sigma}_S \ceq \prod_{p \in S}{\sigma}^{n/p} \in G,
    \end{equation*}
    so that 
    \begin{equation*}
        \mf{D} = \sum_{S \subset \calp} (-1)^{\norm{n}-\#S} {\sigma}_SP \; \in \Div^0(X_{L}).
    \end{equation*}
    Each element ${\sigma}^{n/p} \in G$ is of order $p$, so each element ${\sigma}_S$ is of order $\prod_{p \in S} p$. This implies that the $2^{\norm{n}}$ elements $\{{\sigma}_S \in G \}_{S \subset \calp}$ are pair-wise distinct. Then, \ref{item:P-neq-iotaP} implies that $\mf{D} \neq \iota \mf{D}$, and that $\mf{D}$ is reduced with support consisting of the $2^{\norm{n}}$ points $\{{\sigma}_SP\}_{S \subset \calp}$. But \Cref{lem:cast-sec-ineq} implies that $[\mf{D}] \neq 0 \in J(L)$, a contradiction.\qedhere
\end{enumerate}
\end{proof}


\begin{emp}\label{emp:transc}
In \Cref{sec:curves} we define certain hyperelliptic curves over rational function fields, and in \Cref{tab:data} we exhibit tuples of rational points on these curves. \Cref{lem:Vn-phin} turns out to be a robust tool for determining the corresponding subgroups generated in the respective Mordell--Weil groups-- we use it to dispense with all the cases in which the curves in question are of genus at least $2$ (see the proof of \Cref{thm:dim-Vn}). This is a consequence of two pleasant features of the Mestre--Shioda construction (cf. \ref{emp:mestre} and \ref{emp:strategy}). The first is that these points enjoy a high degree of symmetry, namely, they generate torsors under cyclic groups (exactly as in \Cref{lem:mbV-direct-sum} above). The second is that they are \defi{transcendental points} of the curves (we call a point $P:\spec L \too X$ transcendental  if the image of $P$ is the generic point of $X$). We conclude this section with two lemmas which exploit this last feature-- \Cref{lem:alg-closed-tors} is used to obtain the assumption \ref{item:tors} in \Cref{lem:Vn-phin}, and \Cref{lem:dP-D-inf-order} is used in the proof of \Cref{thm:dim-Vn} to deal with the genus one cases. 
\end{emp}

\begin{lem}\label{lem:alg-closed-tors}
If $K$ is algebraically closed in $L$, then we have an equality of torsion subgroups $$J(L)_{\tors} = J(K)_{\tors}.$$
\end{lem}

\begin{proof}
Suppose $P:\spec L \too J$ is a torsion point of order $n$, and let $x \in J$ be the image of $P$. Then $x$ lies in $J[n]$, the kernel of multiplication by $n$. Since $J[n]$ is a finite $K$-scheme, the residue field $\kappa(x)$ is a finite extension of $K$. The morphism $P:\spec L \too\spec \kappa(x)$ gives a $K$-algebra inclusion  $\kappa(x)\inj L$. So, if $K$ is algebraically closed in $L$ then $\kappa(x) = K$ and we conclude that $P \in J(K)_{\tors}$ (more precisely, that $P$ is the base-change to $L$ of the point $x \in J(K)_{\tors}$).
\end{proof}

\begin{lem}\label{lem:dP-D-inf-order}
Let $P\in X(L)$ be a transcendental point, and let $D \in \Div X$ be an effective divisor of degree $d\geqslant 1$. Then, the point $[dP-D] \in J(L)$ is of infinite order.
\end{lem}

\begin{proof}
Put $\mf{D} \ceq dP - D$, and assume for a contradiction that $[\mf{D}]$ is a torsion point of $J(L)$. If $n$ denotes the order of this point, then by replacing $d$ (resp. $D$) with $nd$ (resp. $nD$), we may assume for a contradiction simply that $[\mf{D}] = 0 \in J(L)$. Recall that the {$d$-th symmetric power} of $X/K$ is a variety $X^{(d)}/K$ whose $F$-points (for any extension $F/K$) are naturally in bijection with the effective degree $d$ divisors on $X_F$. Consider the $K$-morphism $j:X \too J$ defined as the composition of the $K$-morphisms
    \begin{equation*}
        \begin{tikzcd}
            X \arrow[r, hook, "\Delta"] & X^d \arrow[r] & X^{(d)}  \arrow[r] & J,
        \end{tikzcd}
    \end{equation*}
    where the first arrow is the diagonal embedding, the second is the natural quotient by the symmetric group $S_{d}$, and the last arrow is defined on $F$-points (for any $L/K$) by $E \mtoo [E - D]$. Thus, for any point $Q \in X(L)$, we have $j(Q) = [dQ - D]$. Since the morphism $j:X \too J$ is non-constant, the image of the composition $j\comp P: \spec L \too J$ is not a closed point. But this is absurd since $j(P) = [\mf{D}] = 0 \in J(L)$.
\end{proof}

\section{The specialization method}\label{sec:specialization}

In this section, $k$ denotes an infinite, finitely generated field (\defi{i.f.g.} field). 

\begin{emp}\label{emp:families} 
\textbf{Non-isotrivial curves.} 
Let $\kk/{\k}$ be the function field of a $\k$-variety $\scrv$ of positive dimension, and let $\X/\kk$ be a nice curve of positive genus. We say that $\X/\kk$ is \defi{constant} if there is an $\kk$-isomorphism $\X \iso X_{\kk}$ for some curve $X/{\k}$. We say that $\X/\kk$ is \defi{isotrivial} if $\X_{\L}/\L$ is constant, for some finite extension $\L/\kk$. By replacing $\scrv$ with an open subset if necessary, we can spread out $\X/\kk$ to a smooth, projective family of curves $\scrx \longrightarrow \scrv$ with geometrically integral fibers (cf.~\cite{Qpoints}*{Theorem~3.2.1 (i)}). For ${\g}\geqslant 2$, this family induces a $\k$-morphism $\scrv \longrightarrow M_{g,k}$, where $M_{g,k}$ denotes the coarse moduli scheme of genus ${\g}$ curves over $\k$ (cf. ~\cite{mumfordGIT}*{Page 103, Proposition~5.4}). It is an exercise to see that $\X/\kk$ is non-isotrivial if and only this morphism is non-constant, or equivalently, if the geometric fibers $\{\scrx_P/\kbar\}_{P\in \scrv(\kbar)}$ comprise infinitely many $\kbar$-isomorphism classes. If $\X/\kk$ is an elliptic curve (i.e. of genus one and possessing a $\kk$-rational point), then it is isotrivial if and only if its $j$-invariant is in $\kk \setminus k$. 
\end{emp}

\begin{emp}\textbf{The specialization map.}\label{emp:specialization}
With $\scrx/\mathscr{V}$ and $\X/\kk$ as above, for a point $s \in \scrv$, the fiber $\scrx_s / \kappa(s)$ is called a \defi{specialization} of $\X$. The \defi{relative Jacobian} of $\scrx /\scrv$ is an abelian $\scrv$-scheme $\scrj/\scrv$~\cite{neronmodels}*{Proposition~9.4.4} such that, for any point $s \in \scrv$, the fiber $\scrj_s / \kappa(s)$ is isomorphic to the Jacobian of $\scrx_s / \kappa(s)$; in particular, the generic fiber of $\scrj/\scrv$ is isomorphic to the Jacobian $\J/\kk$ of $\X/\kk$. 

If $\scrv$ is regular, then the natural injective map $\scrj(\scrv) \inj \J(\kk)$  is also surjective (cf. \cite{movinglemma}*{Proposition~6.2}), so we have an isomorphism of groups $\J(\kk) \cong \scrj(\scrv)$. For any $P \in \scrv({\k})$, we also have a natural map $\scrj(\scrv) \longrightarrow \scrj_P({\k})$; pre-composing this with the isomorphism $\J(\kk) \cong \scrj(\scrv)$, we get the \defi{specialization map} $\sigma_P: \J(\kk) \longrightarrow \scrj_P({\k}),$ which is a group homomorphism. 
\end{emp}

\begin{emp} \textbf{Thin sets.}\label{emp:thin sets}
Recall from ~\cite{serremw}*{Page~121, Section~9.5} or~\cite{CTrankjump}*{Section 2} that a set $T\subset \scrv({\k})$ is \defi{thin} if there is a generically finite $\k$-morphism $\phi: \calz \too \scrv$, admitting no rational section $\scrv \DashedArrow \calz$, such that $\phi(\calz({\k})) \supset T$. We say that $\scrv/{\k}$ satisfies the \defi{Hilbert Property} if $\scrv({\k})$ is not thin; note that this property depends only on the birational isomorphism class of $\scrv/{\k}$. The usefulness of thin sets consists in the fact that if $\scrv/{\k}$ satisfies the Hilbert Property, then $\{\scrv({\k})\backslash T\st T \subset \scrv({\k}) \textup{ is thin}\}$ is a collection of Zariski-dense subsets of $\scrv({\k})$ that is closed under finite intersections. We say that $\k$ is \defi{Hilbertian} if $\p^n({\k})$ is not thin for every $n\geqslant 1$. It is known (see \cite{serremw}*{Section 9.5, Page~129}) that {ifg} fields (in particular, number fields) are Hilbertian.
\end{emp}

\begin{thm}[{\cite{CTrankjump}*{Theorem~2.4}, \cite{serremw}*{Section~11.1}}] \label{thm:specialization}
Let $\scrv$ be a smooth $\k$-variety with function field $\kk$, and let $\scrj/\scrv$ an abelian scheme with generic fiber $\bb{J}/\kk$. Then, the set 
\begin{equation*}
    T \ceq \{\, P \in \scrv({\k}) \st  \sigma_P : \bb{J}(\kk) \too \scrj_P({\k}) \textup{ is not injective } \}
\end{equation*}
is thin. In particular, if $\scrv/k$ satisfies the Hilbert Property, then $\scrv(k)\setminus T$ is Zariski-dense in $\scrv$. 
\end{thm}


Recall the constants defined in \eqref{eq:constants} and \eqref{eq:Rkg-B}.

\begin{cor}\label{cor:family-spec}
Let $k$ be an \textup{i.f.g.} field, let $R$ and $N$ be positive integers, and let $B$ be a finite abelian group. Let $\kk/{\k}$ be a rational function field, and $\X/\kk$ a non-isotrivial curve of genus $g\geqslant 1$. 
\begin{enumerate}[(a)]
    \item If $\#\X(\kk) \geqslant N$, then $N({\k},{\g}) \geqslant N$.\label{part:speccor a}
    \item If $\J(\kk)$ contains a subgroup isomorphic to $B\oplus \Z^R$, then $R({\k},{\g};B) \geqslant R$.\label{part:speccor c}
\end{enumerate}
\end{cor}

\begin{proof}
Choose a smooth $\k$-variety $\scrv$ with function field  $\kk$ which admits a smooth, projective family of curves $\scrx/\scrv$ with generic fiber $\X/\kk$ (as in \ref{emp:families}). Since $\scrv$ is a $k$-rational variety, it satisfies the Hilbert Property (see \cite{serremw}*{Section 9.5, Page~129}), so $\scrv({\k}) \setminus S$ is Zariski-dense in $\scrv$ for any thin set $S\subset \scrv({\k})$. The non-isotriviality of $\X/\kk$ then implies that: \emph{if $S \subset \scrv({\k})$ is thin, then the set $\{\scrx_P/{\k} \}_{P \in \scrv({\k})\setminus S}$ contains infinitely many nice genus ${\g}$ curves over $\k$.}
Now, part \ref{part:speccor a} follows because, after replacing $\scrv$ with a suitable open subset, any $N$ points of $\X(\kk)$ can be spread out to $N$ disjoint sections $\scrv \too \scrx$ (i.e. they specialize to $N$ distinct $\k$-points on any $\k$-fiber of $\scrx/\scrv$). Part \ref{part:speccor c} follows by applying \Cref{thm:specialization} to the relative Jacobian $\scrj/\scrv$ of $\scrx/\scrv$.
\end{proof}

\section{The Mestre--Shioda construction}\label{sec:MS}

\emph{In this section, and continuing until the end of \ref{emp:field-inclusions}, $k$ denotes a field of characteristic different from $2$}. The Mestre--Shioda construction  is predicated on the following lemma. For the rest of the article, we fix a positive integer $d$ and we set $$r=2d+2.$$ 

\begin{lem}\label{lem:sqroot}
Let $\mbm,\mbq,$ and $\mbh$ denote the tuples of indeterminates $m_0,\dotsc,m_{2d+1}, q_0,\dotsc,q_{d},$ and $h_0,\dotsc,h_{d},$ respectively. 
Define the polynomials
\begin{equation}\label{def:polys}
    \begin{alignedat}{3}
        \bb{m}(x)    & \ceq x^{2d+2} + m_{2d+1}x^{2d+1} && + \dotsb  + m_0    \; && \in  {\k}[\mbm][x],\\
        \bb{q}(x)    & \ceq \qquad \quad\;\; x^{d+1} + q_{d}x^{d} && + \dotsb + q_0     \;   && \in  {\k}[\mbq][x],\\
        \bb{h}(x) & \ceq \qquad \qquad \qquad \quad h_{d}x^{d} && + \dotsb + h_0   \; && \in  {\k}[\mbh][x].
    \end{alignedat}
\end{equation}
Then, the map ${\k}[\mbm] \too {\k}[\mbh,\mbq]$ defined by equating coefficients in the identity 
\begin{equation}\label{eq:sqroot-identity}
    \bb{m}(x) = \bb{q}(x)^2 - \bb{h}(x)
\end{equation}
is a $\k$-algebra isomorphism. 
\end{lem}

\begin{proof}
One can express each $q_i,h_i$ as a polynomial in the variables $m_0,\dotsc,m_{2d-1}$ by comparing coefficients (in descending order) in the identity $\bb{m}(x) = \bb{q}(x)^2 - \bb{h}(x)$ and sequentially solving for the variables $q_{d},\dotsc,q_0,h_{d},\dotsc,h_0$; we leave the details to the reader.
\end{proof}

\begin{emp}\label{emp:mestre}
Identifying ${\k}[\mbm]$ with ${\k}[\mbh,\mbq]$ via the isomorphism in \Cref{lem:sqroot}, we view $\bb{h}(x), \bb{q}(x)$ and $\bb{m}(x)$ as polynomials in ${\k}[\mbm][x]$ satisfying the identity \eqref{eq:sqroot-identity}. Let $\mbu$ denote the tuple of indeterminates $u_1,\dotsc,u_{r}$, let  $s_i(\mbu)$ denote the $i$-th elementary symmetric function in $\mbu$, and define the $\k$-algebra homomorphism
\begin{equation}\label{eq:km-to-ku}
    {\k}[\mbm] \inj {\k}[\mbu], \qquad m_i \mapsto (-1)^is_{r-i}(\mbu), \; \textup{ for } i=0,\dotsc,r-1.
\end{equation}
This homomorphism is injective and gives the factorization
\begin{equation}\label{eq:bbm-factor}
    \bb{m}(x) = (x - u_{1})\dotsb (x- u_{r}) \in {\k}[\mbu][x].
\end{equation}
Put $\mm \ceq k(\mbm)$ and $\kk \ceq k(\mbu)$, and define the $(d+1)$-dimensional rational function field
\begin{equation*}
    K \ceq {\k}(\mbh) = {\k}(h_0,\dotsc,h_d),
\end{equation*}
so that we have a chain of (transcendental) field extensions $k \inj K \inj \kk$. Now, 
Mestre's construction~\cites{mestre11,mestre12} extracts from the identities \eqref{eq:sqroot-identity} and \eqref{eq:bbm-factor} the following geometric consequence: the hyperelliptic curve 
\begin{equation*}
    \X/K:\, y^2 = \bb{h}(x)
\end{equation*}
 has the $\kk$-rational points
\begin{equation}\label{eq:Pi}
   Q_i \ceq (u_i,\bb{q}(u_i)), \quad \textup{ for } \quad i=1,\dotsc,r.   
\end{equation}
These points, together with their conjugates $\iota Q_i = (u_i,-\bb{q}(u_i)),$ contribute $2r = 4d+4$ $\kk$-points to $\X$. For $i=1,\dotsc,r$, define the divisor $\mf{Q}_i \ceq 2Q_i - D_{\infty} \in \Div^0(\X_{\kk})$ (cf. \ref{emp:hyperelliptic}). Then, we have
\begin{equation*}
    2\textup{div}(y - \bb{q}(x)) = \mf{Q}_1 + \dotsb + \mf{Q}_{r} \in \Div^0(\X_{\kk}),
\end{equation*}
giving the relation $[\mf{Q}_1] + \dotsb + [\mf{Q}_{2d+2}] = 0 \in \J(\kk)$, where $\J/\kk$ denotes the Jacobian of $\X$. 
Shioda showed in~\cite{shiodasymmetry}*{Theorems~4, 5} that if ${\g}$ is positive (i.e. if $d\geqslant 3$), then the points $[\mf{Q}_i]$ do not satisfy any other relation and hence generate a subgroup of rank $r-1 = 2d+1$ in $\J(\kk)$.

Now, if $d$ is even then the genus ${\g}$ of $\X$ satisfies $d=2g+2$, so we have the bounds
\begin{equation}\label{eq:X-kk-rk}
    \# \X (\kk)  \geqslant 2r= 8g+12 \qquad \textup{ and } \qquad \rk \X/\kk \geqslant 4g+5.
\end{equation}
Thus, if $k$ is an {i.f.g.} field (cf. \Cref{sec:specialization}) then we can apply  \Cref{cor:family-spec} to get the bounds $N({\k},{\g})\geqslant 8{\g}+12$ and $R({\k},{\g})\geqslant 4{\g}+5$. 
\end{emp}

\section{Strategy behind improvement}\label{sec:strategy}


\begin{emp}\label{emp:strategy}
Let $f(x) \in {\k}[x]$ be a polynomial of degree $n\geqslant 1$. Our basic idea, generalizing \Cref{exmp:shioda-x2} below, is to consider (variants of) the hyperelliptic curve
\begin{equation*}
    \X_f/K: \, y^2 = \bb{h}(f(x)).    
\end{equation*}
This curve fits into a commutative diagram:
\begin{equation}\label{diag:pullback-f}
    \begin{tikzcd}
    \X_f \arrow[r, "\phi_f"] \arrow[d, "\pi",swap] & \X \arrow[d, "\pi"] & & (x,y) \arrow[r, mapsto] \arrow[d, mapsto] & (f(x),y) \arrow[d, mapsto]\\
    \p^1_{K} \arrow[r, "f"] & \p^1_{K},  && x \arrow[r, mapsto] & f(x).
    \end{tikzcd}
\end{equation}
To simplify the exposition, suppose for the moment that $nd$ is even. Then, since $\bb{h}(f(x))$ is separable (cf. \Cref{lem:disc-g-f}), the genus ${\g}_f$ of $\X_f$ is determined by the formula
\begin{equation}\label{eq:gf}
    nd = 2g_f + 2,
\end{equation}
Let $\kk_f/\kk$ be a splitting field extension for $\bb{m}(f(x))$. One checks, by replacing $x$ with $f(x)$ in \eqref{eq:sqroot-identity}, that $\X_f$ has the $\kk_f$-points
\begin{equation}\label{eq:Pt}
    \{ Q_t \ceq (t,\bb{q}(f(t))) \st \bb{m}(f(t)) = 0 \}.
\end{equation}
Since $\bb{m}(f(x))$ has $rn$ distinct roots $t \in \kk_f$ (cf. \Cref{lem:disc-g-f}), the points \eqref{eq:Pt} and their hyperelliptic conjugates $\iota Q_t \ceq (t,-\bb{q}(f(t)))$ yield
\begin{equation}\label{eq:XfX-kkf}
    \#\X_f(\kk_f) \geqslant 2rn =  8{\g}_f + 4n + 8.
\end{equation}
Fixing any root $t$ of $\bb{m}(f(x))$ and letting $t'$ vary over the remaining roots, we have $rn - 1$ points $[Q_t - Q_{t'}] \in \J_f(\kk_f)$. If these are linearly independent, then we have
\begin{equation}\label{eq:rank-Xf-kkf}
    \rk \X_f/\kk_f \geqslant rn-1 =  4{\g}_f + 2n + 3.
\end{equation}
The key observation now is that, relative to ${\g}_f$, the number of $\kk_f$-points constructed grows linearly in $n$. Indeed, if $f(x) = x$ (so that $n=1$ and $d = 2g_f + 2$), then we have $\X_f = \X$ and $\kk_f = \kk$, and \eqref{eq:XfX-kkf} and \eqref{eq:rank-Xf-kkf} specialize to the bounds in  \eqref{eq:X-kk-rk}. At the other extreme, if $d=1$ and $n=2g_f + 2$, then \eqref{eq:XfX-kkf} and \eqref{eq:rank-Xf-kkf} would give lower bounds 
\begin{equation*}
    \# \X_f(\kk_f) \geqslant 16g_f + 16 \quad \textup{ and } \quad \rk \X_f/\kk_f \geqslant 8g_f + 7.\footnote{It turns out, in all the cases we consider, that the classes $[Q_t - Q_{t'}] \in \J_f(\kk)$ associated to the points \eqref{eq:Pt} never generate rank greater than $8g_f$. This occurs because the $K$-morphism $\phi_f:\X_f \to \X \cong \p^1_K$ introduces linear dependence relations of the form $[\phi_f^*E] = 0 \in \J_f(\kk)$, where $E$ is a degree zero divisor on $\p^1_{\kk}$.}
\end{equation*}
Supposing we are able to get a good rank bound, we would like to then specialize $\X_f/\kk_f$ (by applying \Cref{thm:specialization}) to obtain curves over $k$ with large rank. Thus, the usefulness of our strategy is tempered by the requirement that \emph{$\kk_f$ must be the function field of a $k$-variety which satisfies the Hilbert Property} (cf. \ref{emp:thin sets}). We examine when this holds in the next section.  
\end{emp}

\begin{exmp}\label{exmp:shioda-x2}
The observation that curves of the form $\X_f$ have more points relative to genus (than Mestre's ``original'' curve $\X$) was inspired by the following antecedents in the literature.
\begin{enumerate}[(a),leftmargin=0pt,itemindent=1cm]
    \item  In~\cite{shiodasymmetry}, Shioda worked with the case $f(x)=x^2$, so that $n=2$. The curve $\X_f/K: y^2 = \bb{h}(x^2)$ is non-constant of genus ${\g}_f = d-1$, and  Shioda proved that the inequality \eqref{eq:rank-Xf-kkf} holds if ${\g}_f\geqslant 2$. The splitting field $\kk_f$ of $\bb{m}(x^2) \in \mm[x]$ is a rational function field over $\k$ (see \Cref{exmp:quad-mono}), so the Specialization \Cref{thm:specialization} yields the bounds $R({\k},{\g}_f)\geqslant 4g_f+7$ and $N({\k},{\g}_f)\geqslant 8{\g}_f+16$ (which constitute the current published records for large enough ${\g}$).
    \item Mestre (see \cite{elkiesgenus4}*{Page 4}) considered the case $d=1$ and $f(x) = x^{2g+2}$. He observed that $\X_f/K: y^2 = h_1x^{2g+2} + h_0$ is a genus ${\g}$ curve having at least $16(g+1)$ rational points over $\kk_f$, which is again a rational function field over $\k$ (in fact, $\kk_f$ is the field $\kkcyc{n}$ defined in \ref{emp:fields-Ln}).  
\end{enumerate}
\end{exmp}


\section{Constructing the splitting field \texorpdfstring{$\kk_f$}{} of \texorpdfstring{$\bb{m}(f(x))$}{} over \texorpdfstring{$\kk$}{}}\label{sec:splitting-field}

We fix for this section a polynomial $f(x) \in k[x]$ of degree $n\geqslant 1$. In \ref{emp:const-spl-fld} below, we construct an abstract splitting field extension $\kk_f$ for the composition $\bb{m}(f(x)) \in \kk[x]$ (under a simplifying assumption on $f$). 
In \Cref{prop:spl-flt-rat} we give necessary and sufficient conditions (in terms of a certain curve ${C}_f/{k}$ associated to $f$) for $\kk_f$ to be a rational function field over $\k$. \Cref{lem:alg-closed} furnishes a key ingredient in determining the ranks of our curves. 


\begin{lem}\label{lem:disc-g-f}
Let $f(x) \in {k}[x]$ be a polynomial of degree $n\geqslant 1$ whose derivative $f'(x)$ is not identically zero. Let $\mbf{g}$ denote the tuple ${\g}_0,\dotsc,{\g}_{r-1}$, with $r\geqslant 1$, and define
    \begin{equation*}
        \bb{g}(x) \ceq x^r + g_{r-1}x^{r-1}+\dotsb+g_0 \in {k}[\mbf{g}][x].
    \end{equation*}
Then, the composition $\bb{g}(f(x)) \in k(\mbf{g})[x]$ is irreducible and separable. Moreover, the set of irreducible factors of $\disc (\bb{g}\circ f)$ in $\ol{k}[\mbf{g}]$ (up to scaling by $\units{\ol{K}}$) is
    \begin{equation}\label{eq:factors}
        \{\disc \, \bb{g}\} \cup \{\bb{g}(f(\alpha)) \st f'(\alpha) = 0 \}.
    \end{equation}
\end{lem}

\begin{proof}
We assume without loss of generality that $k$ is algebraically closed. Let $k(\mbf{g}) \inj k(\mbv) = K(v_1,\dotsc,v_r)$ be the injection which gives the factorization $\bb{g}(x) = (x-v_1)\dotsb (x-v_r) \in k(\mbv)[x]$, so that we have the factorization $\bb{g}(f(x)) = (f(x) - v_1)\dotsb (f(x) - v_r) \in k(\mbv)[x]$. Each $f(x) - v_i$ is irreducible over $k(v_i)$, and hence, also over $k(\mbv)$ (since $k(\mbv)/k(v_i)$ is purely transcendental). By the Fundamental Theorem of Symmetric Functions, $k(\mbv)/k(\mbf{g})$ is Galois, and its Galois group is the symmetric group $S_r$. This group acts transitively on the irreducible factors $f(x) - v_i \in k(\mbv)[x]$, which implies that their product $\bb{g}(f(x))$ is irreducible in $k(\mbf{g})[x]$.

We verify now that \eqref{eq:factors} indeed enumerates all the irreducible factors of $\disc(\bb{g}\comp f)$ (which also implies the separability of $\bb{g}\comp f$). Write $f'(x) = (x-\alpha_1)\dotsb (x-\alpha_{n-1}) \in {\k}[x]$. It follows from~\cite{odoni}*{Lemma~8.1 (i)} that $\disc \, \bb{g}$ is irreducible in ${\k}[\mbf{g}]$, and that
    \begin{equation*}
    \disc ( \bb{g}\comp f) = (\disc \,\bb{g})^n \prod_{i=1}^r \disc (f(x) - v_i) \in {\k}[\mbv].
    \end{equation*} 
    The $r$-term product above lies in ${\k}[\mbf{g}]$, indeed:
    \begin{equation*}
    \prod_{i=1}^r \disc (f(x) - v_i)  = \prod_{i=1}^r\prod_{j=1}^{n-1} (f(\alpha_j) - v_i) = \prod_{j=1}^{n-1} \bb{g}(f(\alpha_j)) \in {\k}[\mbf{g}].
    \end{equation*}
    Each $\bb{g}(f(\alpha_j)) \in {\k}[\mbf{g}]$, being linear in ${\g}_0,\dotsc,{\g}_{r-1}$, is irreducible in ${\k}[\mbf{g}]$,  so we are done.
\end{proof}

\begin{emp}\label{emp:Cf}
\textbf{The curve ${C}_f/k$ and group $G_f$.} 
Let $K_f/k(u)$ denote a splitting field for $f(x) - u\in k(u)[x]$, so that we have a factorization $f(x) - u = (x -t_1) \dotsb (x - t_n) \in K_f[x]$ in which the roots are pairwise distinct (by \Cref{lem:disc-g-f}). The \defi{monodromy group of $f$} is the Galois group 
\begin{equation}
    G_f \ceq \gal(f(x)-u/ k(u)) = \gal(K_f/k(u)).
\end{equation} 
Writing $R_f$ for the integral closure of $k[u]$ in $K_f$, the $G_f$-action on $K_f$ restricts to $R_f$, and $k[u]$ is the ring of invariants $(R_f)^{G_f}$. The \defi{Galois cover} associated to $f$ is the finite $k$-morphism 
\begin{equation}\label{eq:phif}
    \phi_f:C_f \ceq \spec R_f \too \A^1_{k}(u)
\end{equation}
corresponding to the inclusion $k[u] \subset R_f$. It is the quotient by the (right) action of $G_f$ on $C_f/k$. 
Note that $C_f/k$ is \emph{geometrically} integral if and only if $k$ is algebraically closed in $K_f$. 
\end{emp}

\begin{emp}\label{emp:const-spl-fld}
\textbf{Constructing $\kk_f/\kk$.}
Identify the $k$-algebra $k[u]^{\otimes r}$ with the polynomial algebra $k[\mbu]=k[u_1,\dotsc,u_{r}]$ (from \ref{emp:mestre}) by writing $u_i$ for the $i$-th simple tensor $1\otimes \dotsb \otimes u \otimes \dotsb \otimes 1 \in k[u]^{\otimes r}$. By taking the $r$-fold tensor product (over $k$) of the inclusion $k[u]\inj R_f$, we get an inclusion $k[\mbu] \inj R_f^{\otimes r}$. For each $i=1,\dotsc,r$ and $j=1,\dotsc,n$, let $t_{i,j} \in R_f^{\otimes r}$ denote the $i$-th simple tensor $1 \otimes \dotsb \otimes t_j \otimes \dotsb \otimes 1$ (where the $t_j$ is in the $i$-the place), so that $t_{i,1},\dotsc,t_{i,n}$ are the roots of $f(x) - u_i \in k[\mbu][x]$. By construction, we have  the factorization
\begin{equation}\label{eq:bbm-f-factor}
    \bb{m}(f(x)) = \prod_{i=1}^{r} (f(x) - u_i) = \prod_{i=1}^{r} \prod_{j=1}^n (x - t_{i,j}) \in R_f^{\otimes r}[x]. 
\end{equation}
The action of each monodromy group 
\begin{equation*}
    G_{i,f} \ceq \gal(f(x) - u_i/k(u_i)), \qquad i=1,\dotsc,r
\end{equation*}
extends naturally to an action on $R_f^{\otimes r}$ in which $G_{i,f}$ fixes $t_{e,j}$ for every $e \neq i$. These together give an action of $G_f^r = \prod_{i=1}^r G_{i,f}$ on $R_f^{\otimes r}$ with ring of invariants $k[\mbu]$. The cover of $k$-schemes $\phi_f^r: C_f^r \too \A^{r}_{k}(\mbu)$  is the quotient by the associated (right) action of $G_f^r$ on $C_f^r$. If $k$ is algebraically closed in $K_f$ then $C_f^r$ is a $k$-variety with function field $K_f^{\otimes r}$. So, the extension of function fields
\begin{equation}\label{eq:kk-inj-kCfr}
    \kk = k(\mbu) \inj k(C_f^r) = K_f^{\otimes r}
\end{equation}
induced by $\phi_f^r$ is Galois and we have a natural group isomorphism $\gal(k(C_f^r)/\kk) \cong G_f^r.$ This group acts faithfully on the roots $t_{i,j}$ of $\bb{m}(f(x))$ in $k(C_f^r)$, which implies the following. 
\end{emp}

\begin{lem}\label{lem:g-f-spl-fld}
If $K_f/k$ is regular, then $k(C_f^r)$ is a splitting field for $\bb{m}(f(x)) \in \kk[x]$, and we have $$\gal(\kk_f/\kk) = \prod_{i=1}^r G_{i,f} \cong G_f^r.$$
\end{lem}

\begin{exmp}\label{exmp:quad-mono}
If $f(x)$ is quadratic then the inclusion ${\k}(u) \subset k(t)$ defined by $u \ceq f(t)$ makes $k(t)$ into a splitting field for  $f(x) - u \in k(u)[x]$. So, we have  $C_f = \A^1_k(t), K_f= k(t),$ and $G_f = Z_2$. Since $k(t)/k$ is regular, \Cref{lem:g-f-spl-fld} says that the extension $\kk \inj k(\mbf{t}) \ceq k(t_1,\dotsc,t_r)$ defined by $u_i \ceq f(t_i), \, i=1,\dotsc,r$, makes $k(\mbf{t})$ into a splitting field for $\bb{m}(f(x))\in \kk[x]$. 
\end{exmp}


\begin{prop}\label{prop:spl-flt-rat}
$\kk_f$  is a rational function field over $\k$ if and only if $C_f/k$ is birationally isomorphic to $\p^1_k$ (equivalently, $K_f$ is a one-dimensional rational function field over $k$).
\end{prop}

\begin{proof}
The ``if'' direction is obvious. For the converse, suppose that ${\k}_f(C_f^{r})$ is a rational function field over $\k$, or equivalently, that $C_f^{r}$ is $\k$-birational to $\p^{r}_k$; this implies in particular that  $K_f/k$ is regular. We can find inside $C_f^{r}$ an  integral rational curve $U/{\k}$ which possesses a $\k$-point $P\in U({\k})$. For some $i=1,\dotsc,{r}$, the $i$-th projection $\pi_i: C_f^r \surj C_f$ restricts to a non-constant $\k$-morphism $\pi_i: U \too C_f$. Then, $C_f/{\k}$ is of genus $0$ and has the point $\pi_i(P) \in C_f({\k})$, which implies that $\ol{C}_f$ is $\k$-isomorphic to $\p^1_k$. 
\end{proof}

\begin{rem}
The proof above implies in fact that $\kk_f$ is \emph{contained} in a $2d$-dimensional rational function field over $\k$ if and only if ${C}_f/k$ is birational to $\p^1_{k}$.
\end{rem}

\begin{emp}\label{emp:field-inclusions}
Recall from \ref{emp:mestre} that $K$ denotes the field $k(\mbh)$, and $\mm$ the field $k(\mbm) = K(\mbq)$. Consider now the chain of field extensions
\begin{equation}
    \begin{tikzcd}
        k \arrow[r, hook] & K \arrow[r, hook] & \mm \arrow[r, hook] & \kk \arrow[r, hook] & \kk_f,
    \end{tikzcd}
\end{equation}
in which the first two extensions are each purely transcendental of degree $d+1$ and the last two are finite. The curve $\X_f$ (and all its variants which we consider) will be defined over $K$, and they will be equipped with $\kk_f$-points that are \emph{transcendental points} (cf. \ref{emp:transc}). The next lemma furnishes (via \Cref{lem:alg-closed-tors}) the key assumption \ref{item:tors} of \Cref{lem:Vn-phin}, which is used in the proof of \Cref{thm:dim-Vn} to determine the subgroups (of the respective Mordell--Weil groups) generated by these $\kk_f$-points.
\end{emp}

\begin{lem}\label{lem:alg-closed}
If $k$ is algebraically closed and $\chr k = 0$, then $K$ is algebraically closed in $\kk_f$.
\end{lem}

\begin{proof}
Suppose for a contradiction that there exists a finite extension $F/K$ which is contained in $\kk_f$. Let $S$ denote the integral closure of ${\k}[\mbh]$ in $F$, and let $\delta \in {\k}[\mbh]$ denote the discriminant of the ${\k}[\mbh]$-module $S$. Since $\chr k=0$, there are no non-trivial connected \'etale covers of $\A^{d+1}_k(\mbh)$, so $\delta$ is not a unit in ${\k}[\mbh]$. The coordinate ring $R_f^{\otimes r}$ of $C_f^{r}$ is integrally closed in its fraction field $\kk_f$, so $S$ is a subring of $R_f^{\otimes r}$. So, we have inclusions ${\k}[\mbm] = {\k}[\mbh][\mbq] \inj S[\mbq] \inj R_f^{\otimes r}$ which correspond to finite $\k$-morphisms
\begin{equation*}
        C_f^{r} \too \spec S[\mbq] \too \A^{r}_k(\mbm).
\end{equation*}
The branch locus $B'$ of the second morphism is cut out by $\delta \in {\k}[\mbm]$ (since it is obtained from the inclusion ${\k}[\mbh] \inj S$ by tensoring with the polynomial ring ${\k}[\mbq]$). On the other hand, the branch locus $B$ of $C_f^{r} \too \A^{r}_k(\mbm)$ is cut out by $\disc (\bb{m}\circ f) \in {\k}[\mbm]$. Since $B' \subset B$, $\delta$ must divide $\disc (\bb{m}\circ f)$ in ${\k}[\mbm]$. But \Cref{lem:disc-g-f} says that the irreducible factors of $\disc (\bb{m}\circ f)$ in ${\k}[\mbm]$ (up to scaling by constants in $\units{k}$) are $\disc \, \bb{m} = \disc (\bb{q}(x)^2 - \bb{h}(x))$, and the elements of the form $\bb{m}(f(\alpha)) = \bb{q}(f(\alpha))^2 - \bb{h}(\alpha)$, where $\alpha \in k$ is a root of $f'(x)$. None of these can equal $\delta$, since they do not lie in the subring ${\k}[\mbh] \subset {\k}[\mbm]$, so we have obtained a contradiction.
\end{proof}

\begin{rem}\label{rem:char-k-0}
The proof above makes critical use of the fact that $\chr k = 0$ (namely, we make use of the fact that $\A_k^{d+1}(\mbh)$ does not admit any non-trivial connected \'etale covers in characteristic zero). We fully expect, nevertheless, that the statement of \Cref{lem:alg-closed} continues to hold even if $\chr k > 0$. For example, if $f(x) = x^n$ then we can modify the argument as follows. We have in this situation $R_f^{\otimes r} \ceq k[\mathbold{t}] = k[t_1,\dotsc,t_r]$, with the inclusion $k[\mbu] \inj k[\mathbold{t}]$ defined by $\mbu = \mathbold{t}^n$. Then, the cover $C_f^n = \spec k[\mathbold{t}] \too \spec k[\mbm]$ is \emph{totally ramified} at the origin $P = (0,\dotsc,0)$. This implies that any finite cover of $\spec S \too \spec k[\mbh]$ arising as in the proof above must be ramified, and the rest of the proof carries through. 
\end{rem}

\section{Polynomials having cyclic and dihedral monodromy groups}\label{sec:mono}

From now, until the end of the article, \emph{$k$ denotes a field of characteristic zero.}

\begin{emp}\label{emp:fin-subgp-pgl2}
\Cref{prop:spl-flt-rat} invites consideration of the following problem: \emph{Produce polynomials $f(x) \in {\k}[x]$ for which $K_f$ is a one-variable rational function field $k(z)$}. We take as our starting point the observation that the monodromy group of such a polynomial must be a finite subgroup of $\aut_k(k(z)) = \pgl_2({\k})$. 

The determination of the finite subgroups of $\pgl_2(\ol{k})$, for $k$ of characteristic zero, is classical. Apart from the three sporadic groups $A_4, S_4,$ and $A_5$\footnote{Here, $A_n$ denotes the alternating group on $n$ letters}, which we will not consider in this article, $\pgl_2(\ol{k})$ contains the family of cyclic groups $\{Z_n \}$ and the family of dihedral groups $\{D_n \}$\footnote{Here, $D_n$ denotes the dihedral group of order $2n$ which admits the presentation $D_n = \langle \delta, \sigma \mid \sigma^n, \delta^2, (\sigma \delta)^2 \rangle$. In particular, $D_1 \cong Z_2$ and $D_2 \cong V_4$ are both abelian groups.}. For a given number field $k$, however, only finitely many of these occur in $\pgl_2(k)$. Indeed, Beauville~\cite{beauville} showed that $\pgl_2(k)$ contains $Z_n$ if and only if $k$ contains $\lambda_n \ceq \zeta_n + \zeta_n\inv$; this is the source of the requirement in \Cref{thm:main} that $k$ contain $\lambda_n, \lambda_{2n},$ or $\zeta_n$. 
\end{emp}

\begin{emp}\label{emp:cyc-mono}
\textbf{Cyclic monodromy.} If $k$ contains $\zeta \ceq \zeta_n$, then $\pgl_2(k)$ contains ${\sigma}_{n} \ceq \mattwo{\zeta}{0}{0}{1}$. 
Identifying $\pgl_2(k)$ with $\aut_k(k(\ucyc{n}))$, this corresponds to the automorphism $\ucyc{n} \mapsto \zeta \ucyc{n}$, and the cover 
\begin{equation}\label{eq:fcyc-p1}
    \fcyc{n}:\A^1_{k}(\ucyc{n}) \too \A^1_{k}(u), \qquad u \ceq (\ucyc{n})^n, 
\end{equation}
is a quotient by $\langle {\sigma}_{n} \rangle$. The inclusion $k(u) \inj k(\ucyc{n})$ yields the factorization
\begin{equation*}
    x^n - u = \prod_{j=1}^n (x - \sigma_n^j(\ucyc{n})) =  \prod_{j=1}^n \big(x - \zeta^j \ucyc{n}\big) \in k(\ucyc{n})[x].
\end{equation*}
So, \emph{if $k$ contains $\zeta_n$ and $f(x) \ceq x^n$, then $C_f = \A^1_k(\ucyc{n})$, $K_f = k(\ucyc{n})$ and $G_f = \langle {\sigma}_{n} \rangle \cong Z_n$}. 
\end{emp}

\begin{emp}\label{emp:dih-mono}
\textbf{Dihedral monodromy with $\zeta_n$.}  
Consider the element $\delta \ceq \mattwo{0}{1}{1}{0} \in \pgl_2(k)$. Identifying $\pgl_2(k)$ with $\aut_k(k(\udih{n}))$, $\delta$ corresponds to the involution $\udih{n} \mapsto (\udih{n})\inv$. If $k$ contains $\zeta \ceq \zeta_n$, we have again the automorphism ${\sigma}_{n}:\udih{n} \mapsto \zeta \udih{n}$, which together with $\delta$ generates a copy of $D_n$ inside $\aut_k(k(\udih{n}))$. To construct a quotient by this group, define 
\begin{equation*}
\tr(x) \ceq x + x\inv, \qquad  \fdih{n}(x)  \ceq x^n + x^{-n}.
\end{equation*}
Observe that $\fdih{n}(x) \in k[x,x\inv]$ is invariant under the automorphism $x \mapsto x\inv$, whose ring  of invariants is $k[\tr(x)]$. So, there exists some $\fRe{n}(x) \in k[x]$ which satisfies the identity
\begin{equation}\label{eq:dih-func-eq}
   \fRe{n}(\tr(x)) = \fdih{n}(x) = \tr(x^n) \in k[x,x\inv].
\end{equation}
The rational function $\fdih{n}$ has exactly two poles, each of order $n$, at the points $0$ and $\infty$. So, if we put $\ure{n} \ceq \tr(\udih{n})$ and $\udih{1} \ceq (\udih{n})^n$, and write  $\gmk(\udih{n}) \ceq \spec k[\udih{n},(\udih{n})\inv] = \p^1_k(\udih{n}) - \{0,\infty\}$, then we have a cover
\begin{equation}\label{eq:fdih-n}
    \fdih{n} : \gmk(\udih{n}) \too \A^1_{k}(u), \qquad u \ceq \fdih{n}(\udih{n}) = \fRe{n}(\ure{n}) = \tr(\udih{1}),
\end{equation}
which is a quotient by $\langle {\sigma}_{n}, \delta \rangle \cong D_n$. 
Applying $\sigma_n^j$ to \eqref{eq:dih-func-eq} has the effect of replacing $x$ with $\zeta ^j x$, and it follows that \eqref{eq:fdih-n} gives rise to the factorization
\begin{equation}\label{eq:fre-n-fact-udih}
    \begin{alignedat}{1}
        \fRe{n}(x) - u & = \prod_{j=1}^n \big(x - \sigma_n^j(\ure{n})) \\
        & = \prod_{j=1}^n \big(x - \big(\zeta^j \udih{n} + \big(\zeta^j \udih{n}\big)\inv\big)\big) \in k\big(\udih{n}\big)[x].
    \end{alignedat}
\end{equation}
So, \emph{for $n\geqslant 3$, if $\zeta_n \in k$ and $f\ceq \fRe{n}$, then $C_f = \gmk(\udih{n}), K_f = k(\udih{n})$ and $G_f = \langle \dih{\sigma}{n}, \delta \rangle \cong D_n$}. Note that if $n=1$ or $2$ then $k(\udih{n})/k(u)$ is  a quadratic extension of $k(\ure{n})$, which is the splitting field of $\fRe{n}(x) - u$; see \Cref{exmp:quad-mono}). 
\end{emp}

\begin{emp}
\textbf{Descent.}
In general, if $T/k$ is a commutative group scheme with a point $P\in T(k)$ of order $n$, then the translation-by-$P$ map $T \iso T$ induces an order $n$ automorphism $\ol{\sigma} \in \aut_k(k(T))$. This generates a dihedral group of order $2n$ with the involution $\delta \in \aut_k(k(T))$ corresponding to the inversion homomorphism $T \iso T$. The construction in \ref{emp:dih-mono} is the special case when $T = \gmk$, so the requirement that $k$ contain $\zeta_n$ boils down to the fact that the $n$-th division field of $\gmk$ (i.e. the extension of $k$ generated by its points of order $n$) is $k(\zeta_n)$.  

In the next subsection, we construct the splitting field for $\fRe{n}(x) - u$ with the more relaxed assumption $\lambda_n \in \units{k}$. The basic idea is to replace the $\gmk$ with an algebraic torus $T_n/k$ whose $n$-th division field is $k(\lambda_n) = k$. 
\end{emp}


\begin{emp}\label{emp:dih-descent}
\textbf{Dihedral monodromy with $\lambda_n$.} Fix $n\geqslant 3$, and put $\lambda \ceq \lambda_n$ and $\zeta \ceq \zeta_n$. For simplicity, let $k$ be the real cyclotomic field $\Q(\lambda)$. The assumption $n\geqslant 3$ is equivalent to saying that $\ell \ceq k(\zeta)$ is a quadratic extension of $k$. Let $\tau \in \gal(\ell/k)$ denote complex conjugation, and let $N_{\ell/k}:\units{\ell} \to \units{k}$ denote the usual norm homomorphism defined by  $N_{\ell/k}(\alpha) \ceq \alpha \tau(\alpha)$. 
Define
\begin{equation*}
    \dih{\ol{u}}{n} \ceq s - \zeta t \in \ell[s,t].
\end{equation*}
The \defi{norm-torus} associated to $\ell/k$ is the affine plane curve
\begin{equation}
    T_n/k:\; N_{\ell/k}(\dih{\ol{u}}{n}) = s^2 - \lambda st + t^2 = 1.
\end{equation}
It is a one-dimensional algebraic $k$-torus (i.e. a $k$-group scheme which becomes isomorphic to $\mbf{G}_{m,\ol{k}}$ upon base-change to $\ol{k}$). We have a natural group isomorphism
\begin{equation*}
\begin{alignedat}{1}
    T_n(k) & \iso N^1_{\ell/k} \ceq \{ \alpha \in \units{\ell} \mid N_{\ell/k}(\alpha) = 1\},\\
    (s_0,t_0) & \mtoo s_0 -\zeta t_0.
\end{alignedat}
\end{equation*}
The product of two $k$-points $\alpha$ and $\beta$ is the point $(s_0,t_0) \in T_n(k)$ determined by the identity $\alpha \beta = s_0 - \zeta t_0 \in N^1_{\ell/k}$. In particular, $(0,-1) \in T_n(k)$ corresponds to $\zeta \in N^1_{\ell/k}$, so we have $\zeta (s_0 - \zeta t_0) = t_0 - \zeta(\lambda t_0 - s_0)$. So, the multiplication-by-$\zeta$ map on $T_n$ induces the order $n$ automorphism $\ol{\sigma}_n \in \aut_k(k(T_n))$ defined by
\begin{equation}\label{eq:ol-sigma}
    \ol{\sigma}_n(s) = t \qquad \textup{ and } \qquad  \ol{\sigma}_n(t) = \lambda t - s. 
\end{equation}
By construction, the powers of $\ol{\sigma}_n$ satisfy 
\begin{equation}\label{eq:ol-sigma-powers}
    \ol{\sigma}_n^j(\dih{\ol{u}}{n}) = \zeta^j \dih{\ol{u}}{n}  \in \ell(T_n), \quad \textup{ for } j=1,\dotsc,n.
\end{equation}
Since $N_{\ell/k}(\dih{\ol{u}}{n}) = \dih{\ol{u}}{n} (s - \zeta\inv t) = 1 \in \ell(T_n)$, we have 
\begin{equation*}
    (\dih{\ol{u}}{n}) \inv = s-\zeta \inv t = \tau(\dih{\ol{u}}{n}) \in \ell(T_n),
\end{equation*}
and it follows that $\tr(\dih{\ol{u}}{n}) = \tr_{\ell/k}(\dih{\ol{u}}{n})$ is an element in $k(T_n)$. We put 
\begin{equation*}
    \Re{\ol{u}}{n} \ceq \tr(\dih{\ol{u}}{n}) = 2s - \lambda t \in k(T_n),
\end{equation*}
so that, for $j=1,\dotsc,n,$ we have the identity 
\begin{equation*}
    \ol{\sigma}_n^j(\Re{\ol{u}}{n}) = \tr (\zeta^j (\dih{\ol{u}}{n})) = \tr(\zeta^j)s - \tr(\zeta^{j+1})t \in k(T_n).
\end{equation*}
Now, the functional equation \eqref{eq:dih-func-eq} implies that $\fRe{n}(\Re{\ol{u}}{n}) = \fdih{n}(\dih{\ol{u}}{n}) \in k(T_n)$. Similar to the situation in \ref{emp:dih-mono}, each automorphism $\ol{\sigma}_n^j$ fixes $\fdih{n}(\dih{\ol{u}}{n}) \in k(T_n)$, so applying $\ol{\sigma}_n^j$ to this relation, we find that that the cover 
\begin{equation}\label{eq:Tn-A1}
\begin{alignedat}{1}
    \dih{{f}}{n}: T_n &\too \A^1_k(u) \\
    u & \ceq \fRe{n}(\Re{\ol{u}}{n}) = \fdih{n}(\dih{\ol{u}}{n})
\end{alignedat}
\end{equation}
yields the factorization
\begin{equation}\label{eq:fRe-real-fact}
    \fRe{n}(x) - u = \prod_{j=1}^n (x - \ol{\sigma}_n^j(\Re{\ol{u}}{n})) =  \prod_{j=1}^n \big(x - (\tr(\zeta^j)s - \tr(\zeta^{j+1})t) \big) \in k(T_n)[x].
\end{equation}
\emph{We conclude that for $n\geqslant 3$, if $k$ contains $\lambda_n$ and  $f = \fRe{n}(x)$, then $C_f = T_n$ and $K_f = k(T_n)$}. We will later need the following more precise statement. 
\end{emp}

\begin{lem}
Let $n\geqslant 3$, assume that $k$ contains $\lambda_n$, and put $\ell = k(\zeta_n)$. Then,
\begin{equation}\label{eq:iso-phi-u-ol}
    \phi: \ell(\udih{n}) \iso \ell(T_n), \qquad \udih{n} \mtoo \dih{\ol{u}}{n}
\end{equation}
is an isomorphism of $\ell(u)$-algebras which sends $\sigma_n^j(\ure{n}) \in \ell(\udih{n})$ to the element $\ol{\sigma}_n^j(\Re{\ol{u}}{n}) \in \ell(T_n)$, for $j=1,\dotsc,n$. Thus, $\phi$  transforms the factorization \eqref{eq:fRe-real-fact} into the factorization \eqref{eq:fre-n-fact-udih}.
\end{lem}

\begin{proof}
We have $\phi((\udih{n})\inv)=(\dih{\ol{u}}{n})\inv = s-\zeta\inv t \in \ell(T_n)$. Using this, one checks that the map
\begin{equation*}
    s \mtoo \dfrac{\zeta \inv \udih{n} - \zeta (\udih{n})\inv}{\zeta - \zeta \inv}, \qquad t \mtoo \dfrac{(\udih{n})\inv - \udih{n}}{\zeta - \zeta \inv}.
\end{equation*}
defines the inverse $\phi\inv: \ell(T_n) \iso \ell(\udih{n})$. That $\phi$ is an isomorphism of $\ell(u)$-algebras follows because $\phi(u) = \phi(\fdih{n}(\udih{n})) = \fdih{n}(\dih{\ol{u}}{n}) = u \in \ell(T_n)$. The $Z_n$-actions on both fields are compatible with $\phi$ in the sense that $\phi(\sigma_n(\udih{n})) = \zeta_n \phi(\udih{n}) = \zeta (\dih{\ol{u}}{n}) = \ol{\sigma}_n(\dih{\ol{u}}{n})$. Thus, $\phi$ transforms each element $\sigma_n^j(\ure{n}) = \tr(\zeta^j \udih{n})$ into the corresponding element $\ol{\sigma}_n^j(\Re{\ol{u}}{n}) = \tr( \zeta^j (\dih{\ol{u}}{n}))$.
\end{proof}

\begin{exmp}
Let us illustrate the factorization \eqref{eq:fRe-real-fact} in the cases $n=3,4$ and $6$ (namely, the cases used to obtain \Cref{cor:bounds-Q}). Since $\lambda_3 = -1,\lambda_4 = 0,$ and $\lambda_6 = 1$, we may work over $k\ceq \Q$. Then, the polynomials $\fRe{n}(x)$ and the associated tori $T_n/\Q$ are given by
\begin{equation*}
    \begin{alignedat}{3}
        \fRe{3}(x) & =  x^3 - 3x, &&\qquad   T_3/\Q &&:\; s^2 + st + t^2 = 1,\\
        \fRe{4}(x) & = x^4 - 4x^2 + 2, && \qquad T_4/\Q &&:\; s^2 + t^2 = 1,\\
        \fRe{6}(x) & = x^6 - 6x^4 + 9x^2 - 2,  &&\qquad  T_6/\Q &&:\; s^2 - st + t^2 = 1.
    \end{alignedat}
\end{equation*}
Defining the inclusions $\Q(u) \inj \Q(T_n)$ as in \eqref{eq:Tn-A1}, the factorizations in \eqref{eq:fRe-real-fact} become:
\begin{equation*}
\begin{alignedat}{1}
    \fRe{3}(x) -u &= (x + s + t)(x - s + 2t)(x - t + 2s),\\
    \fRe{4}(x) -u &= (x - 2s)(x + 2s)(x-2t)(x+2t),\\
    \fRe{6}(x) -u &= (x + s + t)(x-s-t)(x - s + 2t)(x+s-2t)(x - t + 2s)(x+t-2s).
\end{alignedat}
\end{equation*}
\end{exmp}

\begin{rem}
We can construct a birational map between $\p^1_k(z)$ and $T_n$ by stereographic projection from the point $(0,1) \in T_n$. That is, we take $z$ to be the slope of a line passing through $(0,1)$ and the generic point $(s,t) \in T_n$, so that $z = (t-1)/s \in k(T_n)$. Replacing $t$ with $sz+1$ in the equation defining $T_n$ and solving for $s$ and $t$ in terms of $z$, we obtain
\begin{equation*}
    s = \dfrac{\lambda - 2z}{z^2 - \lambda z +1}, \qquad t = \dfrac{1-z^2}{z^2 -\lambda z + 1},
\end{equation*}
so $k(\mbf{T}_n) = k(z)$. One can check that $\ol{\sigma}_n \in \aut_k(k(T_n))$ corresponds to the order $n$ automorphism
\begin{equation*}
    z \mtoo \dfrac{(\lambda + 1)z-1}{z+1} \in \aut_k(k(z)).
\end{equation*}
This is in fact the automorphism mentioned by Beauville in \cite{beauville} to conclude that $\pgl_2(k)$ contains $Z_n$ if $\lambda_n \in k$. Finally, we note that  
\begin{equation*}
    \dih{\ol{u}}{n} = {s-\zeta t} = \dfrac{\zeta z - 1}{z - \zeta} \in \ell(T_n) = \ell(z).
\end{equation*}
Thus, the automorphism $\theta \ceq \mattwo{\zeta}{-1}{1}{-\zeta} \in \pgl_2(\ell)$ changes the coordinate on $\p^1_{\ell}$ from $z$ to $u$,  and $\ol{\sigma}_n \ceq \mattwo{\lambda +1}{-1}{1}{1}$ and $\sigma_n \ceq \mattwo{\zeta}{0}{0}{1}$ are conjugate to each other by $\theta$, i.e. we have ${\sigma}_n = \theta \comp \ol{\sigma}_n \comp \theta \inv \in \pgl_2(\ell)$. The reader can moreover verify that the involution $\delta \ceq \mattwo{0}{1}{1}{0}$ commutes with $\theta$, so we have $\langle \sigma_n, \delta \rangle = \theta \comp \langle \ol{\sigma}_n, \delta \rangle \theta \inv \subset \pgl_2(\ell)$.
\end{rem}

\begin{rem}\label{rem:char-k-neq-2n}
If $\chr k > 0$, then the entire discussion in \ref{emp:cyc-mono}, \ref{emp:dih-mono} and \ref{emp:dih-descent} holds verbatim so long as $\chr k$ does not divide $2n$. 
\end{rem}

 \section{The rational function fields \texorpdfstring{$\kkcyc{n}, \kkdih{n},$}{} and \texorpdfstring{$\kkre{n}$}{}.}\label{sec:kkn}

 \begin{emp}\label{emp:sec-kkn}
\textbf{Conventions.} We assume for this section that $\chr k = 0$. In what follows, the phrase ``for each $i$'' means ``for $i=1,\dotsc,r$''. If $\mbf{t} = t_1,\dotsc,t_r$ and $\mbf{s} = s_1,\dotsc,s_r$ are $r$-tuples of indeterminates and $q(x) \in k(x)$ is a rational function, then we write $\mbf{t} = q(\mbf{s})$ to abbreviate the  collection of identities $t_i = q(s_i),$ for each $i$. We use the notation from \ref{emp:mestre}; in particular $r=2d+2$ for a fixed positive integer $d$, and $\kk$ denotes the rational function field $k(\mbu) = k(u_1,\dotsc,u_r)$. 
\end{emp}

\begin{emp}\label{emp:fields-Ln}
\textbf{The fields $\kkcyc{n}, \kkdih{n},$ and $\kkre{n}$.} 
Let $n\in \N$. Let $\mbucyc{n}$ and $\mbudih{n}$ denote the tuples of indeterminates $\ucyc{1,n}, \dotsc , \ucyc{r,n}$ and $\udih{1,n}, \dotsc , \udih{r,n}$. Let $\mbure{n}$ denote the $r$-tuple $\ure{1,n},\dotsc,\ure{r,n}$ defined by
\begin{equation*}
    \mbure{n} \ceq \tr (\mbudih{n} ) = (\mbudih{n}) + (\mbudih{n})\inv
\end{equation*}
Next, if $n\geqslant 3$ and $k$ contains $\lambda_n$, let $\mbf{s}$ and $\mbf{t}$ denote the tuples $s_1,\dotsc,s_r$ and $t_1,\dotsc,t_r$, and put $\dih{\ol{\mbu}}{n} \ceq \mbf{s} - \zeta_n \mbf{t}$, so that
\begin{equation*}
    \mbf{T}_n/k: \; \{ \mbf{s}^2 - \lambda_n \mbf{s} \mbf{t} + \mbf{t}^2 = 1\} \; \subset \A^{2r}_k(\mbf{s},\mbf{t})
\end{equation*}
can naturally be identified with the $r$-fold product of the torus $T_n/k$ from \ref{emp:dih-descent}. Let $\Re{\ol{\mbu}}{n}$ denote the $r$-tuple $\Re{\ol{u}}{1,n} , \dotsc , \Re{\ol{u}}{r,n}$ defined by
\begin{equation*}
    \Re{\ol{\mbu}}{n} \ceq \tr(\dih{\ol{\mbu}}{n}) = 2\mbf{s} - \lambda_n \mbf{t}.
\end{equation*}
Now, define the following extensions of rational function fields, (the third extension is defined below only for $n\geqslant 3$ and $k$ containing $\lambda_n$; in \ref{emp:direct-system}, we complete the picture by defining also $\kkre{2}$): 
\begin{equation}\label{eq:kk-inj-kkbul}
\begin{alignedat}{5}
    \kk & \inj \kkcyc{n}   && \ceq k(\mbucyc{n}), \qquad && \mbu && \ceq (\mbucyc{n})^n,\\
    \kk & \inj \kkdih{n}   && \ceq k(\mbudih{n}), \qquad && \mbu && \ceq \fRe{n}(\mbure{n}) = \fdih{n}(\mbudih{n}),\\
    \kk & \inj \kkre{n}    && \ceq k(\mbf{T}_n), \qquad && \mbu && \ceq \fRe{n}(\Re{\ol{\mbu}}{n}) = \fdih{n}(\dih{\ol{\mbu}}{n}).
\end{alignedat}
\end{equation}
Note in particular that $\mbu = \mbucyc{1} = \mbure{1}$, and hence, we have $\kk = \kkcyc{1}= k(\mbure{1})$. 
\end{emp}

\begin{emp}\label{emp:aut-kkbul}
\textbf{Automorphisms.}
If $k$ contains $\zeta_n$ and $\kk_n = k(\mbu_n)$ denotes $\kkcyc{n} = k(\mbucyc{n})$ or $\kkdih{n} = k(\mbudih{n})$, then define for each $i$ the automorphism
\begin{equation}\label{eq:sigma-i-n}
    {\sigma}_{i,n} : u_{e,n} \mtoo 
    \ccases{
        \zeta_n u_{e,n} & \textup{ if }\; e = i,\\
        u_{e,n} & \textup{ if } \; e \neq i.
    }  \quad \in \aut_k(\kk_n).
    \end{equation}    
For the field $\kkdih{n}$, define also for each $i$ the automorphism
\begin{equation}\label{eq:delta-i-n}
    {\delta}_{i,n} : \udih{e,n} \mtoo 
    \ccases{
        \big(\udih{e,n}\big)\inv & \textup{ if }\; e = i,\\
        \udih{e,n} & \textup{ if } \; e \neq i.
    }  \quad \in \aut_k(\kkdih{n}).
    \end{equation}
Finally, if $n\geqslant 3$ and $k$ contains $\lambda_n$, then define for each $i$ the automorphism
\begin{equation}\label{eq:sigma-i-n-re}
        \ol{\sigma}_{i,n} : \dih{\ol{u}}{e,n} \mtoo \ccases{
        \zeta_n \dih{\ol{u}}{e,n} & \textup{ if }\; e = i,\\
        \dih{\ol{u}}{e,n} & \textup{ if } \; e \neq i.
        }  \quad \in \aut_k(\kkre{n}).
    \end{equation}  
\end{emp}

\begin{prop}\label{prop:Ln-fact}
Let $n\in \N$.
\begin{enumerate}[(a),leftmargin=0pt,itemindent=1cm]
    \item If $k$ contains $\zeta_n$, then we have the factorizations
    \begin{equation}\label{eq:m-fcyc-fRe-fact}
    \begin{alignedat}{3}
        \bb{m}(x^n) & = \prod_{i=1}^r \big(x^n - u_i \big) && =  \prod_{i=1}^r \prod_{j=1}^n \big(x - {\sigma}_{i,n}^j (\ucyc{i,n} ) \big) && \in \kkcyc{n}[x],\\
        \bb{m}(\fRe{n}(x)) & = \prod_{i=1}^r \big(\fRe{n}(x) - u_i \big) && = \prod_{i=1}^r \prod_{j=1}^n \big(x - {\sigma}_{i,n}^j (\ure{i,n})\big) && \in \kkdih{n}[x], \\
    \end{alignedat}
    \end{equation}
    \item If $n\geqslant 3$ and $k$ contains $\lambda_n$, then we have the factorization
    \begin{equation}\label{eq:m-fRe-fact}
        \bb{m}(\fRe{n}(x)) = \prod_{i=1}^r \big(\fRe{n}(x) - u_i \big) = \prod_{i=1}^r \prod_{j=1}^n \big(x - \ol{\sigma}_{i,n}^j (\Re{\ol{u}}{i,n} )\big) \in \kkre{n}[x].
    \end{equation}
    If moreover $k$ contains $\zeta_n$, then the map
    \begin{equation}\label{eq:dihLn-iso-ReLn}
        \phi_n: \kkdih{n} \iso \kkre{n}, \qquad \mbudih{n} \mtoo \dih{\ol{\mbu}}{n},
    \end{equation}
    is a $\kk$-algebra isomorphism which sends 
    \begin{equation*}
        {\sigma}_{i,n}^j (\ure{i,n}) \mtoo \ol{\sigma}_{i,n}^j(\Re{\ol{u}}{i,n}), \qquad \textup{ for } i=1,\dotsc,r.
    \end{equation*}
    So, \eqref{eq:dihLn-iso-ReLn} transforms the factorization \eqref{eq:m-fRe-fact} into the second line of \eqref{eq:m-fcyc-fRe-fact}.
    \item If $k$ is algebraically closed, then $K = k(\mbh)$ is algebraically closed in $\kkcyc{n}$ and $\kkdih{n}$. Consequently, for any abelian variety $\J/K$ we have an equality of torsion subgroups
    \begin{equation*}
        \J(K)_{\tors} = \J(\kkcyc{n})_{\tors} = \J(\kkdih{n})_{\tors}.
    \end{equation*}
\end{enumerate}
\end{prop}

\begin{proof}
As in \ref{emp:const-spl-fld}, identify the three extensions \eqref{eq:kk-inj-kkbul} with the $r$-th tensor power of the extensions $k(u) \inj k(\ucyc{n})$ from \ref{emp:cyc-mono}, $k(u) \inj k(\udih{n})$ from \ref{emp:dih-mono}, and $k(u) \inj k(T_n)$ from \ref{emp:dih-descent}, respectively. \Cref{lem:g-f-spl-fld} then implies that these give splitting fields for the respective $\bb{m}(f_n(x))$'s, and the factorizations \eqref{eq:m-fcyc-fRe-fact} and \eqref{eq:m-fRe-fact} hold by construction. The isomorphism \eqref{eq:dihLn-iso-ReLn} is simply the $r$-th tensor power (over $k$) of the isomorphism \eqref{eq:iso-phi-u-ol}. For part (c), assume that $k$ is algebraically closed. Then, the extension $\kkcyc{n}/\kk$ is nothing but the extension $\kk_f/\kk$ for $f= \fcyc{n}$, so \Cref{lem:alg-closed} implies that $K$ is algebraically closed in $\kkcyc{n}$.  Similarly, for $n\geqslant 3$, $\kkdih{n}/\kk$ is the extension $\kk_f/\kk$ for $f \ceq \fRe{n}$, and since $\kk \subset \kkdih{2} \subset \kkdih{4}$, \Cref{lem:alg-closed} implies that $K$ is algebraically closed in $\kkdih{n}$ for all $n\in \N$. The statement about the torsion subgroups then follows from \Cref{lem:alg-closed-tors}.
\end{proof}

\begin{emp}\label{emp:direct-system}
\textbf{Directed systems.} Arrange the fields $\{\kkcyc{n}\}_{n \in \N}$ and $\{\kkdih{n}\}_{n \in \N}$ into directed systems by defining for all $n\in \N$ and all $m$ dividing $n$ the inclusion
\begin{equation}\label{eq:kkm-kkn}
    \begin{alignedat}{3}
        \kkcyc{m} & \inj \kkcyc{n}, \qquad && \mbucyc{m} && \ceq (\mbucyc{n})^{n/m},\\
        \kkcyc{m} & \inj \kkdih{n}, \qquad && \mbudih{m} && \ceq (\mbudih{n})^{n/m}.
    \end{alignedat}
\end{equation}
Using the functional equation \eqref{eq:dih-func-eq}, one checks that the above inclusions give (for $m$ dividing $n$):
\begin{equation}\label{eq:mbure-m-n}
    \mbure{m} = \fRe{n/m}(\mbure{n}) \in \kkdih{n}. 
\end{equation}
We will have need later for the $r$-tuple $\mbuim{n} \ceq \uim{1,n},\dotsc,\uim{r,n}$ defined by 
\begin{equation*}
    \mbuim{n} \ceq (\mbudih{n}) -(\mbudih{n})\inv \in \kkdih{n}.
\end{equation*}
If $\fIm{n}(x) \in \Z[x]$ is the polynomial  determined by the identity
\begin{equation}\label{eq:fim-func-eq}
    \fIm{n}(\tr(x)) = \dfrac{x^n - x^{-n}}{x - x\inv} \in \Z[x,x\inv],
\end{equation}
then the following identities hold for every $m\mid n \in \N$:
\begin{equation}\label{eq:uim}
{\mbuim{m}} = ({\mbuim{n}}) \fIm{n/m}(\mbure{n}), \qquad (\mbuim{n})^2 = (\mbure{n})^2 - 4.
\end{equation}

Next, for $n\geqslant 3$ and $k$ containing $\lambda_{2n}$, we embed $\kkre{n}$ in $\kkre{2n}$ by sending the tuples $\mbf{s}$ and $\mbf{t}$ in $\kkre{n}$ to the tuples $\mbf{s}'$ and $\mbf{t}'$ in $\kkre{2n}$, respectively, which give the identity
\begin{equation*}
(\dih{\ol{\mbu}}{2n})^2 = \mbf{s}' - \zeta_n \mbf{t}' \in \kkre{2n}(\zeta_n).
\end{equation*}
In this way, we identify $\kkre{n}$ with the subfield of $\kkre{2n}$ generated over $k$ by these $r$-tuples $\mbf{s}'$ and $\mbf{t}'$. To complete the picture with the case $n=2$, we similarly \emph{define} the $r$-tuple of elements  $\dih{\ol{\mbu}}{2} \ceq \mbf{s}' - \zeta_4 \mbf{t}'  \in \kkre{4}$ by the identity $\dih{\ol{\mbu}}{2} \ceq (\dih{\ol{\mbu}}{4})^2 \in \kkdih{4}$, and we put 
\begin{equation*}
\kkre{2} \ceq k(\mbf{s}',\mbf{t}') \subset \kkre{4}.
\end{equation*}
Defining an isomorphism $\phi_2: \kkdih{2} \iso \kkre{2}$ by $\mbudih{2} \mtoo \dih{\ol{\mbu}}{2}$, the upshot of defining these inclusions is that we obtain for all $n\geqslant 2$ a commutative diagram
\begin{equation*}
\begin{tikzcd}
\kkdih{2n} \arrow[r, "\phi_{2n}"] & \kkre{2n}\\
\kkre{2n} \arrow[u, hook] \arrow[r, "\phi_n"] & \kkre{n} \arrow[u, hook],
\end{tikzcd}
\end{equation*}
which in turn gives for all $n\geqslant 2$ (using the functional equation \eqref{eq:dih-func-eq}) the identity $\fRe{2}(\Re{\ol{\mbu}}{2n}) = \Re{\ol{\mbu}}{n} \in \kkre{2n}$, equivalently, 
\begin{equation}\label{eq:re-ol-mbu-2}
(\Re{\ol{\mbu}}{2n})^2 = \Re{\ol{\mbu}}{n} + 2 \in \kkre{2n}. 
\end{equation}
\end{emp}

\begin{rem}
Exactly as in \Cref{rem:char-k-neq-2n}, if $\chr k$ is positive (and distinct from $2$), then essentially everything in this section continues to hold if we restrict our attention to $n$ not divisible by $\chr k$. A notable (possible) exception to this is the statement of \Cref{prop:Ln-fact} (c) for the field  $\kkdih{n}$, since the latter statement relies on \Cref{lem:alg-closed}, which in turn uses that $\chr k = 0$ (note:  \Cref{rem:char-k-0} for an alternative argument which works for the cyclic case). 
\end{rem}



\section{Defining the curves and stating the point and rank bounds}\label{sec:curves}

\begin{emp}\label{emp:define-curves}
Recall the notation from \ref{emp:fields-Ln}. Following the convention in \ref{emp:sec-kkn}, we put 
\begin{equation}\label{eq:mbv}
    \mbv \ceq \bb{q}(\mbu) \in \kk,
\end{equation}
so that the identity $\bb{m}(x) = \bb{h}(x)^2 - \bb{r}(x)$ from \Cref{lem:sqroot} yields $\mbv^2 = \bb{h}(\mbu) \in \kk$. Consider now the following table of data; below, $\left[ nd \right] \in \{0,1\}$ denotes the remainder of $nd$ mod $2$ (so that $nd+\left[ nd \right]$ is always even). 
\begin{table}[H]
\centering
\renewcommand{\arraystretch}{1.7}
\setlength{\tabcolsep}{5pt}
\adjustbox{max width=\textwidth}{
\caption{Main table of data}\label{tab:data} 
\begin{tabular}{|c||r|l|c|c|c|c|c|}
\hline
& & & & \multicolumn{2}{c|}{$nd$ is even} & \multicolumn{2}{c|}{$nd$ is odd}\\
\hline
$\bullet$  & {$\bbfbul{n}(x)$} & {$\mbQbul{n}$} & $\kk_n$ & $\bul{b}{n}$ & $nd$  & $\bul{b}{n}$ & $nd$\\
\hhline{|====|==|==|}
\textsf{cyc} & $\bb{h}(x^n)$ & $\big(\mbucyc{n},\mbv \big)$ & $\kkcyc{n}$ & $0$ & $2\gcyc{n}+2$  & $1$ & $2\gcyc{n}+1$\\ 
\hline
{\textsf{re}} & {$\bb{h}(\fRe{n}(x))$} & $\big(\mbure{n},\mbv\big)$ & \multirow{3}{*}{$\kkdih{n}$} & {$0$} & {$2\gre{n}+2$}  & {$1$} & {$2\gre{n}+1$}\\
\cline{1-3} \cline{5-8}
\textsf{im} & $(x^2 - 4)\bb{h}(\fRe{n}(x))$ & $\big(\mbure{n},\mbv \mbuim{n}\big)$ & & $2$ & $2\gim{n}$  & $3$ & $2\gim{n}-1$\\
\cline{1-3} \cline{5-8}
\textsf{dih} & $x^{nd+\left[ nd \right]}\bb{h}(\fdih{n}(x))$ & $\big(\mbudih{n},\mbv(\mbudih{n})^{\frac{nd+\left[ nd \right]}{2}}\big)$ & & $0$ & $\gdih{n}+1$  & $2$ & $\gdih{n}$\\
\hline \hline
$\bullet$  & {$\wtbbfbul{n}(x)$} & {$\wtmbQbul{n}$} & $\kk_{2n}$ & $\bul{\wt{b}}{n}$ & $nd$  & $\bul{\wt{b}}{n}$ & $nd$\\
\hhline{|====|==|==|}
\textsf{cyc} & $x\bb{h}(x^n)$ & $\big(\mbucyc{n},\mbv \mbucyc{2n}\big)$ & $\kkcyc{2n}$ & $2$ & $2\gwtcyc{n}$  & $1$ & $2\gwtcyc{n} + 1$ \\ 
\hline
{\textsf{re}} & {$(x+2)\bb{h}(\fRe{n}(x))$} & $\big(\mbure{n},\mbv \mbure{2n}\big)$ & \multirow{3}{*}{$\kkdih{2n}$} & {$2$} & {$2\gwtre{n}$} & {$1$} & {$2\gwtre{n}+1$}\\
\cline{1-3} \cline{5-8}
\textsf{im} & $(x-2)\bb{h}(\fRe{n}(x))$ & $\big(\mbure{n},\mbv \mbuim{2n}\big)$ & & $2$ & $2\gwtim{n}$  & $1$ & $2\gwtim{n}+1$\\
\cline{1-3} \cline{5-8}
\textsf{dih} & $x^{nd-\left[ nd \right]+1}\bb{h}(\fdih{n}(x))$ & $\big(\mbudih{n},\mbv  \mbudih{2n}(\mbudih{n})^{\frac{nd-\left[ nd \right]}{2}}\big)$ & & $2$ & $\gwtdih{n}$ & $0$ & $\gwtdih{n}+1$\\
\hline
\end{tabular}}
\end{table}
The four values of $\bullet$ are short-forms, respectively, for the labels \textbf{cyc}lic, \textbf{re}al, \textbf{im}aginary, and \textbf{dih}edral. For each of these labels (and $n\in \N$), we define a pair of hyperelliptic curves using the polynomials $\bbfbul{n}(x)$ and $\wtbbfbul{n}(x)$ in the second column:
\begin{equation}\label{eq:Xbul-Xwtbul}
        \Xbul{n}/K: \; y^2 = \bbfbul{n}(x), \qquad \textup{and} \qquad \Xwtbul{n}/K: \; y^2 = \wtbbfbul{n}(x).
\end{equation}
It is an exercise (using \Cref{lem:disc-g-f}) to see that, in each case, $\bbfbul{n}$ and $\wtbbfbul{n}$ are separable polynomials. Using this fact and the genus-degree formula for hyperelliptic curves, we have expressed in the last and third-last column of \Cref{tab:data} the integer $nd$ in terms of the genera $\gbul{n}$ and $\gwtbul{n}$ of these curves. 
We explain the remaining data of the table in the next proposition.
\end{emp}


\begin{prop}\label{prop:table-data}
Fix $n\in \N$ and $\bullet \in \{\textup{\textsf{cyc}, \textsf{re}, \textsf{im}, \textsf{dih}}\}$. 
\begin{enumerate}[(a),itemindent=1cm, leftmargin=0pt]
    \item The integers $\bul{b}{n}$ and $\bul{\wt{b}}{n}$ denote the number of $K$-rational Weierstrass points of $\Xbul{n}$ and $\Xwtbul{n}$, respectively.
    \item The third column of \Cref{tab:data} defines an $r$-tuple of points
    \begin{equation}
        \begin{alignedat}{4}
            \mbQbul{n} & \ceq \Qbul{1,n} , \dotsc , \Qbul{r,n} && \in \;&& \Xbul{n}(\kk_n), \quad && \textup{ if } \zeta_n \in k,\\
            \wtmbQbul{n} & \ceq \wtQbul{1,n} , \dotsc , \wtQbul{r,n} && \in \; && \Xwtbul{n}(\kk_{2n}), \quad && \textup{ if } \zeta_{2n} \in k.
        \end{alignedat}
    \end{equation}
    \item If $n\geqslant 2$, then we have $r$-tuples of points
    \begin{equation}\label{eq:Re-mbP-points}
    \begin{alignedat}{5}
        &\Re{\mbP}{n} && \ceq (\Re{\ol{\mbu}}{n}, \mbv) && \in \Xre{n}(\kkre{n}), \quad && \textup{ if } \lambda_n \in k,\\
        &\Re{\wt{\mbP}}{n} && \ceq (\Re{\ol{\mbu}}{n}, \mbv \Re{\ol{\mbu}}{2n}) && \in \Xwtre{n}(\kkre{2n}), \quad && \textup{ if } \lambda_{2n} \in k.
    \end{alignedat}
    \end{equation}
   Moreover, if $k$ contains $\zeta_n$ (resp. $\zeta_{2n}$), the isomorphism $\kkdih{n} \iso \kkre{n}$ (resp. $\kkdih{2n} \iso \kkre{2n}$) from  \eqref{eq:dihLn-iso-ReLn} induces bijections (of sets of points and $r$-tuples):
    \begin{equation*}
    \begin{alignedat}{5}
        \Xre{n}(\kkdih{n}) & \iso \Xre{n}(\kkre{n}), \qquad && (\mbure{n},\mbv) && \coloneqq \Re{\mbP}{n} && \mtoo \mbQre{n} && = (\Re{\ol{\mbu}}{n},\mbv) ,\\
        \Xwtre{n}(\kkdih{2n}) & \iso \Xwtre{n}(\kkre{2n}), \qquad && (\mbure{n},\mbv \mbure{2n}) && \coloneqq \Re{\wt{\mbP}}{n} && \mtoo \wtmbQre{n} && = (\Re{\ol{\mbu}}{n},\mbv \Re{\ol{\mbu}}{2n}).
    \end{alignedat}
    \end{equation*}
\end{enumerate}
\end{prop}

\begin{proof}
Part (a) follows by verifying in each case that $\bul{b}{n}$ (resp. $\bul{\wt{b}}{n}$) equals the number of $K$-rational roots of the polynomial $\bbfbul{n}(x)$ (resp. $\wtbbfbul{n}(x)$), plus one if $\deg \, \bbfbul{n}(x)$ (resp. $\deg \wtbbfbul{n}(x)$) is odd; we leave this to the reader.

We turn now to parts (b) and (c). By \eqref{eq:kk-inj-kkbul}, we have the identities
    \begin{equation*}
        \mbu = \fcyc{n}(\mbucyc{n}) = \fRe{n}(\mbure{n}) =\fRe{n}(\Re{\ol{\mbu}}{n}) = \fdih{n}(\mbudih{n}).
    \end{equation*}
    Plugging these into $\bb{h}(x)$ and using the fact that $\bb{h}(\mbu) = \mbv^2$ immediately yields
    \begin{equation}\label{eq:bbfbul-u-v2}
        \mbv^2 = \bbfcyc{n}(\mbucyc{n}) = \bbfRe{n}(\mbure{n}) = \bbfRe{n}(\Re{\ol{\mbu}}{n}) = \bb{h}(\fdih{n}(\mbudih{n})),
    \end{equation}
    which proves that the formulas for $\mbQcyc{n},\mbQre{n},\Re{\mbP}{n},$ and $\mbQdih{n}$ indeed define $r$-tuples of points in $\Xcyc{n}(\kkcyc{n}), \Xre{n}(\kkdih{n}), \Xre{n}(\kkre{n})$, and $\Xdih{n}(\kkdih{n})$,  respectively. Recalling from \eqref{eq:uim} the identity $(\mbuim{n})^2 = (\mbure{n})^2 - 4$, we have that
    \begin{equation*}
        (\mbv\mbuim{n})^2 = \bbfRe{n}(\mbure{n}) ((\mbure{n})^2 - 4) = \bbfIm{n}(\mbure{n}) \in \kkdih{n},
    \end{equation*}
    which verifies that $\mbQim{n} \ceq (\mbure{n},\mbv \mbuim{n})$ defines an $r$-tuple in $\Xim{n}(\kkdih{n})$. With this, we have dealt with the first half of \Cref{tab:data}. The analogous statements for the second half of the table and the tuple $\Re{\wt{\mbP}}{n}$ then follow by combining the identities \eqref{eq:bbfbul-u-v2} with the following identities (see \eqref{eq:uim} and \eqref{eq:re-ol-mbu-2}):
    \begin{equation*}
        (\mbure{2n})^2 = \mbure{n}+2, \quad (\Re{\ol{\mbu}}{2n})^2 = \Re{\ol{\mbu}}{n}+2, \quad \textup{ and } \quad (\mbuim{2n})^2 = \mbure{n}-2
    \end{equation*}
    This proves that the formulas in \Cref{tab:data} and \eqref{eq:Re-mbP-points} define $r$-tuples of points in the respective curves, as claimed. The final claim in part (c) follows now from \Cref{prop:Ln-fact} together with  the fact that the isomorphism $\kkdih{2n} \iso \kkre{2n}$ sends each $\ure{i,2n}$ to $\Re{\ol{u}}{i,2n}$ and each $\ure{i,n} = (\ure{i,2n})^2 - 2$ to $\Re{\ol{u}}{i,n} = (\Re{\ol{u}}{i,2n})^2 - 2$ (see \ref{emp:direct-system}).
\end{proof}

The automorphisms from  \ref{emp:aut-kkbul} act naturally on these $r$-tuples. In the next lemma we count the total number of Galois conjugates arising from these actions (the proof is a straightforward verification which is left to the reader). 

\begin{prop}\label{prop:points-conjugates}
Fix $n \in \N$ and the label $\bullet$.  Let $\kk_n$ denote the corresponding rational function field in \Cref{tab:data}. Recall the automorphisms of $\kk_n$ defined in \ref{emp:aut-kkbul}. Then: 
\begin{enumerate}[(a),itemindent=1cm, leftmargin=0pt]
\item If $\bullet \in \{\textup{\textsf{cyc}, \textsf{re}, \textsf{im}}\}$, and $\zeta_n \in k$ (resp. $\zeta_{2n} \in k$),  then the points 
\begin{equation*}
\begin{alignedat}{2}
    \iota^{e}\sigma_{i,n}^j (\Qbul{i,n}) & \in \Xbul{n}(\kk_n), \qquad &&  i=1,\dotsc,r, j=1,\dotsc,n, e=0,1,\\
    \sigma_{i,2n}^j (\wtQbul{i,n}) & \in \Xwtbul{n}(\kk_{2n}), \qquad &&  i=1,\dotsc,r, j=1,\dotsc,2n,
\end{alignedat}
\end{equation*}
are pair-wise distinct and hence contribute $2rn$ points to the respective curves. 
\item The curves $\Xdih{n}/K$ and $\Xwtdih{n}/K$ admit the extra involution $\delta: (x,y) \mapsto (x\inv, y/(x^{g+1}))$ (where $g$ denotes the genus $\gdih{n}$ or $\gwtdih{n}$ of the curves, respectively), and the points
\begin{equation*}
\begin{alignedat}{2}
    \iota^{e}\tau^\ell\sigma_{i,n}^j (\Qdih{i,n}) & \in \Xdih{n}(\kkdih{n}), \qquad &&  i=1,\dotsc,r,\, j=1,\dotsc,n, \;\;\ell=0,1,\, e=0,1, \\
    \tau^\ell\sigma_{i,2n}^j (\wtQdih{i,n}) & \in \Xwtdih{n}(\kkdih{2n}), \qquad && i=1,\dotsc,r,\, j=1,\dotsc,2n, \ell=0,1,
\end{alignedat}
\end{equation*}
are pair-wise distinct and hence contribute $4rn$ points to the respective curves.
\item If $n\geqslant 3$ and $k$ contains $\lambda_n$, then the points 
\begin{equation*}
	\iota^{e}\ol{\sigma}_{i,n}^j (\Re{P}{i,n})  \in \Xre{n}(\kkre{n}), \qquad   i=1,\dotsc,r, j=1,\dotsc,n, e=0,1,
\end{equation*}
are distinct and contribute $2rn$ points to $\Xre{n}$. Similarly, if $n\geqslant 2$ and $k$ contains $\lambda_{2n}$, then the points
\begin{equation*}
	\ol{\sigma}_{i,2n}^j (\Re{P}{i,n})  \in \Xwtre{n}(\kkre{2n}), \qquad  i=1,\dotsc,r, j=1,\dotsc,2n,
\end{equation*}
are distinct and contribute $2rn$ points to $\Xwtre{n}$. 
\end{enumerate}
\end{prop}

\begin{emp}\label{emp:Wbul}
Fix $n\in \N$ and $\bullet \in \{\textup{\textsf{cyc}, \textsf{re}, \textsf{im}, \textsf{dih}}\}$. Assume that $\zeta_n \in k$, so that we may consider the tuple of points $\mbQbul{n} \in \Xbul{n}(\kk_n)$. Viewing $\Xbul{n}(\kk_n)$ as a module over the Galois group 
\begin{equation*}
    \bb{G}_n \ceq \gal(\kk_n/\kk),
\end{equation*}
we define the degree zero divisors $2\Qbul{i,n} - D_{\infty}$ (cf. \ref{emp:hyperelliptic}) and we define
\begin{equation*}
    \Wbul{i,n} \ceq \spn_{\Z \bb{G}_n} \{ \bul{\mf{Q}}{i,n} \}, \qquad \mbWbul{n} \ceq  \Wbul{1,n} \oplus \dotsb \oplus \Wbul{r,n} \subset \Div^0(\Xbul{n})_{\kk_n}.
\end{equation*}
Similarly, if $\zeta_{2n} \in k$ then we have the divisors $\bul{\wt{\mf{Q}}}{i,n} \ceq 2\wtQbul{i,n} - D_{\infty}$  and we can define
\begin{equation*}
    \wtWbul{i,n} \ceq \spn_{\Q \bb{G}_{2n}} \{ \bul{\wt{\mf{Q}}}{i,n} \}, \qquad \wtmbWbul{n} \ceq  \wtWbul{1,n} \oplus \dotsb \oplus \wtWbul{r,n} \subset \Div^0(\Xwtbul{n})_{\kk_{2n}}.
\end{equation*}
For each of these spaces $W$, we denote by $[W]$ the image of $W$ under the natural map from $\Div^0$ to $\pic^0$ which sends a degree zero $\Q$-divisor $E$ to its class $[E]$. Thus, the integers
\begin{equation}
\bul{R}{n} \ceq \dim_{\Q} [\mbWbul{n}], \qquad \bul{\wt{R}}{n}  \ceq \dim_{\Q} [\wtmbWbul{n}]
\end{equation}
serve as lower bounds for the rank of $\Xbul{n}/\kk_n$ (resp. $\Xwtbul{n}/\kk_{2n}$).

Similarly,  for $\bullet \in \{\textup{\textsf{cyc}, \textsf{re}, \textsf{im}}\}$, we define the integers 
\begin{equation}
\bul{N}{n} \ceq 2rn + \bul{b}{n}, \qquad \bul{\wt{N}}{n} \ceq 2rn + \bul{\wt{b}}{n},
\end{equation}
and for the dihedral curves, we define
\begin{equation}
\dih{N}{n} \ceq 4rn + \dih{b}{n}, \qquad \dih{\wt{N}}{n} \ceq 4rn + \dih{\wt{b}}{n}.
\end{equation}
\end{emp}

\begin{lem}\label{lem:lower-bd}
For $n\in \N$ and each $\bullet$:
\begin{enumerate}[(a)]
    \item We have the lower bounds
    \begin{equation}
    \begin{alignedat}{6}
        \# \Xbul{n}(\kk_n) & \geqslant \bul{N}{n} && \quad \textup{ and } \quad && \rk \Xbul{n}/\kk_n && \geqslant \bul{R}{n}, \quad && \textup{ if } \zeta_n && \in k,\\
        \# \Xwtbul{n}(\kk_{2n}) & \geqslant \bul{\wt{N}}{n} && \quad \textup{ and } \quad && \rk \Xwtbul{n}/\kk_{2n} && \geqslant \bul{\wt{R}}{n}, \quad && \textup{ if } \zeta_{2n} && \in k.
    \end{alignedat}
    \end{equation}
    \item If $n\geqslant 3$, then we have the lower bounds
    \begin{equation}
    \begin{alignedat}{6}
        \# \Xre{n}(\kkre{n}) & \geqslant \bul{N}{n} && \quad \textup{ and } \quad && \rk \Xre{n}/\kkre{n} && \geqslant \bul{R}{n}, \quad && \textup{ if } \lambda_n && \in k,\\
        \# \Xwtre{n}(\kkre{2n}) & \geqslant \bul{\wt{N}}{n} && \quad \textup{ and } \quad && \rk \Xwtre{n}/\kkre{2n} && \geqslant \bul{\wt{R}}{n}, \quad && \textup{ if } \lambda_{2n} && \in k.
    \end{alignedat}
    \end{equation}
\end{enumerate}
\end{lem}

\begin{proof}
For (a), the lower bound on the number of points follows by counting the points coming from \Cref{prop:points-conjugates} together with the $\bul{b}{n}$ (or $\bul{\wt{b}}{n}$) Weierstrass points, and the statement about the rank is obvious. Part (b) is an immediate consequence of \Cref{prop:table-data} (c) and the fact that the isomorphism $\phi: \kkdih{n} \iso \kkre{n}$ from \Cref{prop:Ln-fact} induces a bijection from the $\gal(\kkdih{n}/\kk)$-conjugates of $\mbQre{n} \in \Xre{n}(\kkdih{n})$ to the $\gal(\kkre{n}/\kk)$-conjugates of $\Re{\mbP}{n}$. The proof for $\Xwtre{n}$ is similar.
\end{proof}


Now, the bulk of the remainder of the paper is devoted to proving the following.

\begin{thm}[Proved in \ref{emp:proof-pt-rank-bounds}]\label{thm:point-rank-bounds}
The integers $\bul{R}{n}$ and $\bul{\wt{R}}{n}$ satisfy the values in the following table. 
\begin{table}[H]
\centering
\renewcommand{\arraystretch}{1.5}
\setlength{\tabcolsep}{5pt}
\adjustbox{max width=\textwidth}{
\caption{Rank bounds}\label{tab:ranks} 
\begin{tabular}{|c||c|c||c|c||c|c|} 
    \hline
    &  \multicolumn{2}{c||}{$d=1$} & \multicolumn{2}{c||}{$d=2$} & \multicolumn{2}{c|}{$d\geqslant 3$} \\
    \hhline{|-||--||--||--|}
    $\bullet$ & $\bul{R}{n}$ & $\bul{\wt{R}}{n}$ & $\bul{R}{n}$ & $\bul{\wt{R}}{n}$ & $\bul{R}{n}$ & $\bul{\wt{R}}{n}$ \\
    \hhline{|=||=|=||=|=||=|=|}
    \textup{cyc} &  {$8\gcyc{n}$} & {$8\gwtcyc{n}$} & {$6\gcyc{n}$} & {$6\gwtcyc{n}$} & {$rn-1$} & {$rn$}\\
    \hhline{|-||-|-||-|-||-|-|}
    \textup{re} &  {$8\gre{n}$} & {$8\gwtre{n}$} & {$6\gre{n}$} & {$6\gwtre{n}$} & $rn-1$ & $rn$\\
    \hhline{|-||-|-||-|-||-|-|}
    \textup{im} &  {$8\gim{n}-4$} & {$8\gwtim{n}$} & {$6\gim{n}$} & {$6\gwtim{n}$} & $rn$ & $rn$\\
    \hhline{|-||-|-||-|-||-|-|}
    \textup{dih} &  {$8\gdih{n}-4$} & {$8\gwtdih{n}$} & {$6\gdih{n}$} & {$6\gwtdih{n}$} & $2rn-1$ & $2rn$ \\ 
    \hline
\end{tabular}}
\end{table}
\end{thm}

\begin{emp}\label{emp:rank-strategy}
\textbf{Proof strategy.}
We will proceed in three steps:
\begin{enumerate}[(i), leftmargin=0pt,itemindent=1cm]
    \item First, we show in \Cref{lem:reduction-1} that \Cref{thm:point-rank-bounds}  holds in full generality if it holds for the integers $\cyc{R}{n}, \Re{R}{n},$ and $\Im{R}{n}$. This follows by exploiting the fact that the curves $\Xbul{n}/K$ and $\Xwtbul{n}/K$ are related by various covers which induce various isogeny relations between the respective Jacobians.
    \item Next, in \Cref{prop:Wn-Vn-decomp} we reduce the task of determining the integers $\cyc{R}{n}, \Re{R}{n}$ and $\Im{R}{n}$ to determining, in each case (and for all $n\in \N$), the dimension of a certain distinguished subspace $[\bul{\mbV}{n}] \subset [\mbWbul{n}]$. This again follows by exploiting the fact that each family $\{ \Xbul{n} \}_{n \in \N}$ forms an inverse system with very nice properties (see \Cref{sec:zhat}).
    \item Finally, in \Cref{thm:dim-Vn} we apply the rank-bounding tools from \Cref{sec:tools} to determine the dimension of each space $[\bul{\mbV}{n}]$, and we complete the proof of \Cref{thm:point-rank-bounds}.
\end{enumerate}
\end{emp}

We conclude this section by showing how \Cref{thm:main} follows from \Cref{thm:point-rank-bounds}, using the fact (proved in \Cref{prop:non-iso}) that any positive genus curve in the real, imaginary, and dihedral families is non-isotrivial over its field of definition $K$.

\begin{emp}\label{emp:proof-thm-main}
\textit{Proof of \Cref{thm:main}}. There are five cases to consider, depending on the value of $dn$ (in terms of $g$); we use below the fact that $2rn = 4dn + 4n$, and the fact that for odd $n$ we have $\Q(\lambda_n) = \Q(\lambda_{2n})$ and $\Q(\zeta_n) = \Q(\zeta_{2n})$.
\begin{enumerate}[$\bullet$,leftmargin=10pt]
    \item \uline{$2g+2$}: If $nd$ is even then $nd=2\gre{n}+2$. If $\lambda_n \in k$ then \Cref{thm:point-rank-bounds} gives $$\# \Xre{n}(\kkre{n}) \geqslant 2rn+\Re{b}{n} = 8 \gre{n} + 4n + 8.$$ 
    \item \uline{$2g+1$}: If $nd$ is odd then $nd = 2\gwtre{n} + 1$. If $\lambda_n \in k$, then also $\lambda_{2n} \in k$, so \Cref{thm:point-rank-bounds} gives $$\# \Xwtre{n}(\kkre{2n}) \geqslant 2rn+\Re{\wt{b}}{n} = 8 \gre{n} + 4n + 5.$$ 
    \item \uline{$2g$}: If $nd$ is even then $nd = 2\gwtre{n}$ and $\Re{\wt{b}}{n} = 2$, which implies that $\Jwtre{n}$ has a $K$-rational point of order $2$. If $\lambda_{2n} \in k$ then \Cref{thm:point-rank-bounds} gives $$\# \Xwtre{n}(\kkre{2n}) \geqslant 2rn+\Re{\wt{b}}{n} = 8 \gre{n} + 4n + 2.$$
    \item \uline{$2g-1$}: If $nd$ is odd then $nd = 2\gim{n} - 1$ and $\Im{b}{n} = 3$, so $\Jim{n}$ has three $K$-points of order $2$ which generate the group $V_4$ inside $\Jim{n}(K)$. If $\zeta_n \in k$ then \Cref{thm:point-rank-bounds} gives $$\# \Xim{n}(\kkdih{n}) \geqslant 2rn+ \Im{b}{n} = 8 \gim{n} + 4n -1.$$ 
    \item \uline{$g+1$}: If $nd$ is even then $nd = \gdih{n}+1$, so if $\zeta_n \in k$ then \Cref{thm:point-rank-bounds} gives $$\# \Xdih{n}(\kkdih{n}) \geqslant 4rn+\dih{b}{n} = 8 \gdih{n} + 8n+8.$$ If $nd$ is odd then $nd = \gwtdih{n}+1$, and if $\zeta_n \in k$, then also $\zeta_{2n} \in k$; so, \Cref{thm:point-rank-bounds} gives $$\# \Xwtdih{n}(\kkdih{2n}) \geqslant 4rn+\dih{\wt{b}}{n} = 8 \gwtdih{n} + 8n+8.$$
\end{enumerate}
One checks, in each of the five cases above, that the rank bounds in \Cref{tab:ranks} coincide with the corresponding bounds given in \Cref{thm:main}. In all cases considered, the points are defined over a rational function field over $k$, and  the curves in question are non-isotrivial over $K = k(\mbf{h})$ (this is proved in \Cref{prop:non-iso}). The desired result now follows by specialization (more precisely, \Cref{cor:family-spec}).
\end{emp}


\section{First reduction: partner curves}\label{sec:partner}



\begin{lem}\label{lem:hyp-partner}
Let $X/K$ be a hyperelliptic curve which admits double covers
\begin{equation*}
    \begin{tikzcd}
        & X \arrow[dl, "\phi"'] \arrow[dr, "\wt{\phi}"] \\
        Y & & \wt{Y}.
    \end{tikzcd}
\end{equation*}
Let $\tau$ and $\wt{\tau} \in \aut_K(X)$ be the involutions with quotients $\phi$ and $\wt{\phi}$, respectively. Assume that $\wt{\tau} = \iota \tau \in \aut_K(X)$ (cf. \ref{emp:hyperelliptic}). Then, the homomorphism
\begin{equation*}\label{eq:hyp-jac-pushf}
    (\phi_*,{\wt{\phi}}_*): J_{X} \too J_{Y} \times J_{\tilde{Y}}
\end{equation*}
is a $K$-isogeny, with dual isogeny defined by 
\begin{equation*}\label{eq:hyp-jac-pullb}
    \phi^* + \wt{\phi}^* : J_{Y} \times J_{\tilde{Y}} \too J_{X}.
\end{equation*}
\end{lem}
Following Shioda~\cite{shiodasymmetry}, we call the curves $Y$ and $\wt{Y}$ above \defi{partner curves} relative to $X$. 
\begin{proof}
The key point is that $\iota P = -P$ for any geometric point $P\in J_X(\ol{K})$ (because if $P$ is represented by a divisor class $[D]$, then $P+\iota P = \pi^*[\pi_*D] = 0 \in J_X(\ol{K})$). So, we have $\wt{\tau}P = \iota \tau P = - \tau P \in J_X(\ol{K})$, giving
\begin{equation*}
    (\phi^* + \wt{\phi}^*)\comp (\phi_*,\wt{\phi}_*)(P) = \phi^*\phi_*(P) + \wt{\phi}^*\wt{\phi}_*(P) = (1+\tau)P + (1 - {\tau})P = 2P,
\end{equation*}
as desired. 
\end{proof}

\begin{emp}\label{emp:zhat-formula}
\textbf{Coordinates and double covers.} Let $n\in \N$. We label the coordinates $(x,y)$ on each of the curves $\Xbul{n}$ and $\Xwtbul{n}$ from \eqref{eq:Xbul-Xwtbul} as follows. For $\bullet \in \{ \textup{\textsf{cyc}, \textsf{re}, \textsf{dih}} \}$, we write $(\xbul{n},\ybul{n})$ and $(\xbul{n},\ywtbul{n})$ for the coordinates on $\Xbul{n}$ and $\Xwtbul{n}$, respectively. On the other hand, we write $(\xRe{n},\yim{n})$ and $(\xRe{n},\ywtim{n})$ for the coordinates on $\Xim{n}$ and $\Xwtim{n}$, respectively. In the next table, we define various formulas relating the coordinates which define double covers between the respective curves; below, $\trm(x)$ denotes the rational function $x - x\inv$. We note also that $\fRe{2}(x) = x^2 - 2 \in \Z[x]$. 
\begin{table}[H]
\centering
\renewcommand{\arraystretch}{2.0}
\setlength{\tabcolsep}{5pt}
\adjustbox{max width=\textwidth}{
\caption{Formulas defining double covers}\label{tab:formulas}
\begin{tabular}{|l|ll||l|ll|} 
        \hline
        $\phicyc{n}$ & $(\xcyc{n},\ycyc{n})$ & $\ceq \left((\xcyc{2n})^{2},\ycyc{2n} \right)$ & $\wtphicyc{n}$ & $(\xcyc{n},\ywtcyc{n})$ & $\ceq \left((\xcyc{2n})^{2},\ycyc{2n}\xcyc{2n} \right)$\\ \hline
        $\phire{n}$ & $(\xRe{n},\yRe{n})$ & $\ceq \left(\fRe{2}(\xRe{2n}),\yRe{2n} \right)$ & $\wtphire{n}$ & $(\xRe{n},\ywtRe{n})$ & $\ceq \left(\fRe{2}(\xRe{2n}),\yRe{2n}\xRe{2n} \right)$\\ \hline
        $\phiim{n}$ & $(\xRe{n},\yim{n})$ & $\ceq \left(\fRe{2}(\xRe{2n}),\yim{2n} \xRe{2n}\right)$ & $\wtphiim{n}$ & $(\xRe{n},\ywtim{n})$ & $\ceq \left(\fRe{2}(\xRe{2n}),\yim{2n}\right)$\\ \hline
        $\phidih{n}$ & $(\xdih{n},\ydih{n})$ & $\ceq \left((\xdih{2n})^{2}, \ydih{2n}(\xdih{2n})^{\left[ nd \right]}\right)$ & $\wtphidih{n}$ & $(\xdih{n},\ywtdih{n})$ & $\ceq \left((\xdih{2n})^{2}, \ydih{2n} (\xdih{2n})^{1-\left[ nd \right]}\right)$\\ \hhline{|======|}
        $\Re{\psi}{n}$ & $(\xRe{n},\yRe{n})$ & $\ceq \left(\tr(\xdih{n}),\dfrac{\ydih{n}}{(\xdih{n})^{\frac{nd+ \left[ nd \right]}{2}}} \right)$ & $\Re{\wt{\psi}}{n}$ & $(\xRe{n},\ywtRe{n})$ & $\ceq \left(\tr(\xdih{n}),\dfrac{\ywtdih{n}\tr(\xdih{n})}{(\xdih{n})^{\frac{nd -\left[ nd \right]}{2}}} \right)$\\ \hline
        $\Im{\psi}{n}$ & $(\xRe{n},\yim{n})$ & $\ceq \left(\tr(\xdih{n}),\dfrac{\ydih{n}\trm(\xdih{n})}{(\xdih{n})^{\frac{nd+\left[ nd \right]}{2}}}\right)$  & $\Im{\wt{\psi}}{n}$ & $(\xRe{n},\ywtim{n})$ & $\ceq \left(\tr(\xdih{n}),\dfrac{\ywtdih{n}\trm(\xdih{n})}{(\xdih{n})^{\frac{nd -\left[ nd \right]}{2}}} \right)$\\\hline
\end{tabular}}
\end{table}
All of the double covers defined above fit into the following pair of commutative diagrams of curves over $K$ (as usual, each morphism $\pi$ is the natural projection-to-$x$ double cover):
\begin{equation}\label{diag:best-diag}
    \begin{tikzcd}
        & & & \Xdih{2n} \arrow[ddll, "\Re{\psi}{2n}"'] \arrow[ddrr, "\Im{\psi}{2n}"] \arrow[dd, "\phidih{n}\;\;"',shift left] \arrow[dd, "\,\;\;\dih{\tilde{\phi}}{n}", shift right] \\ \\
        \Xcyc{2n} \arrow[dd, "\phicyc{n}\;\;"',shift left] \arrow[dd, "\,\;\;\wtphicyc{n}", shift right] & \Xre{2n} \arrow[dd, "\phidih{n}\;\;"',shift left] \arrow[dd, "\,\;\;\wtphire{n}", shift right] & & (\Xdih{n},\Xwtdih{n}) \arrow[ddll, "\Re{\psi}{n}\;\;\;"',shift left] \arrow[ddll, "\;\;\;\Re{\wt{\psi}}{n}",shift right] \arrow[ddrr, "\Im{\wt{\psi}}{n}", shift left] \arrow[ddrr, "\Im{\psi}{n}"', shift right]   \arrow[dd, "\pi\;\;"',shift left] \arrow[dd, "\,\;\;\pi", shift right] & & \Xim{2n} \arrow[dd, "\phiim{n}\;\;"',shift left] \arrow[dd, "\,\;\;\wtphiim{n}", shift right] \\ \\
        (\Xcyc{n},\Xwtcyc{n}) \arrow[dd, "\pi\;\;"',shift left] \arrow[dd, "\,\;\;\pi", shift right] & (\Xre{n},\Xwtre{n}) \arrow[ddrr, "\pi\;\;"',shift left] \arrow[ddrr, "\,\;\;\pi", shift right] & & \p^1_K(\xdih{n}) \arrow[dd, "\tr"'] & & (\Xim{n},\Xwtim{n}) \arrow[ddll, "\pi\;\;"',shift left] \arrow[ddll, "\,\;\;\pi", shift right] \\ \\ 
        \p^1_K(\xcyc{n}) & & & \p^1_K(\xRe{n}). 
    \end{tikzcd}
\end{equation}
The next proposition records that fact that for each pair of double covers in a row of \Cref{tab:formulas}, the two targets are partner curves relative to the source (cf. \Cref{lem:hyp-partner}). We write $\Jbul{n}/K$ and $\Jwtbul{n}/K$ for the Jacobians of $\Xbul{n}/K$ and $\Xwtbul{n}/K$, respectively. 
\end{emp}

\begin{prop}\label{prop:partner-isog}
For $n \in \N$, the following are $K$-isogenies:
\begin{equation}\label{eq:isog-partner-curves}
    \begin{alignedat}{1}
        ((\Re{\psi}{n})_*,(\Im{\psi}{n})_*): \Jdih{n} & \too \Jre{n} \times \Jim{n},\\
        ((\Re{\wt{\psi}}{n})_*,(\Im{\wt{\psi}}{n})_*): \Jwtdih{n} & \too \Jwtre{n} \times \Jwtim{n},\\
        ((\phibul{n})_*,(\wtphibul{n})_*): \Jbul{2n} & \too \Jbul{n} \times \Jwtbul{n},\qquad  \textup{ for } \bullet \in \{ \textup{\textsf{cyc}, \textsf{re}, \textsf{im}, \textsf{dih}} \}.
        \end{alignedat}
    \end{equation}
\end{prop}

\begin{proof}
For $\bullet = \textsf{re}$ and $\textsf{cyc}$, the double covers $\phibul{n}$ and $\wtphibul{n}$ are respectively quotients by the involutions 
\begin{equation*}
    \tau_n: (\xbul{2n},\ybul{2n}) \mapsto (-\xbul{2n},\ybul{2n}), \qquad {\wt{\tau}}_{n}:(\xbul{2n},\ybul{2n}) \mapsto (-\xbul{2n},-\ybul{2n}).
\end{equation*}
The covers $\phiim{n}$ and $\wtphiim{n}$, on the other hand, are respectively quotients by
\begin{equation*}
    \tau_n:(\xRe{2n},\yim{2n}) \mapsto (-\xRe{2n},-\yim{2n}), \qquad \wt{\tau}_n:(\xRe{2n},\yim{2n}) \mapsto (-\xRe{2n},\yim{2n}),
\end{equation*}
and $\phidih{n}$ and $\wtphidih{n}$ are respectively quotients by 
\begin{equation*}
    \tau_n:(\xdih{2n},\ydih{2n}) \mapsto (-\xdih{2n},(-1)^{\left[ nd \right]}\ydih{2n}), \qquad {\wt{\tau}}_{n}:(\xdih{2n},\ydih{2n}) \mapsto (-\xdih{2n},(-1)^{1-\left[ nd \right]}\ydih{2n}).
\end{equation*}
In all four of these cases, we have the relation ${\wt{\tau}}_{n} = \iota \tau_n \in \aut_K(\Xbul{2n})$, so \Cref{lem:hyp-partner} implies that $((\phibul{n})_*,(\wtphibul{n})_*):\Jbul{2n} \too \Jbul{n}\times \Jwtbul{n}$ is a $K$-isogeny. 

Next, $\Re{\psi}{n}$ and $\Im{\psi}{n}$ are quotients, respectively, by the involutions 
\begin{equation*}
    \delta_n:(\xdih{n},\ydih{n}) \mapsto \left(\frac{1}{\xdih{n}},\dfrac{\ydih{n}}{(\xdih{n})^{nd+\left[nd\right]}}\right), \qquad \gamma_n:(\xdih{n},\ydih{n}) \mapsto \left(\frac{1}{\xdih{n}},\dfrac{-\ydih{n}}{(\xdih{n})^{nd+\left[nd\right]}}\right).
\end{equation*}
Similarly, $\Re{\wt{\psi}}{n}$ and $\Im{\wt{\psi}}{n}$ are quotients, respectively, by the involutions 
\begin{equation*}
    {\wt{\delta}}_{n}:(\xdih{n},\ywtdih{n}) \mapsto \left( \frac{1}{\xdih{n}},\dfrac{\ywtdih{n}}{(\xdih{n})^{nd-\left[nd\right]}} \right),\qquad {\wt{\gamma}}_{n}:(\xdih{n},\ywtdih{n}) \mapsto \left(\frac{1}{\xdih{n}},\dfrac{-\ywtdih{n}}{(\xdih{n})^{nd-\left[nd\right]}} \right).
\end{equation*}
That the first two homomorphisms in \eqref{eq:isog-partner-curves} are $K$-isogenies now follows from \Cref{lem:hyp-partner} and the fact that $\gamma_{n} =\iota {\delta}_{n} \in \aut_K(\Xdih{n})$ and ${\wt{\gamma}}_{n} =\iota {\wt{\delta}}_{n} \in \aut_K(\Xwtdih{n})$. This concludes the proof. 
\end{proof}

Following the strategy in \ref{emp:rank-strategy}, we reduce the proof of \Cref{thm:point-rank-bounds} as follows.

\begin{lem}\label{lem:reduction-1}
For all $n \in \N$, we have $\dih{R}{n} = \Re{R}{n} + \Im{R}{n},$
and for all four values of $\bullet$, we have $\bul{\wt{R}}{n} = \bul{R}{2n} - \bul{R}{n}.$ Consequently, if the values for $\cyc{R}{n},\Re{R}{n}$ and $\Im{R}{n}$ in \Cref{tab:ranks} are correct for all $n\in \N$, then \Cref{thm:point-rank-bounds} holds.
\end{lem}

\begin{proof}
First note that $\Re{\psi}{n}$ and $\Im{\psi}{n}$ map the $r$-tuple $\mbQdih{n} \in \Xdih{n}(\kk_n)$ to the $r$-tuples $\mbQre{n} \in \Xre{n}(\kk_n)$ and $\mbQim{n} \in \Xim{n}(\kk_n)$, respectively. Thus, they give rise to surjections $\Re{\psi}{n}: \mbWdih{n} \surj \mbWre{n}$ and $\Im{\psi}{n}: \mbWdih{n} \surj \mbWim{n}$. The first isogeny from \eqref{eq:isog-partner-curves} then induces an isomorphism of $\Q$-vector spaces
\begin{equation*}
    ((\Re{\psi}{n})_*,(\Im{\psi}{n})_*): \mbWdih{n} \iso \mbWre{n} \oplus \mbWim{n},
\end{equation*}
which proves that $\dih{R}{n} = \Re{R}{n} + \Im{R}{n}$. An identical argument serves to show that for each $\bullet$ the third isogeny from \eqref{eq:isog-partner-curves} induces an isomorphism of $\Q$-vector spaces
\begin{equation*}
    ((\phibul{n})_*,(\wtphibul{n})_*): \mbWbul{2n} \iso \mbWbul{n} \oplus \wtmbWbul{n},
\end{equation*}
which proves that $\bul{\wt{R}}{n} = \bul{R}{2n} - \bul{R}{n}$. In particular, we recover $\dih{\wt{R}}{n}$ as
$\dih{\wt{R}}{n} = (\Re{R}{2n} + \Im{R}{2n}) - (\Re{R}{n} + \Im{R}{n})$. 
The desired result now follows by checking that these identities are satisfied by the values in \Cref{tab:ranks}; we leave this calculation to the reader. 
\end{proof}

\section{\texorpdfstring{$\zhat$}{}-systems and new parts decompositions}\label{sec:zhat}

\begin{emp}\label{emp:zhat-systems}
\textbf{$\zhat$-systems.} 
View the set $\N$ of positive integers as a poset via the ``inverted'' divisibility relation on $\N$ in which $n \leqslant m$ if and only if $m$ divides $n$ (denoted $m \mid n$). We write $m \| n$ if $m \mid n$ and $m<n$. By a \defi{$\zhat$-system} (of nice curves) over $K$, we mean an inverse system $((\X_n/K)_{n\in\N},(\phi_{n,m})_{m \mid n \in \N})$ (abbreviated by writing $(\X_n,\phi_{n,m})$) of nice curves over $K$ in which each morphism $\phi_{n,m}:\X_n \too \X_m$ is a cover which of degree $n/m$.

The projective lines $\p^1_K(\xbul{n})$ in \eqref{diag:best-diag} give rise to three $\zhat$-systems $(\p^1_K(\bul{x}{n}),\bul{F}{n,m})$ in which each cover is defined, for $n=me$, as follows:
\begin{equation*}
\begin{alignedat}{4}
	& \cyc{F}{n,m} : \p^1_K(\xcyc{n}) && \too \p^1_K(\xcyc{m}), \qquad && \xcyc{m} && \ceq (\xcyc{n})^{e},\\
        & \Re{F}{n,m} : \p^1_K(\xRe{n}) && \too \p^1_K(\xRe{m}), \qquad && \xRe{m} && \ceq \fRe{e}(\xRe{n}),\\
        & \dih{F}{n,m} : \p^1_K(\xdih{n}) && \too \p^1_K(\xdih{m}), \qquad && \xdih{m} && \ceq (\xdih{n})^{e}.
\end{alignedat}
\end{equation*}

The curves $\Xbul{n}$ are all double covers of these projective lines, and they similarly give rise to four $\zhat$-systems of hyperelliptic curves $(\Xbul{n},\phibul{n,m})$, for $\bullet \in \{ \textup{\textsf{cyc}, \textsf{re}, \textsf{im}, \textsf{dih}} \}$; for $n=me \in \N$, the covers $\phibul{n,m}$ are defined as follows:
\begin{equation}
    \begin{alignedat}{4}
        & \phicyc{n,m} : \Xcyc{n} && \too \Xcyc{m}, \qquad && (\xcyc{m},\ycyc{m}) && \ceq ((\xcyc{n})^{e},\ycyc{n}),\\
        & \phire{n,m} : \Xre{n} && \too \Xre{m}, \qquad && (\xRe{m},\yRe{m}) && \ceq (\fRe{e}(\xRe{n}),\yRe{n}),\\
        & \phiim{n,m} : \Xim{n} && \too \Xim{m}, \qquad && (\xRe{m},\yim{m}) && \ceq (\fRe{e}(\xRe{n}),\yim{n} \fIm{e}(\xRe{n})),\\
        & \phidih{n,m} : \Xdih{n} && \too \Xdih{m}, \qquad && (\xdih{m},\ydih{m}) && \ceq ((\xdih{n})^e,\ydih{n}(\xdih{n})^{\frac{e \left[nd \right] - \left[ md\right]}{2}}).
    \end{alignedat}
\end{equation}
Note in particular that each cover $\phibul{2n,n}$ is nothing but the double cover $\phibul{n}$ from \eqref{diag:best-diag}.
\end{emp}

\begin{emp}\label{emp:new-parts}
\textbf{The new parts $\mathbold{\bb{A}_n}/K$.} 
Let $(\X_n,\phi_{n,m})$ be a $\zhat$-system of curves over $K$, and for each $n\in \N$, let $\J_n/K$ denote the Jacobian of $\X_n/K$. Then each $\phi_{n,m}$ gives rise to two homomorphisms of abelian varieties over $K$, namely 
\begin{equation*}
\begin{alignedat}{3}
    (\phi_{n,m})_*:\J_n &\too \J_m,  && \qquad [D] &&\mtoo [(\phi_{n,m})_*D],\\
    (\phi_{n,m})^*:\J_m &\too \J_n,  && \qquad [E] &&\mtoo [(\phi_{n,m})^*E].
\end{alignedat}
\end{equation*}
For each $n \geqslant 2$, we define the following  abelian sub-variety of $\J_n/K$:
\begin{equation}
    \bb{A}_n/K \ceq \textup{The identity component of } \bigcap_{m \| n} \ker \left( \J_n \xlongrightarrow{(\phi_{n,m})_*} \J_m \right).
\end{equation}
\end{emp}

\begin{lem}\label{lem:fib-prod-deg}
Let $d$ and $m$ be positive integers, neither of which divides the other, and put $\ell \ceq \lcm(d,m)$. Then, we have
\begin{equation*}
    (\phi_{\ell,m})_* \comp (\phi_{\ell,d})^* = 0 \in \Hom(\bb{A}_d,\J_m). 
\end{equation*}
\end{lem}

\begin{proof}
Put $h \ceq \gcd(d,m)$. The desired result follows from the claim that
\begin{equation}\label{eq:psi-c-b-a}
    (\phi_{\ell,m})_*(\phi_{\ell,d})^* = (\phi_{m,h})^*(\phi_{d,h})_* \in \Hom(\J_d,\J_m).
\end{equation}
Let $B \subset \X_h$ denote the set of branch points of $\phi_{d,h} \comp \phi_{\ell,d}$. Put $Y_{h} \ceq  \X_{h} \setminus B$, and let $Y_d,Y_m,$ and $Y_{\ell}$ denote the open subsets of $X_{d},X_{m},$ and $X_{\ell}$ obtained by removing the pre-image of $B$ under the respective morphism to $\X_{h}$. Any geometric point of $\J_{d}$ can be represented by the class of some divisor $D \in \Div^0(Y_{d})_{\ol{K}}$, so to prove \eqref{eq:psi-c-b-a}, it suffices to show that 
\begin{equation*}
    (\phi_{\ell,m})_*(\phi_{\ell,d})^* = (\phi_{m,h})^*(\phi_{d,h})_* \in \Hom(\Div Y_{d}, \Div Y_{m}).
\end{equation*}
For this, it suffices to show that for any geometric point $P \in Y_{d}(\ol{K})$, we have  
\begin{equation*}
    (\phi_{\ell,m})_*(\phi_{\ell,d})^*(P) = (\phi_{m,h})^*(\phi_{d,h})_*(P).
\end{equation*}
This in turn would follow from the claim that the following diagram of affine curves is \emph{Cartesian}:
\begin{equation*}
    \begin{tikzcd}
        & Y_{\ell} \arrow[dl, "\phi_{\ell,d}",swap] \arrow[dr, "\phi_{\ell,m}"] &\\
        Y_{d} \arrow[dr, "\phi_{d,h}",swap] & & Y_{m} \arrow[dl, "\phi_{m,h}"]\\
        & Y_{h}. &
    \end{tikzcd}
\end{equation*}
That is, we want to show that the natural induced morphism of $Y_{h}$-schemes 
\begin{equation}\label{eq:Yc-fib}
    Y_{\ell} \too Y_{d} \times_{Y_h} Y_m
\end{equation}
is an isomorphism. To see this, let $L_i$ denote the function field of each $Y_i$, and note that the restriction of \eqref{eq:Yc-fib} to the generic point of $Y_{h}$ gives rise to a homomorphism of $L_{h}$-algebras
\begin{equation*}
    L_{d} \otimes_{L_{h}} L_{m} \too L_{\ell}.
\end{equation*}
This must be an isomorphism because the existence of a non-trivial kernel would imply the existence of an intermediate field $L_{h} \subset L \subset L_{d}$ that admits an embedding (over $L_{h}$) into $L_{m}$, which is impossible since the degrees $[L_{d}:L_{h}] = d/h$ and $[L_m:L_{h}] = m/h$ are co-prime. Thus, \eqref{eq:Yc-fib} is at least a \emph{birational} isomorphism. Since the source and target are both finite \'etale of degree $\ell/h = dm/h^2$ over $Y_{h}$, \eqref{eq:Yc-fib} is also an isomorphism. This concludes the proof.
\end{proof}

To streamline notation, we put $\bb{A}_1/K \ceq \J_1/K$. The next lemma makes precise the notion that $\bb{A}_n$ is the ``new part'' of $\J_n$, with the ``old part'' being the abelian sub-variety generated by the images $(\phi_{n,m})^*\J_m$, for $m \| n$. We note that the previous lemma furnishes the commutativity of the diagram \eqref{diag:Bn-Bm-Jn-Jm} below, which constitutes a crucial step in the proof.

\begin{lem}\label{lem:Jn-decomp-An}
Let $n\in \N$, and put $\B_n \ceq \prod_{m \mid n} \bb{A}_m$. Then, the homomorphism
\begin{equation}\label{eq:isog-An-Jn}
    \Theta_n \ceq \sum_{m \mid n} \phi_{n,m}^* \in \Hom({\B}_n,\J_n)
\end{equation}
is a $K$-isogeny. Consequently, the genus ${\g}_n$ of $\X_n$ and dimension $a_n$ of $\bb{A}_n$ satisfy the identities
\begin{equation}\label{eq:gn-an}
    g_n = \sum_{m \mid n} a_m \qquad \textup{and} \qquad a_n = \sum_{m \mid n} \mu(m)g_{n/m},
\end{equation}
where $\mu(-)$ denotes the M\"obius function.
\end{lem}


We call the isogeny $\Theta_n: \B_n \too \J_n$ the \defi{new parts decomposition} of $\J_n$.

\begin{proof}
The second identity in \eqref{eq:gn-an} follows by applying M\"obius inversion to the first identity, which itself follows from the claim that $\Theta_n$ is a $K$-isogeny, so we need only prove the latter. 
Recall that the \defi{isogeny category} $\Isog_K$ is the abelian category in which the objects are abelian varieties over $K$ and the hom-sets, for objects $A,B \in \Isog_K$, are the $\Q$-vector spaces $\Hom^0(A,B) = \Hom(A,B) \otimes_{\Z} \Q$. The Poincar\'e Reducibility Theorem implies that $\Isog_K$ is a semi-simple category, i.e. that any $J \in \Isog_K$ admits a decomposition $J \cong B_1^{e_1} \times \dotsb \times B_r^{e_r}$, in which each $B_i/K$ is simple. We note here that giving an isomorphism $A \too B$ in $\Isog_K$ is equivalent to giving a $K$-isogeny $A \too B$. So, in the present situation our goal is to show that {$\Theta_n \in \Hom^0(\B_n, \J_n)$ is an isomorphism}. To that end, we define:
\begin{equation*}
    \wt{\Theta}_n  \ceq \sum_{m \mid n} \dfrac{(\phi_{n,m}^*)}{\deg \phi_{n,m}} =  \sum_{m \mid n} \dfrac{m}{n}(\phi_{n,m}^*)  \in \Hom^0(\B_n,\J_n).
\end{equation*}
Let $(m)_{m\mid n} \in \Hom(\B_n,\B_n)$ denote the $K$-isogeny given by multiplication-by-$m$ on each $\bb{A}_m$. Then
\begin{equation*}
    \Theta_n = \wt{\Theta}_n \cdot (m)_{m\mid n}\in \Hom^0(\B_n,\J_n),
\end{equation*}
so it suffices to prove that: \emph{$\wt{\Theta}_n\in \Hom^0(\B_n,\J_n)$ is an isomorphism}. 
To see this, fix $n\in \N$,  take any $m\mid n$, let $\pi_{n,m}:\B_n \too \B_m$ be the natural projection, and consider the following diagram of $\Isog_K$:
\begin{equation}\label{diag:Bn-Bm-Jn-Jm}
    \begin{tikzcd}
        \B_n  \arrow[d, "\wt{\Theta}_{n}",swap] \arrow[rr, "\pi_{n,m}"] && \B_m  \arrow[d, "\wt{\Theta}_{m}"]\\
        \J_n \arrow[rr, "(\phi_{n,m})_*"] && \J_m.
    \end{tikzcd}
\end{equation}
Suppose $d$ is any divisor of $n$. If $d$ also divides $m$, then we have 
\begin{equation*}
    \wt{\Theta}_m \cdot \pi_{n,m} = \dfrac{d}{m} (\phi^*_{m,d}) = (\phi_{n,m})_* \cdot \wt{\Theta}_n \in \Hom^0(\bb{A}_d,\J_m). 
\end{equation*}
One the other hand, if $d$ does \emph{not} divide $m$, then \Cref{lem:fib-prod-deg} says that 
\begin{equation*}
    \wt{\Theta}_m \cdot \pi_{n,m} = 0 = (\phi_{n,m})_* \cdot \wt{\Theta}_n \in \Hom^0(\bb{A}_d,\J_m).  
\end{equation*}
Thus, we have $\wt{\Theta}_m \pi_{n,m} =(\phi_{n,m})_* \wt{\Theta}_n \in \Hom^0(\bb{A}_d,\J_m)$ for every divisor $d$ of $n$, and it follows that the diagram \eqref{diag:Bn-Bm-Jn-Jm} commutes. The projections $\pi_{n,m}$ arrange the $\B_n$'s into an inverse system, and the commutativity of \eqref{diag:Bn-Bm-Jn-Jm} implies that the $\wt{\Theta}_m$'s extend to a morphism of the respective inverse limits:
\begin{equation}\label{eq:varinj-isog}
    \varprojlim \wt{\Theta}_m : \prod_{m \| n} \bb{A}_m = \varprojlim_{m \| n} {\B}_m \too \varprojlim_{m \| n} \J_m.
\end{equation}
This morphism fits into the following commutative diagram in $\Isog_K$:
\begin{equation}\label{diag:Bn-Jn-exact-seq}
    \begin{tikzcd}
        1 \arrow[r]& \bb{A}_n \arrow[d, equals] \arrow[r] & \B_n  \arrow[rr, "\varprojlim \pi_{n,m}"] \arrow[d, "\wt{\Theta}_n",swap] && \displaystyle{\prod_{m \| n} \bb{A}_m} \arrow[d, "\varprojlim \wt{\Theta}_m"] \arrow[r] & 1\\
        1 \arrow[r] & \bb{A}_n \arrow[r] & \J_n \arrow[rr, "\varprojlim (\phi_{n,m})_*"] && \displaystyle{\varprojlim_{m \| n} \J_m},
    \end{tikzcd}
\end{equation}
in which the rows are exact. Now, we show that $\wt{\Theta}_n$ is an isomorphism by (strong) induction on $n$. The base-case $n=1$ is trivial, since $\bb{A}_1 = \B_1 = \J_1$ and $\wt{\Theta}_1$ is simply the identity. Next, suppose we have some $n\geqslant 2$ such that  $\wt{\Theta}_m$ is an isomorphism for every $m<n$. Then, the right vertical arrow in \eqref{diag:Bn-Jn-exact-seq} is an isomorphism, so we may apply the Snake Lemma to conclude that the middle arrow $\wt{\Theta}_n$ has trivial kernel and co-kernel, and hence, is an isomorphism. 
\end{proof}

\begin{rem}
We note in passing that $\wt{\Theta}_n$ being an isomorphism in $\Isog_K$  implies after the fact that the bottom right arrow of \eqref{diag:Bn-Jn-exact-seq} is surjective.
\end{rem}

\section{Second reduction: From $[\mbWbul{n}]$ to $[\bul{\mbV}{n}]$}\label{sec:Wn-to-Vn}

\begin{emp}
In this section, we assume for simplicity that $k$ is algebraically closed; this does not pose any issues because the main result of this section, \Cref{prop:Wn-Vn-decomp} below, continues to hold provided that $k$ contains the appropriate root of unity (needed to consider the spaces $\mbWbul{n}$ and $\wtmbWbul{n}$). 

 For each $n\in \N$ and  $\bullet \in \{ \textup{\textsf{cyc}, \textsf{re}, \textsf{im}, \textsf{dih}} \}$, we denote the new part of $\Jbul{n}/K$ (the Jacobian of $\Xbul{n}/K$) by $\Abul{n}/K$. We put $\Bbul{n} \ceq \prod_{m \mid n} \Abul{m}$, and we denote by 
\begin{equation*}
    \bul{\Theta}{n} : \Bbul{n} \too \Jbul{n}
\end{equation*}
the new parts decomposition provided by \Cref{lem:Jn-decomp-An}. 

Fix $n\geqslant 2$ and $\bullet \in \{ \textup{\textsf{cyc}, \textsf{re}, \textsf{im}}\}$. Let $\kk_n$ denote the corresponding field in \Cref{tab:data}. Observe from the proof of \Cref{thm:point-rank-bounds} (a) that for $i=1,\dotsc,r$, the point $\Qbul{i,n} \in \Xbul{n}(\kk_n)$ generates a torsor (cf. \ref{emp:tools-setup}) under the cyclic group
\begin{equation*}
    G_{i,n} \ceq \langle \sigma_{i,n} \rangle \subset  \gal(\kk_n/\kk_1). 
\end{equation*}
So, for each $i$ the space $\Wbul{i,n}$ from \ref{emp:Wbul} is free of rank one over the group algebra $\Q G_{i,n} = \spn_{\Q} \{1, \sigma_{i,n} , \dotsc , \sigma_{i,n}^{n-1} \}$, which we identify with the commutative ring $\Q[z]/(z^n - 1)$ in the obvious way. Recall the inverse cyclotomic polynomial $\Psi_n(x)$ from \eqref{eq:Psi}. 
Since $\sigma_{i,n}$ fixes $D_{\infty}$, we have $\Psi_n(\sigma_{i,n})D_{\infty} = \Psi_n(1) D_{\infty} = 0$. It follows that the divisor
\begin{equation*}
    \bul{\mf{P}}{i,n} \ceq \Psi_n({\sigma}_{i,n})\Qbul{i,n} = \Psi_n(\sigma_{i,n}) \bul{\mf{Q}}{i,n} \in \Div^0 (\X_n)_{\kk_n}
\end{equation*}
lies in $\Wbul{i,n}$. We define for each $i$ the subspace (as in \eqref{eq:V-[V]}):
\begin{equation*}
    \bul{V}{i,n} \ceq \spn_{\Q {G}_{i,n}} \{ \bul{\mf{P}}{i,n} \} \, \subset \Wbul{i,n}.
\end{equation*}
Adding these subspaces up for $i=1,\dotsc, r$, we obtain the subspace
\begin{equation*}
    \bul{\mbV}{n} \ceq \bul{V}{1,n}  \oplus \dotsb \oplus \bul{V}{r,n} \subset \mbWbul{n}.
\end{equation*}
To complete the picture, we define
\begin{equation*}
    \bul{\mbV}{1} = \mbWbul{1} \subset \Jbul{1}(\kk_1)_{\Q} = \J(\kk)_{\Q}. 
\end{equation*}
In the next proposition, we observe that there is a decomposition of $[\mbWbul{n}]$ into the direct sum $ \oplus_{m \mid n} [\bul{\mbV}{m}]$ which mirrors the new parts decomposition of $\Jbul{n}$. 
\end{emp}

\begin{prop}\label{prop:Wn-Vn-decomp}
For any $n\in \N$ and $\bullet \in \{ \textup{\textsf{cyc}, \textsf{re}, \textsf{im}} \}$, we have the following.
\begin{enumerate}[(a),itemindent=1cm, leftmargin=0pt]
    \item The embedding $[\bul{\mbV}{n}] \inj  \Jbul{n}(\kk_n)_{\Q}$ has image contained in $\Abul{n}(\kk_n)_{\Q}$.
    \item The new parts decomposition $\bul{\Theta}{n}: \Bbul{n} \too \Jbul{n}$  induces a $\Q$-isomorphism
    \begin{equation*}
        \bul{\Theta}{n}: \bigoplus_{m \mid n} \,[\bul{\mbV}{m}] \iso [\mbWbul{n}] \subset \Jbul{n}(\kk_n)_{\Q}. 
    \end{equation*}
    Consequently, if we write $\bul{\rho}{n}$ for the dimension of $[\bul{\mbV}{n}]_{\Q}$, we have
    \begin{equation*}
        \bul{R}{n} = \sum_{m \mid n} \bul{\rho}{m}.
    \end{equation*}
\end{enumerate}
\end{prop}

\begin{proof}
For part (a), start by fixing $i \in \{1,\dotsc,r\}$. Let $m$ be any divisor of $n$, and observe that $\phibul{n,m}(\Qbul{i,n}) = \Qbul{i,m} \in \Xbul{m}(\kk_n)$. In fact, $\Qbul{i,m}$ is a $\kk_m$-point of $\Xbul{m}$, and since the restriction of  $\sigma_{i,n} \in \aut_K(\kk_n)$ to the subfield $\kk_m$ is the automorphism $\sigma_{i,m} \in \aut_K(\kk_m)$, we have $\sigma_{i,n} \Qbul{i,m} = \sigma_{i,m} \Qbul{i,m}$. Thus, for any $j \in \Z$ we have
\begin{equation*}
\begin{tikzcd}
    \sigma_{i,n}^j\Qbul{i,n} \arrow[r, mapsto, "{\phibul{n,m}}"] & \sigma_{i,n}^j\Qbul{i,m} = \sigma_{i,m}^j\Qbul{i,m} \in \Xbul{m}(\kk_n).
\end{tikzcd}
\end{equation*}
It follows that the pushforward map $(\phibul{n,m})_*$ on divisors (over $\kk_n$) satisfies
\begin{equation*}
\begin{tikzcd}
    \bul{\mf{P}}{i,n} \arrow[rr, mapsto, "{(\phibul{n,m})_*}"] && {\Psi_m(\sigma_{i,m})\Qbul{i,m}} \in \Div^0 \Xbul{m}({\kk_n}).
\end{tikzcd}
\end{equation*}
Now, if $m<n$, then $x^m-1$ divides $\Psi_{n}(x)$ in $\Q[x]$ (cf. the definition of $\phi_n(x)$ in \eqref{eq:Psi}). Thus, we have $\Psi_n(\sigma_{i,m}) = 0 \in \Z G_{i,m}$, giving
\begin{equation*}
\begin{tikzcd}
    {[\bul{\mf{P}}{i,n}]} \arrow[rr, mapsto, "{(\phibul{n,m})_*}"] && {[\Psi_m(\sigma_{i,m})\Qbul{i,m}]} = 0 \in \Jbul{m}({\kk_n}).
\end{tikzcd}
\end{equation*}
By letting $m$ range over the divisors of $n$ less than $n$, we find that $[\bul{\mf{P}}{i,n}]$ is the kernel of each push-forward map $(\phibul{n,m})_*$, and since $\Abul{n}$ is the connected component of the intersections of these kernels, we conclude that some multiple of  $[\bul{\mf{P}}{i,n}]$ is contained in $\Abul{n}(\kk_n)$. This implies that $[\bul{V}{i,n}]$ is contained in $\Abul{n}(\kk_n)_{\Q}$, and finally, by allowing $i$ to range over $\{1,\dotsc,r \}$, we conclude that $[\bul{\mbV}{n}] = \sum_{i} [\bul{V}{i,n}]_{\Q}$ is also contained in $\Abul{n}(\kk_n)_{\Q}$, which concludes the proof of (a). 

We turn now to the proof of (b). Observe that the inverse system $(\J_n, (\phibul{n,m})_*)$ of abelian varieties gives rise an inverse system of $\Q$-vector spaces $([\mbWbul{m}], (\phibul{n,m})_*)$. For $n\in \N$, consider the natural representation of $\prod_{i=1}^r G_{i,n} = \gal(\kk_n/\kk_1)$ on $[\mbWbul{n}] = \sum_{i=1}^r [\Wbul{i,n}]$. For each $i$, the subspace  $[\bul{V}{i,n}] \subset [\Wbul{i,n}]$ is an invariant subspace on which  $\bul{\sigma}{i,n}$ acts by multiplication by a primitive $n$-th root of unity on $\bul{\mf{Q}}{i,n}$, i.e. $[\bul{V}{i,n}]_{\C}$ decomposes as the product $\oplus_{\chi} V_{\chi}$, where $\chi$ ranges over the $\phi(n)$-many embeddings $Z_N \inj \units{\C}$.  For any $m \mid n$, the homomorphism $(\phibul{n,m})^* : [\bul{V}{i,m}] \inj [\bul{W}{i,n}]$ is an embedding whose image is similarly the invariant subspace on which $\bul{G}{i,n}$ acts by multiplication by a primitive $m$-th root of unity. As we range over all $m$ dividing $n$, the direct sum of the subspaces $[\bul{V}{i,m}]$ equals $[\Wbul{i,n}]$. Summing over $i=1,\dotsc,r$, we conclude that 
\begin{equation*}
    \bigoplus_{m \mid n} [\bul{\mbV}{m}] = \bigoplus_{m \mid n} \sum_{i=1}^r [\bul{V}{i,n}] \iso  \sum_{i=1}^r [\Wbul{i,n}] = [\mbWbul{n}]. \qedhere
\end{equation*}
\end{proof}

We conclude this section with a last lemma to simplify the task of determining the integers $\bul{\rho}{n}$. 

\begin{lem}\label{lem:mbV-direct-sum}
Let $n\geqslant 2$, and let $\bullet \in \{ \textup{\textsf{cyc}, \textsf{re}, \textsf{im}} \}$. Then, the natural addition map
\begin{equation*}
    [\bul{V}{1,n}] \oplus \dotsb \oplus [\bul{V}{r,n}] \surj  [\bul{\mbV}{n}]
\end{equation*}
is an isomorphism of $\Q$-vector spaces. Consequently, $\bul{\rho}{n} = r\phi(n)$ if and only if $[\bul{\mf{P}}{i,n}] \in \Jbul{n}(\kk_n)$ is a point of infinite order for each $i=1,\dotsc,r$.
\end{lem}

\begin{proof}
Note first that for each $i$, since $[\bul{V}{i,n}]$ is annihilated by $\Phi_n(\sigma_{i,n})$, \Cref{lem:V-[V]-basic} implies that it is either trivial or else free of rank one over $\Q(\zeta_n) = \Q[\sigma_{i,n}]/\Phi_n(\sigma_{i,n})$, with the latter being the case if and only if $[\bul{\mf{P}}{i,n}] \in \Jbul{n}(\kk_n)$. Suppose for a contradiction, then, that we have elements $m_i \in [\bul{V}{i,n}]$, with $m_e \neq 0$ for some $e$, such that $m_1 + \dotsb + m_r = 0 [\bul{\mbV}{n}]$. Observe that $\sigma_{e,n}$ acts trivially on $m_i$ for $i \neq e$, since it acts trivially on the point $\Qbul{i,n}$, for $i \neq e$. Thus, multiplying the relation $\sum_i m_i = 0$ by $(\sigma_{e,n} - 1)$ results in the relation $(\sigma_{e,n} - 1) m_e = 0 \in [\bul{V}{e,n}]$. This forces $[\bul{V}{e,n}] = 0$, a contradiction. 
\end{proof}

\section{Determining $\cyc{\rho}{n},\Re{\rho}{n},$ and $\Im{\rho}{n}$, and finishing the proof of \Cref{thm:point-rank-bounds}}\label{sec:rank-proof}

In this section, we finish up the proof of \Cref{thm:point-rank-bounds} (b) by determining the integers $\bul{\rho}{n} = \dim_{\Q} [\bul{\mbV}{n}]$, for $\bullet \in \{\textup{\textsf{cyc}, \textsf{re}, \textsf{im}}\}$. Since this integer is geometric in nature (i.e. it does not change upon extensions of the ground field $k$), we again assume without loss of generality for this section that $k$ is algebraically closed. 

There is a trivial case when $\bul{\rho}{n} = 0$, namely, if $n$ and $d$ are small enough that $\dim \Abul{n} = 0$. It turns out that if $d$ and $n$ are large enough that $\dim \Abul{n} > 0$, then there are essentially only two distinct scenarios:
\begin{enumerate}[$\bullet$, itemindent=1cm,leftmargin=0pt]
    \item For $n=1$, we have $\Xcyc{1} = \Xre{1} = \X$ (the original curve in \ref{emp:mestre}), and for $i=1,\dotsc,r$, we have $\Qcyc{i,1} = \Qre{i,1} =  (u_i,v_i)$ (the points in \eqref{eq:Pi}). In this case the equality 
    \begin{equation*}
        \dim_{\Q} [\cyc{\mbV}{1}] = \dim_{\Q} [\Re{\mbV}{1}] = r - 1
    \end{equation*}
    was proved by Shioda in~\cite{shiodasymmetry}*{Theorems~4, 5} (as we mentioned in \ref{emp:mestre}). We show that similarly that $\dim_{\Q} [\Im{\mbV}{1}] = r$.
    \item For $n\geqslant 2$, we show using the tools in \Cref{sec:tools} that each $[\bul{\mf{P}}{i,n}] \in \J_n(\kk_n)$ is of infinite order, which implies that the maximum possible dimension $\bul{\rho}{n} = r\phi(n)$ is attained. 
\end{enumerate}


\begin{thm}\label{thm:dim-Vn}
Each integer $\bul{\rho}{n} = \dim_{\Q} [\bul{\mbV}{n}]$ is given by the values in the table below.  
\begin{table}[H]
    \centering
    \renewcommand{\arraystretch}{1.3}
    \setlength{\tabcolsep}{5pt}
    \adjustbox{max width=\textwidth}{
    \caption{Dimension of $[\cyc{\mbV}{n}],[\Re{\mbV}{n}],$ and $[\Im{\mbV}{n}]$}\label{tab:re-cyc-im-bounds} 
    \begin{tabular}{|r||c|c|c||c|c||c|c|} 
    \hline
    & \multicolumn{3}{c||}{$d=1$} & \multicolumn{2}{c||}{$d=2$} & \multicolumn{2}{c|}{$d\geqslant 3$} \\
    \hhline{|-||---||--||--|}
    & $n=1$ & $n=2$ & $n\geqslant 3$ & $n=1$ & $n\geqslant 2$ & $n=1$ & $n\geqslant 2$ \\
    \hhline{|=||===||==||==|}
    $\cyc{\rho}{n}= \Re{\rho}{n}= $ &  $0$ & $0$ & $4\phi(n)$ & $0$ & $6\phi(n)$ & $r-1$ & $r\phi(n)$\\
    \hhline{|-||-|-|-||-|-||-|-|}
    $\Im{\rho}{n}= $ & $4$ & $0$ & $4\phi(n)$ & $6$ & $6\phi(n)$ & $r$ & $r\phi(n)$ \\
    \hline
    \end{tabular}}
    \end{table}

\end{thm}

\begin{proof}
We have $\gcyc{n} = \gre{n}$,  and these equal $0$ if and only if $(n,d) = (1,1),(1,2),$ or $(2,1)$ (note: $\gim{n} = \gre{n}+1$, and this never equals $0$). If $d=1$, then $\dim \Aim{2} = \gim{2} - \gim{1} = 1 - 1 = 0$. In all these cases the equality  $\dim_{\Q} [\bul{\mbV}{n}]_{\Q} = 0$ is forced for the trivial reason that the ambient abelian variety $\Abul{n}$ is of dimension zero. 

Next, let us show that $\dim_{\Q} [\Im{\mbV}{1}]_{\Q} = r$ for all $d\geqslant 1$. To ease notation, write $Q_i$ for $\Qim{i,1}$ and put $\mf{Q}_i \ceq \Im{\mf{Q}}{i,1} = 2Q_i - D_{\infty}$  so that,  by definition, $[\Im{\mbV}{1}]_{\Q}$ is spanned by $[\mf{Q}_1],\dotsc,[\mf{Q}_r] \in \Jim{1}(\kk_1)$. We want to show that these points are linearly independent. Suppose for a contradiction that they are not, and that we have a non-trivial linear dependence relation
\begin{equation}\label{eq:dep-rel}
    a_1 [\mf{Q}_1] + \dotsb + a_r [\mf{Q}_r] = 0 \in \Jim{1}(\kk_1), 
\end{equation}
in which $a_s \neq 0$ for some $s \in \{1,\dotsc,r\}$. Here we have $\kk_1 = \kkdih{1} = k(\udih{1,1},\dotsc, \udih{r,1})$. Recall from \eqref{eq:delta-i-n} that we have an involution $\delta_s \ceq \delta_{s,1} \in \aut_K(\kk_1)$ which sends $\udih{s,1} \mapsto (\udih{s,1})\inv$ and which fixes $\udih{i,1}$ for $i\neq s$. Since $\uim{i,1} = (\udih{i,1}) - (\udih{i,1})\inv$, we see that $\delta_{s}$ fixes $\uim{i,1}$ for $i\neq s$, and $\delta_{s}$ sends $\uim{i,1} \mapsto -\uim{i,1}$. Since $Q_i = (\ure{i,1},v_i\uim{i,1})$, we find that
\begin{equation*}
    \delta_{s} (Q_i) = \ccases{
    Q_i & \; \textup{ if } \; i \neq s,\\
    \iota Q_i & \; \textup{ if } \; i=s,
    }
\end{equation*}
which in turn implies that
\begin{equation*}
    \delta_{s} ([\mf{D}_{i}]) = \ccases{
    [\mf{Q}_{i}] & \; \textup{ if } \; i \neq s,\\
    -[\mf{Q}_{i}] & \; \textup{ if } \; i=s,
    }
\end{equation*}
So, multiplying the relation \eqref{eq:dep-rel} by $(\delta_{s} - 1)$, we obtain
\begin{equation*}
    (\delta_{s} - 1) a_s [\mf{Q}_{s}] = -2a_s [\mf{Q}_{s}] = 0 \in \Jim{1}(\kk_1).
\end{equation*}
So, $[\mf{Q}_{s}]$ must be a torsion point. But since $Q_s \in \Xim{1}(\kk_1)$ is a transcendental point of $\Xim{1}/K$, this  contradicts \Cref{lem:dP-D-inf-order}, which asserts that $[\mf{Q}_s]$ is of infinite order.

It remains only to consider the cases with $n\geqslant 2$ in which $\Abul{n}$ is positive dimensional. We proceed now in two parts, depending on the genus $\gbul{n}$ of $\Xbul{n}$. 
\begin{enumerate}[(i),itemindent=1cm,leftmargin=0cm]
    \item \emph{The case $\gbul{n}=1$.} We have $\gim{n} = 1$ if and only if $nd=2$, which holds when $(n,d) = (2,1)$, and we have already dealt with this case. On the other hand, we have $\gcyc{n} = \gre{n} = 1$ if and only if $(n,d) = (3,1),(4,1),$ or $(2,2)$. Suppose we are in one of these three cases, so that $n$ is a power of a prime $p$ (i.e. $p=2$ or $3$), and assume for a contradiction that $\dim_{\Q} [\bul{V}{i,n}] < \phi(n)$. To ease notation, we everywhere drop the subscript ''$i,n$'' and superscript $\bullet$ on the point $\Qbul{i,n}$ and divisors $\bul{\mf{P}}{i,n}$ and $\bul{\mf{Q}}{i,n}$. Put $\tau \ceq \sigma_{i,n}^{n/p} \in \aut_K(\kk_n)$, which is an element of order $p$.  \Cref{lem:V-[V]-basic} says that $[\mf{P}] = [(\tau - 1)Q] =  (\tau - 1)[\mf{\Q}]$ is a torsion point of $\Jbul{n}(\kk_n)$. For $j=2,\dotsc,p$, we have
        \begin{equation}\label{eq:mfD-diff}
            (\tau^j -1)[\mf{Q}] = [2(\tau^j-1)Q] = 2(1+ \tau  \dotsb + \tau^{j-1})[\mf{P}] = 0 \in \Jbul{n}(\kk_n).
        \end{equation}
        The quotient by $\langle \tau \rangle$ is the degree $p$ cover $\phibul{n,n/p}: \Xbul{n} \too \Xbul{n/p}$. The target here is of genus $0$, so putting $\mf{E} \ceq (1+ \tau + \dotsb + \tau^{p-1})\mf{Q} \in \Div^0 (\Xbul{n})_{\kk_n},$ we have
        \begin{equation}\label{eq:mfD-sum}
            2[\mf{E}] = (1+ \tau + \dotsb + \tau^{p-1}) [\mf{Q}] = (\phibul{n,n/p})^*(\phibul{n,n/p})_* [\mf{Q}] = 0 \in \Jbul{n}(\kk_n).
        \end{equation}
        So, \eqref{eq:mfD-diff} and \eqref{eq:mfD-sum} together imply that
        \begin{equation*}
            \spn_{\Q} \{ [\mf{Q}], \tau[\mf{Q}] \dotsc, \tau^{p-1}[\mf{Q}] \}  = \spn_{\Q} \{ [\mf{E}],(\tau-1)[\mf{Q}], \dotsc, (\tau^{p-1}-1)[\mf{Q}]\} = \{0 \}.
        \end{equation*}
        But this implies that $[\mf{Q}] \in \Jbul{n}(\kk_n)$ is a torsion point, which again contradicts \Cref{lem:dP-D-inf-order} since $Q$ is a transcendental point of $\Xbul{n}/K$. This concludes the proof. 
    \item \emph{The case $\gbul{n}\geqslant 2$.} We will be done if we are able to apply \Cref{lem:Vn-phin} to the curve $\Xbul{n}/K$ and point $\Qbul{i,n} \in \Xbul{n}(\kk_n)$. To that end, let us start by noting that in all three cases (real, cyclic, and imaginary), and for all $d$ and $n$, assumption \ref{item:P-neq-iotaP} of the lemma is satisfied- the point $\Qbul{i,n}$ and its hyperelliptic conjugate $\iota \Qbul{i,n}$ each generate distinct $G_{i,n}$-torsors in $\Xbul{n}(\kk_n)$, i.e. the points $\sigma_{i,n}^j \Qbul{i,n}, \iota \sigma_{i,n}^j \Qbul{i,n} \in \Xbul{n}(\kk_n)$ are all distinct. Again in all three cases and for all $d$ and $n$, assumption \ref{item:tors} of the lemma is satisfied- since we have assumed $k$ to be algebraically closed, \Cref{prop:Ln-fact} (c) says that $K$ is algebraically closed in $\kk_n$, and \Cref{lem:alg-closed-tors} then says that we have an equality of torsion subgroups $\Jbul{n}(\kk_n)_{\tors} = \Jbul{n}(K)_{\tors}$. Recall now that the third (and final) assumption \ref{item:g-norm} of \Cref{lem:Vn-phin} requires that $\gbul{n} \geqslant \max\{ 2,2^{\norm{n}-1}\}$, where $\norm{n}$ denotes the number of distinct prime factors of $n$. Assumption \ref{item:g-norm} automatically holds if $n$ is a prime power so we need only consider the case when $\norm{n}\geqslant 2$. In this case, using the fact that $\gre{n} =\gcyc{n} = \gim{n} - 1 = (nd - 1 - \left[nd \right])/2$ (recall that $\left[nd \right]\in \{0,1\}$ denotes the remainder of $nd$ mod $2$), we have the chain of inequalities
        \begin{equation*}
            \begin{alignedat}{1}
            \gbul{n}    & \geqslant ({nd - 1 - \left[ nd \right]})/{2}\\
                        & \geqslant ({nd}-2)/{2}\\
                        & \geqslant (n/2) - 1\\
                        & \geqslant (2^{\norm{n}-1}3)/2 - 1\\
                        & = 2^{\norm{n}-1}(1+1/2) - 1\\
                        & = 2^{\norm{n}-1} + (2^{\norm{n}-2} - 1)\\
                        & \geqslant 2^{\norm{n}-1}.
            \end{alignedat}
        \end{equation*}
        Thus, if $\gbul{n} \geqslant 2$ then assumption \ref{item:g-norm} of \Cref{lem:Vn-phin} is satisfied and the lemma implies that $\dim_{\Q} [\bul{V}{i,n}] = \phi(n)$. \qedhere 
        \end{enumerate}
\end{proof}

\begin{emp}\label{emp:proof-pt-rank-bounds}
\textit{Proof of \Cref{thm:point-rank-bounds}.} Fix the integers $d$ and $n$, and fix the label $\bullet \in \{ \textsf{cyc}, \textsf{re}, \textsf{im} \}$.  Then, one checks that the values of $\bul{\rho}{n}$ proved in \Cref{thm:dim-Vn}, together with the identity $\sum_{m \mid n} \phi(m) = n$, imply that $\bul{R}{n} = \sum_{m \mid n} \dim_{\Q} [\bul{\mbV}{n}]$ equals the value asserted in \Cref{tab:ranks}. Finally, the validity of \Cref{thm:point-rank-bounds} (b) follows from \Cref{lem:reduction-1}. \qedhere 
\end{emp}

\section{Non-isotriviality of the Real, Imaginary, and Dihedral families}\label{sec:non-iso}

The main result of this section is the following proposition (which was already used in the proof of \Cref{thm:main}). 

\begin{prop}[Proved in \ref{emp:proof-non-iso}]\label{prop:non-iso}
Any curve in the real, imaginary, and dihedral families in \Cref{tab:data} is non-isotrivial over $K =k(\mbh)$ (provided, of course, that its genus is positive). 
\end{prop}

We will prove the proposition with the help of two lemmas. For simplicity, and without loss of generality, we assume throughout that $k$ is algebraically closed. 

\begin{lem}\label{lem:defranchis}
Let $\phi: \X \too \Y$ be a cover of nice curves over $K$. If $\Y/K$ is non-isotrivial of genus $g\geqslant 1$, then $\X/K$ is also non-isotrivial. 
\end{lem}

\begin{proof}
Assume for a contradiction that $\Y/K$ is non-isotrivial of genus $g\geqslant 1$ but $\X/K$ is isotrivial. Choose a suitable smooth $k$-variety $\scrv$ with function field $K$ and spread out the cover $\phi:\X \too \Y$ to a finite $\scrv$-morphism $\Phi: \scrx \too \scry$. For any point $P \in \scrv(\ol{k})$, this specializes to a cover $\Phi_P: \scrx_P \too \scry_P$ of curves over $\ol{k}$. By assumption, as we vary the geometric point $P \in \scrv(\ol{k})$, the curves $\scrx_P/\ol{k}$ are all isomorphic to a fixed curve $X/\ol{k}$, whereas $\scry_P/\ol{k}$ varies over infinitely many isomorphism classes of genus $g$ curves. But this is absurd-- if $g \geqslant 2$, then $X/\ol{k}$ would cover infinitely many non-isomorphic curves of genus $g \geqslant 2$, which contradicts a classical theorem of de Franchis (see \cite{kani-defranchis}), whereas if $\Y$ is of genus one then the Jacobian of $X/\ol{k}$ would contain infinitely many elliptic curves over $\ol{k}$ in its isotypic decomposition. 
\end{proof}

Below, for an element $\alpha \in \ol{K}$, we write $[\alpha]$ to denote the corresponding geometric point on $\A^1_{\ol{K}} = \p^1_{\ol{K}} - \{ \infty \}$. 

\begin{lem}\label{lem:non-iso-G}
Let $F(x) \in K[x]$ be a separable polynomial of degree at least $5$, and let $G$ denote its Galois group. Assume that the hyperelliptic curve
\begin{equation*}
    \Y/K:\; y^2 = F(x)
\end{equation*}
is isotrivial. Then, there exists a group embedding
\begin{equation*}
    \theta: G \inj \aut_{\ol{K}}(\p^1_{\ol{K}}), \qquad g \mtoo \theta_g,
\end{equation*}
such that for every root $\alpha$ of $F$ in $\ol{K}$, we have $\theta_g([\alpha]) = [g(\alpha)].$ In particular, if $F(x)$ has a root in $k$ or $\deg F(x)$ is odd, then $G$ is cyclic.
\end{lem}

\begin{proof}
By assumption, there exists a curve $Y/k$ and a $\ol{K}$-isomorphism $\varphi : Y_{\ol{K}} \iso \Y$. Then $Y/k$ is hyperelliptic and admits a double cover $\pi_Y: Y \too \p^1_k$. Thus, there is an isomorphism $f_{\varphi} :\p^1_{\ol{K}} \iso \p^1_{\ol{K}}$ which makes the following diagram commute:
\begin{equation*}
    \begin{tikzcd}
        Y_{\ol{K}} \arrow[r, "\varphi"] \arrow[d, "\pi_{Y}", swap] & \Y \arrow[d, "\pi_{\Y}"],\\
        \p^1_{\ol{K}} \arrow[r, "f_{\varphi}"] & \p^1_{\ol{K}}. 
    \end{tikzcd}
\end{equation*}
We have a natural action of the Galois group $G_K \ceq \gal(\ol{K}/K)$ on $\p^1_{\ol{K}}$ in which every element $g \in G_K$ fixes $\infty$ and sends a point $[\alpha] \in \p^1_{\ol{K}} - \{\infty\}$ to $[g(\alpha)]$. We define the map (of sets)
\begin{equation}
    \theta : G_K \too \aut_{\ol{K}}(\p^1_{\ol{K}}), \qquad g \mtoo g \comp f_{\varphi} \comp g \inv \comp f_{\varphi} \inv.
\end{equation}
Let $\alpha$ be a root of $F(x)$, so that $[\alpha]$ is a branch point of $\pi_{\Y}$ and $P \ceq f_{\varphi}\inv [\alpha]$ is a branch point of $\pi_Y$. Since $P$ is fixed by $G_K$, we have
\begin{equation*}
    \theta_g ([\alpha]) = g \comp f_{\varphi} \comp g \inv(P) = g \comp f_{\varphi}(P) = [g(\alpha)]. 
\end{equation*}
For every $g$ and $h$ in $G_K$, the element $\theta_{gh}$ agrees with the element $\theta_g \theta_h$ when restricted to the branch locus of $\pi_{\Y}$, and since any automorphism of the projective line is determined by its effect on three points, it follows that $\theta_{gh} = \theta_g \theta_h$. Thus, $\theta:G \too \aut_{\ol{K}}(\p^1_{\ol{K}})$ is a group homomorphism. Visibly, $\theta_g$ is the identity if and only if $g$ fixes all the roots of $F(x)$. Thus, it factors through the Galois group of $F(x)$ to give a group embedding $\theta : G \inj \aut_{\ol{K}}(\p^1_{\ol{K}})$ satisfying $\theta_g([\alpha]) = [g(\alpha)]$ for every root $\alpha$ of $F(x)$. 

For the second statement, assume that $F(x)$ has a root $a \in k$, or that $\deg F$ is odd. In the latter case, let $P \ceq \infty$, and in the former case let $P \ceq [a]$. Then the image of $G$ under $\theta$ is a finite subgroup of $\aut_{\ol{K}}(\p^1_{\ol{K}})$ which fixes $P$. So, if $f:\p^1_{\ol{K}} \too \p^1_{\ol{K}}$ denotes the quotient by $G$, then $P$ is totally ramified under $f$ and $G$ is its inertia group. But since $\chr k = 0$, the ramification is necessarily tame and the inertia group of $P$ must be cyclic (see~ \cite[\href{https://stacks.math.columbia.edu/tag/09EE}{Lemma 09EE}]{stacks-project}). 
\end{proof}

\begin{emp}\label{emp:proof-non-iso}
\textit{Proof of \Cref{prop:non-iso}}.
For all $d\geqslant 1$, the curve $\Xim{1}/K$ is non-isotrivial because \emph{any} hyperelliptic curve $X/k$ of genus $\gim{1}$ arises as a specialization because, given a double cover $\pi:X \too \p^1_k$, we can choose the coordinate on $\p^1_k$ in such a way that $2,-2,$ and $\infty$ are all branch points of $\pi$, and such a curve is visibly a specialization of $\Xim{1}$. Since there exist covers $\Im{\psi}{n}: \Xdih{n} \to \Xim{n}$ and  $\phiim{n,1} :\Xim{n} \to \Xim{1}$,  \Cref{lem:defranchis} implies that $\Xim{n}/K$ is non-isotrivial for all $d\geqslant 1$ and $n\in \N$. 

By a similar argument, the curves $\Xwtre{1}/K$ and $\Xwtim{1}/K$ are non-isotrivial if $d\geqslant 2$, and it follows again from \Cref{lem:defranchis}
that $\Xwtre{n}/K,\Xwtim{n}/K,$ and $\Xwtdih{n}/K$ are non-isotrivial for all $n\in \N$ and $d\geqslant 2$. 

Suppose next that $d=1$. If $n=2$ or $3$ then $\Xwtre{n}/K$ and $\Xwtim{n}/K$ are both elliptic curves, and one checks that their $j$-invariant lies in $K \setminus k$. If $n\geqslant 4$, then $\wtbbfRe{n}$ and $\wtbbfIm{n}$ are both separable of degree $n+1 \geqslant 5$ having a root in $k$, and their Galois group over $K$ is dihedral of order $2n$. So in this case \Cref{lem:non-iso-G} implies that $\Xwtre{n}/K$ and $\Xwtim{n}/K$ are non-isotrivial. Again,  \Cref{lem:defranchis} implies that for $d=1$ and $n\geqslant 2$,  $\Xwtdih{n}/K$ is non-isotrivial (note: it is of genus zero when $d=n=1$).

It remains only to prove the statement for the curves $\Xre{n}/K$. If $d\geqslant 3$ then $\Xre{1} = \X$ (the original curve from Mestre's construction \ref{emp:mestre}). This is non-isotrivial over $K$, and so the same is true of $\Xre{n}/K$ for all $n\in \N$ by \Cref{lem:defranchis}. 

Next, suppose that $d=2$. If $n = 2$, then the double cover $\wtphire{2,1}: \Xre{2} \to \Xwtre{1}$ implies by \Cref{lem:defranchis} that $\Xre{2}/K$ is non-isotrivial. If $n\geqslant 3$ then $\bbfRe{n}(x)  = \bb{h}(\fRe{n}(x))$ is of degree  $nd \geqslant 6$. If $a_1$ and $a_2$ denote the roots of $\bb{h}(x)$ in $\ol{K}$, then $\bbfRe{n}(x)$ factors as $h_2(\fRe{n}(x) - a_1)(\fRe{n}(x) - a_2) \in \ol{K}[x]$. The Galois group of each factor contains the cyclic group $Z_n$ (cf. \ref{emp:dih-mono}), which implies that the Galois group  $G$ of $\bbfRe{n}(x)$ contains $Z_n \times Z_n$. Applying \Cref{lem:non-iso-G}, we find that if $\Xre{n}/K$ is isotrivial then $G$ is a subgroup of $\aut_{\ol{K}}(\p^1_{\ol{K}})$. This is impossible because $\aut_{\ol{K}}(\p^1_{\ol{K}})$ does not contains $Z_n \times Z_n$, for $n\geqslant 3$. We conclude that if $d=2$ and $n\geqslant 2$, then $\Xre{n}/K$ is non-isotrivial.

Finally, suppose that $d=1$. If $n=3$ then $\Xre{n}/K$ is the elliptic curve $y^2 = h_1(\fRe{3}(x)) + h_0$, and one can check that its $j$-invariant lies in $K \setminus k$. If $n=4$ then the double cover $\phire{4,2}: \Xre{4} \to \Xwtre{2}$ together with \Cref{lem:defranchis} implies that $\Xre{4}/K$ is non-isotrivial. If $n \geqslant 5$ is odd, then $\bbfRe{n}(x)$ is odd of degree $n$ and has Galois group $D_n$ over $K$, so \Cref{lem:non-iso-G} implies that $\Xre{n}/K$ is non-isotrivial. Now, any integer $n \geqslant 3$ is either a power of two, in which case we have the cover $\phire{n,4}: \Xre{n} \to \Xre{4}$, or else $n$ is divisible by an odd prime $p$, in which case we have the cover $\phire{n,p}: \Xre{n} \to \Xre{p}$. So, we conclude by \Cref{lem:defranchis} that if $d=1$ and $n\geqslant 3$, then $\Xre{n}/K$ is non-isotrivial. This concludes the proof. \qedhere
\end{emp}

\section{Complex multiplication of $\Acyc{n}/K$ and the proof of \Cref{thm:cm-twist}}\label{sec:cm}

For this section, we set $d=1$, so that the polynomial $\bbfcyc{n}(x)$ from \Cref{tab:data} is $\bb{h}(x^n) = h_1x^n + h_0$, and the associated hyperelliptic curve is $\Xcyc{n}/K: y^2 = h_1x^n + h_0.$ For any $P = (a,b) \in \A^2(k)$ such that $ab \neq 0$, the specialized curve 
\begin{equation}
    X_n/k:\; y^2 = a(x^n + b)
\end{equation}
is of the same genus as $\Xcyc{n}$. For all $n$ and all $m$ dividing $n$, we have again the natural covers $\phi_{n,m}: \X_n \to \X_m$ defined by $(x,y) \mapsto (x^{n/m},y)$ which make $(X_n, \phi_{n,m})$ into a $\zhat$-system of hyperelliptic curves over $k$. In particular, each cover $\phi_{n,1}: X_n \to X_1$ is the quotient (provided $k$ contains $\zeta_n$) by the order $n$ automorphism 
\begin{equation*}
    \sigma_n : X_n \too X_1, \qquad (x,y) \mtoo (\zeta_n x,y). 
\end{equation*}
Let $J_n/k$ denote the Jacobian of $X_n$, and let $A_n/k$ denote its new part. Then, by \Cref{lem:Jn-decomp-An}, we have $\dim A_n = 0$ if $n = 1$ or $2$, and
\begin{equation*}
    \dim A_n = \phi(n)/2, \qquad \textup{ if } n\geqslant 3.
\end{equation*}
Define the endomorphism
\begin{equation*}
    \ol{\sigma}_n : J_n \too J_n, \qquad [D] \mtoo [(\sigma_n)_* D].
\end{equation*}
We show below that this endomorphism endows $A_n/k$ with complex multiplication by the number field $\Q(\zeta_n)$.

\begin{prop}\label{prop:complex-mult}
Fix $n\geqslant 3$ (so that $X_n/k$ is of positive genus) and assume that $k$ contains $\zeta_n$. Then, the automorphism $\ol{\sigma}_n \in \aut_k(J_n)$ restricts to an automorphism $\ol{\zeta} : A_{n} \iso A_{n}$ which satisfies 
\begin{equation*}
    \Phi_n(\ol{\zeta}) = 0 \in \End_k^0(A_n).
\end{equation*}
Consequently, the $\Q$-subalgebra $\Q[\ol{\zeta}] \subset \End^0_k(A_n)$ is isomorphic to the cyclotomic field $\Q(\zeta_n)$ via the map $\ol{\zeta} \mapsto \zeta_n$, and $A_n$ has complex multiplication by $\Q(\zeta_n)$. 
\end{prop}

\begin{proof}
For any $m$ dividing $n$, we have the following commutative diagram of curves and an induced diagram of their Jacobians:
\begin{equation*}
    \begin{tikzcd}
        X_n \arrow[r, "\sigma_n"] \arrow[d, "\phi_{n,m}"'] & X_n \arrow[d, "\phi_{n,m}"], & & J_{n} \arrow[r, "\ol{\sigma}_n"] \arrow[d, "(\phi_{n,m})_*"'] & J_{n} \arrow[d, "(\phi_{n,m})_*"]\\
        X_m \arrow[r, "\sigma_m"] & X_m, & & J_{m} \arrow[r, "\ol{\sigma}_m"] & J_{m}.
    \end{tikzcd}
\end{equation*}
It follows at once that $\ol{\sigma}_n(A_n)$ is contained in the kernel of $(\phi_{n,m})_*$. By letting $m$ range over all divisors of $n$, it follows that $\ol{\sigma}_n(A_n) = A_n$, so $\ol{\sigma}_n$ indeed restricts to an automorphism $\ol{\zeta}$ of $A_n$. For every factorization $n=me$, with $m < n$, define the polynomial
\begin{equation*}
    q_m(x) \ceq (x^n - 1)/(x^m - 1) = 1 + x^m + \dotsb + (x^m)^{e-1} \in \Z[x],
\end{equation*}
and observe that 
\begin{equation*}
    (\phi_{n,m})^* \comp (\phi_{n,m})_* = q_m(\ol{\sigma}_n) \in \End_k^0(J_n).
\end{equation*}
It follows that for any $m \| n$, we have $q_m(\ol{\zeta}) = 0 \in \End_k(A_n)$. It is an exercise to see that the gcd of the polynomials $q_m(x)$, as $m$ ranges over the divisors of $n$ less than $n$, is the cyclotomic polynomial $\Phi_n(x)$. So, we have $\Phi_n(x) = \sum_{m \| n} q_m(x) f_m(x)$, for some set of polynomials $f_m(x) \in \Q[x]$. It follows that $\Phi_n(\ol{\zeta}) = 0 \in \End_k(A_n)$. Thus, the $\Q$-subalgebra $\Q[\ol{\zeta}] \subset \End^0_k(A_n)$ is isomorphic to $\Q(\zeta_n)$, and since $\dim A_n = \phi(n)/2$, we conclude that $A_n$ has complex multiplication by $\Q(\zeta_n)$.  
\end{proof}

\begin{rem}
The proof above serves equally well to show that $\Acyc{n}/\kkcyc{n}$ also has complex multiplication by $\Q(\zeta_n)$ (provided $k$ contains $\zeta_n$). 
\end{rem}

We conclude this section (and also the article) with:
\begin{emp}\label{emp:proof-thm-twist}
\textit{Proof of \Cref{thm:cm-twist}.} Fix $N \geqslant 3$, and a number field $k$ which contains $\zeta_N$. Consider the open subset $V \subset \A^4_k(\mbucyc{N})$ defined by the condition $m_0h_0h_1(\disc\, \bb{m}) \neq 0$. Let $P \in V(k)$ be a point, and put $a \ceq h_1(P)$ and $b \ceq h_0(P)$. Then, the entire $\zhat$-system $(\Xcyc{n}, \phicyc{n,m})$ can be specialized to the $\zhat$-system $(X_n,\phi_{n,m})$, where $X_n/k: y^2 = ax^n +b$ is the specialization of $\Xcyc{n}$ at $P$ and the covers $\phi_{n,m}$ are again defined by $(x,y) \mapsto (x^{n/m},y)$. Writing $J_n/k$ for the Jacobian of $X_n/k$, and $A_n/k$ for its new part, it follows that for all $n \in \N$, the specialization homomorphism $\Jcyc{n}(\kk_N) \too J_n(k)$ restricts to a homomorphism $\Acyc{n}(\kk_N) \too A_n(k)$. Note that if $n\geqslant 3$ then $A_n/k$ has complex multiplication by $\Q(\zeta_n)$ by \Cref{prop:complex-mult}.

Now, \Cref{thm:specialization} gives a thin subset $T \subset V(k)$ such that for all $P \in V(k) \setminus T$ and all $m \mid N$, the specialization homomorphism $\Jcyc{m}(\kk_N) \too J_m(k)$, and hence, also its restriction $\Acyc{m}(\kk_N) \inj A_m(k)$, is injective. By \Cref{thm:dim-Vn}, we therefore have $\rk A_m(k) \geqslant  4\phi(n) =8(\dim A_m)$. So, to complete the proof of \Cref{thm:cm-twist}, it remains only to show that: 

\emph{The image of $S \ceq V(k) \setminus T$ under the natural $k$-morphism $f:V \too B\ceq \spec k[h_1,h_0]$ is not thin in $B$.} 

Suppose for a contradiction that this is not the case, i.e. that there exists a generically finite cover $g:U \too B$ admitting no rational section and satisfying $g(U(k)) \supset  f(S)$. Let $\wt{U}$ denote the fibered product $V \times_{B} U$, so that the base-change $\wt{g} :\wt{U} \too V$ is a generically finite cover which satisfies $\wt{g}(\wt{U}(k)) \supset S$. By shrinking $V$ to a smaller open subset if necessary, we may assume that $g$ and $\wt{g}$ are both dominant $k$-morphisms. Now, the generic fiber of $g$ is of the form  
\begin{equation*}
    \sqcup_{i=1}^{\ell}\spec L_i \to \spec K
\end{equation*}
for some positive integer $\ell$, where each $L_i/K$ is an extension of degree at least $2$ (this last condition is equivalent to the statement that $g$ does not admit a rational section $B \DashedArrow U$). So, if we put $\L_i \ceq L_i \otimes_K \kkcyc{N}$, then the generic fiber of $\wt{g}$ is 
\begin{equation*}
    \sqcup_{i=1}^{\ell} \spec \L_i \too \spec \kkcyc{N}.
\end{equation*}
By \Cref{prop:Ln-fact} (c), $K$ is algebraically closed in $\kkcyc{N}$, which implies that each $\L_i/\kkcyc{N}$ is a field extension of the same degree as $L_i/K$. Thus, the cover $\wt{g}:\wt{U} \to V$ does not admit a rational section, which implies that $S$ is a thin subset of $V(k)$. But this implies that $V(k) = T \cup S$ is thin, which is a contradiction because $k$ is a number field (hence satisfies the Hilbert Property) and $V$ is a rational $k$-variety. \qed
\end{emp}

\bibliographystyle{alpha}
\bibliography{high-rank-PTE}

\begin{bibdiv}
\begin{biblist}

\bib{griff-harris}{book}{
      author={Arbarello, E.},
      author={Cornalba, M.},
      author={Griffiths, P.~A.},
      author={Harris, J.},
       title={Geometry of algebraic curves. {V}ol. {I}},
      series={Grundlehren der mathematischen Wissenschaften [Fundamental Principles of Mathematical Sciences]},
   publisher={Springer-Verlag, New York},
        date={1985},
      volume={267},
        ISBN={0-387-90997-4},
         url={https://doi.org/10.1007/978-1-4757-5323-3},
      review={\MR{770932}},
}

\bib{beauville}{incollection}{
      author={Beauville, Arnaud},
       title={Finite subgroups of {${\rm PGL}_2(K)$}},
        date={2010},
   booktitle={Vector bundles and complex geometry},
      series={Contemp. Math.},
      volume={522},
   publisher={Amer. Math. Soc., Providence, RI},
       pages={23\ndash 29},
         url={https://doi.org/10.1090/conm/522/10289},
      review={\MR{2681719}},
}

\bib{neronmodels}{book}{
      author={Bosch, Siegfried},
      author={L\"{u}tkebohmert, Werner},
      author={Raynaud, Michel},
       title={N\'{e}ron models},
      series={Ergebnisse der Mathematik und ihrer Grenzgebiete (3) [Results in Mathematics and Related Areas (3)]},
   publisher={Springer-Verlag, Berlin},
        date={1990},
      volume={21},
        ISBN={3-540-50587-3},
         url={https://doi.org/10.1007/978-3-642-51438-8},
      review={\MR{1045822}},
}

\bib{capratpts}{incollection}{
      author={Caporaso, Lucia},
       title={Counting rational points on algebraic curves},
        date={1995},
      volume={53},
       pages={223\ndash 229},
        note={Number theory, I (Rome, 1995)},
      review={\MR{1452380}},
}

\bib{ulmerpoints}{article}{
      author={Concei\c{c}\~{a}o, Ricardo},
      author={Ulmer, Douglas},
      author={Voloch, Jos\'{e}~Felipe},
       title={Unboundedness of the number of rational points on curves over function fields},
        date={2012},
     journal={New York J. Math.},
      volume={18},
       pages={291\ndash 293},
         url={http://nyjm.albany.edu:8000/j/2012/18_291.html},
      review={\MR{2928577}},
}

\bib{CHM}{article}{
      author={Caporaso, Lucia},
      author={Harris, Joe},
      author={Mazur, Barry},
       title={Uniformity of rational points},
        date={1997},
        ISSN={0894-0347},
     journal={J. Amer. Math. Soc.},
      volume={10},
      number={1},
       pages={1\ndash 35},
         url={https://doi.org/10.1090/S0894-0347-97-00195-1},
      review={\MR{1325796}},
}

\bib{CTrankjump}{article}{
      author={Colliot~Th{\'e}l{\`e}ne, Jean-Louis},
       title={Point g{\'e}n{\'e}rique et saut du rang du groupe de {M}ordell--{W}eil},
        date={2020},
        ISSN={1730-6264},
     journal={Acta Arithmetica},
         url={http://dx.doi.org/10.4064/aa190814-18-3},
}

\bib{Elkies2006}{misc}{
      author={Elkies, Noam~D.},
       title={$\mathbb{Z}^{28}$ in ${E}(\mathbb{Q})$, etc.},
        date={2006},
        note={Email to the \texttt {NMBRTHRY@LISTSERV.NODAK.EDU} mailing list},
}

\bib{elkiesthreelectures}{misc}{
      author={Elkies, Noam~D.},
       title={Three lectures on elliptic surfaces and curves of high rank},
        date={2007},
        note={Preprint. \href{https://arxiv.org/abs/0709.2908}{arXiv: 0709.2908}},
}

\bib{elkiescover}{incollection}{
      author={Elkies, Noam~D.},
       title={About the cover: rational curves on a {$K3$} surface},
        date={2009},
   booktitle={Arithmetic geometry},
      series={Clay Math. Proc.},
      volume={8},
   publisher={Amer. Math. Soc., Providence, RI},
       pages={1\ndash 4},
         url={https://doi.org/10.1109/tc.2008.20},
      review={\MR{2498052}},
}

\bib{elkiesgenus4}{misc}{
      author={Elkies, Noam~D.},
       title={Curves with many points},
        date={1996},
        note={Available at \url{http://people.math.harvard.edu/~elkies/many_pts.pdf}},
}

\bib{faltings}{article}{
      author={Faltings, Gerd},
       title={Endlichkeitss\"{a}tze f\"{u}r abelsche {V}ariet\"{a}ten \"{u}ber {Z}ahlk\"{o}rpern},
        date={1983},
        ISSN={0020-9910},
     journal={Invent. Math.},
      volume={73},
      number={3},
       pages={349\ndash 366},
         url={https://doi.org/10.1007/BF01388432},
      review={\MR{718935}},
}

\bib{movinglemma}{article}{
      author={Gabber, Ofer},
      author={Liu, Qing},
      author={Lorenzini, Dino},
       title={Hypersurfaces in projective schemes and a moving lemma},
        date={2015},
        ISSN={0012-7094},
     journal={Duke Math. J.},
      volume={164},
      number={7},
       pages={1187\ndash 1270},
         url={https://doi.org/10.1215/00127094-2877293},
      review={\MR{3347315}},
}

\bib{kani-defranchis}{article}{
      author={Kani, Ernst},
       title={Bounds on the number of nonrational subfields of a function field},
        date={1986},
        ISSN={0020-9910,1432-1297},
     journal={Invent. Math.},
      volume={85},
      number={1},
       pages={185\ndash 198},
         url={https://doi.org/10.1007/BF01388797},
      review={\MR{842053}},
}

\bib{kulesz}{article}{
      author={Kulesz, Leopoldo},
       title={Jacobiennes de courbes alg\'{e}briques de genre 2 et 3 de grand rang sur $\mathbb{Q}$},
        date={2001},
        ISSN={0024-6107},
     journal={J. London Math. Soc. (2)},
      volume={63},
      number={2},
       pages={288\ndash 298},
         url={https://doi.org/10.1017/S0024610700001824},
      review={\MR{1810130}},
}

\bib{mestre11}{article}{
      author={Mestre, Jean-Fran\c{c}ois},
       title={Courbes elliptiques de rang {$\geq 11$} sur $\mathbb{Q}(t)$},
        date={1991},
        ISSN={0764-4442},
     journal={C. R. Acad. Sci. Paris S\'{e}r. I Math.},
      volume={313},
      number={3},
       pages={139\ndash 142},
      review={\MR{1121576}},
}

\bib{mestre12}{article}{
      author={Mestre, Jean-Fran\c{c}ois},
       title={Courbes elliptiques de rang {$\geq 12$} sur $\mathbb{Q}(t)$},
        date={1991},
        ISSN={0764-4442},
     journal={C. R. Acad. Sci. Paris S\'{e}r. I Math.},
      volume={313},
      number={4},
       pages={171\ndash 174},
      review={\MR{1122689}},
}

\bib{mumfordGIT}{book}{
      author={Mumford, D.},
      author={Fogarty, J.},
      author={Kirwan, F.},
       title={Geometric invariant theory},
     edition={Third},
      series={Ergebnisse der Mathematik und ihrer Grenzgebiete (2) [Results in Mathematics and Related Areas (2)]},
   publisher={Springer-Verlag, Berlin},
        date={1994},
      volume={34},
        ISBN={3-540-56963-4},
         url={https://doi.org/10.1007/978-3-642-57916-5},
      review={\MR{1304906}},
}

\bib{neronhilbert}{article}{
      author={N\'eron, Andr\'e},
       title={Probl\`emes arithm\'etique et g\'eom\'etriques rattach\'es \`a la notion de rang d'une courbe alg\'ebrique dans un corps},
        date={1952},
     journal={Bulletin de la Soci\'et\'e Math\'ematique de France},
      volume={80},
       pages={101\ndash 166},
         url={http://www.numdam.org/item/BSMF_1952__80__101_0},
      review={\MR{0056951}},
}

\bib{neronhilbert2}{inproceedings}{
      author={N\'{e}ron, Andr\'e},
       title={Propri\'{e}t\'{e}s arithm\'{e}tiques de certaines familles de courbes alg\'{e}briques},
        date={1956},
   booktitle={Proceedings of the {I}nternational {C}ongress of {M}athematicians, 1954, {A}msterdam, vol. {III}},
   publisher={Erven P. Noordhoff N.V., Groningen; North-Holland Publishing Co., Amsterdam},
       pages={481\ndash 488},
      review={\MR{0087210}},
}

\bib{odoni}{article}{
      author={Odoni, R. W.~K.},
       title={The {G}alois theory of iterates and composites of polynomials},
        date={1985},
        ISSN={0024-6115},
     journal={Proc. London Math. Soc. (3)},
      volume={51},
      number={3},
       pages={385\ndash 414},
         url={https://doi.org/10.1112/plms/s3-51.3.385},
      review={\MR{805714}},
}

\bib{Qpoints}{book}{
      author={Poonen, Bjorn},
       title={Rational points on varieties},
      series={Graduate Studies in Mathematics},
   publisher={American Mathematical Society, Providence, RI},
        date={2017},
      volume={186},
        ISBN={978-1-4704-3773-2},
      review={\MR{3729254}},
}

\bib{Park-wood}{article}{
      author={Park, Jennifer},
      author={Poonen, Bjorn},
      author={Voight, John},
      author={Wood, Melanie~Matchett},
       title={A heuristic for boundedness of ranks of elliptic curves},
        date={2019},
        ISSN={1435-9855},
     journal={J. Eur. Math. Soc. (JEMS)},
      volume={21},
      number={9},
       pages={2859\ndash 2903},
         url={https://doi.org/10.4171/JEMS/893},
      review={\MR{3985613}},
}

\bib{serremw}{book}{
      author={Serre, Jean-Pierre},
       title={Lectures on the {M}ordell-{W}eil theorem},
     edition={Third Edition},
      series={Aspects of Mathematics},
   publisher={Friedr. Vieweg \& Sohn, Braunschweig},
        date={1997},
        ISBN={3-528-28968-6},
         url={https://doi.org/10.1007/978-3-663-10632-6},
        note={Translated from the French and edited by Martin Brown from notes by Michel Waldschmidt, With a foreword by Brown and Serre},
      review={\MR{1757192}},
}

\bib{shiodasymmetry}{article}{
      author={Shioda, Tetsuji},
       title={Constructing curves with high rank via symmetry},
        date={1998},
        ISSN={0002-9327},
     journal={Amer. J. Math.},
      volume={120},
      number={3},
       pages={551\ndash 566},
         url={http://muse.jhu.edu/journals/american_journal_of_mathematics/v120/120.3shioda.pdf},
      review={\MR{1623420}},
}

\bib{stacks-project}{misc}{
      author={{Stacks Project Authors}, The},
       title={\textit{Stacks Project}},
         how={\url{https://stacks.math.columbia.edu}},
        date={2018},
}

\bib{stichtenoth}{article}{
      author={Stichtenoth, Henning},
       title={Die {U}ngleichung von {C}astelnuovo},
        date={1984},
        ISSN={0075-4102},
     journal={J. Reine Angew. Math.},
      volume={348},
       pages={197\ndash 202},
         url={https://doi.org/10.1515/crll.1984.348.197},
      review={\MR{733931}},
}

\bib{umezu1}{article}{
      author={Shioda, Tetsuji},
      author={Umezu, Yumiko},
       title={On {N}\'{e}ron's construction of curves with high rank. {I}},
        date={1999},
        ISSN={0010-258X},
     journal={Comment. Math. Univ. St. Paul.},
      volume={48},
      number={1},
       pages={35\ndash 47},
      review={\MR{1684766}},
}

\bib{tautz-real-mult}{article}{
      author={Tautz, Walter},
      author={Top, Jaap},
      author={Verberkmoes, Alain},
       title={Explicit hyperelliptic curves with real multiplication and permutation polynomials},
        date={1991},
        ISSN={0008-414X,1496-4279},
     journal={Canad. J. Math.},
      volume={43},
      number={5},
       pages={1055\ndash 1064},
         url={https://doi.org/10.4153/CJM-1991-061-x},
      review={\MR{1138583}},
}

\bib{ulmerranks}{article}{
      author={Ulmer, Douglas},
       title={{$L$}-functions with large analytic rank and abelian varieties with large algebraic rank over function fields},
        date={2007},
        ISSN={0020-9910},
     journal={Invent. Math.},
      volume={167},
      number={2},
       pages={379\ndash 408},
         url={https://doi.org/10.1007/s00222-006-0018-x},
      review={\MR{2270458}},
}

\bib{discursus}{misc}{
      author={Watkins, Mark},
       title={A discursus on $21$ as a bound for ranks of elliptic curves over $\mathbb{Q}$, and sundry related topics},
        date={2015},
        note={Available at \url{http://magma.maths.usyd.edu.au/~watkins/papers/DISCURSUS.pdf}},
}

\end{biblist}
\end{bibdiv}
\end{document}